\documentclass[11pt]{article}
\usepackage{amsfonts}
\long\def\ig#1{\relax}
\ig{Thanks to Roberto Minio for this def'n.  Compare the def'n of
\comment in AMSTeX.}

\newcount \coefa
\newcount \coefb
\newcount \coefc
\newcount\tempcounta
\newcount\tempcountb
\newcount\tempcountc
\newcount\tempcountd
\newcount\xext
\newcount\yext
\newcount\xoff
\newcount\yoff
\newcount\gap%
\newcount\arrowtypea
\newcount\arrowtypeb
\newcount\arrowtypec
\newcount\arrowtyped
\newcount\arrowtypee
\newcount\height
\newcount\width
\newcount\xpos
\newcount\ypos
\newcount\run
\newcount\rise
\newcount\arrowlength
\newcount\halflength
\newcount\arrowtype
\newdimen\tempdimen
\newdimen\xlen
\newdimen\ylen
\newsavebox{\tempboxa}%
\newsavebox{\tempboxb}%
\newsavebox{\tempboxc}%

\makeatletter
\def\settypes(#1,#2,#3){\arrowtypea#1 \arrowtypeb#2 \arrowtypec#3}
\def\settoheight#1#2{\setbox\@tempboxa\hbox{#2}#1\ht\@tempboxa\relax}%
\def\settodepth#1#2{\setbox\@tempboxa\hbox{#2}#1\dp\@tempboxa\relax}%
\def\settokens[#1`#2`#3`#4]{%
     \def\tokena{#1}\def\tokenb{#2}\def\tokenc{#3}\def\tokend{#4}}
\def\setsqparms[#1`#2`#3`#4;#5`#6]{%
\arrowtypea #1
\arrowtypeb #2
\arrowtypec #3
\arrowtyped #4
\width #5
\height #6
}
\def\setpos(#1,#2){\xpos=#1 \ypos#2}

\def\bfig{\begin{picture}(\xext,\yext)(\xoff,\yoff)}
\def\efig{\end{picture}}

\def\putbox(#1,#2)#3{\put(#1,#2){\makebox(0,0){$#3$}}}

\def\settriparms[#1`#2`#3;#4]{\settripairparms[#1`#2`#3`1`1;#4]}%

\def\settripairparms[#1`#2`#3`#4`#5;#6]{%
\arrowtypea #1
\arrowtypeb #2
\arrowtypec #3
\arrowtyped #4
\arrowtypee #5
\width #6
\height #6
}

\def\resetparms{\settripairparms[1`1`1`1`1;500]\width 500}

\resetparms

\def\mvector(#1,#2)#3{
\put(0,0){\vector(#1,#2){#3}}%
\put(0,0){\vector(#1,#2){30}}%
}
\def\evector(#1,#2)#3{{
\arrowlength #3
\put(0,0){\vector(#1,#2){\arrowlength}}%
\advance \arrowlength by-30
\put(0,0){\vector(#1,#2){\arrowlength}}%
}}

\def\horsize#1#2{%
\settowidth{\tempdimen}{$#2$}%
#1=\tempdimen
\divide #1 by\unitlength
}

\def\vertsize#1#2{%
\settoheight{\tempdimen}{$#2$}%
#1=\tempdimen
\settodepth{\tempdimen}{$#2$}%
\advance #1 by\tempdimen
\divide #1 by\unitlength
}

\def\vertadjust[#1`#2`#3]{%
\vertsize{\tempcounta}{#1}%
\vertsize{\tempcountb}{#2}%
\ifnum \tempcounta<\tempcountb \tempcounta=\tempcountb \fi
\divide\tempcounta by2
\vertsize{\tempcountb}{#3}%
\ifnum \tempcountb>0 \advance \tempcountb by20 \fi
\ifnum \tempcounta<\tempcountb \tempcounta=\tempcountb \fi
}

\def\horadjust[#1`#2`#3]{%
\horsize{\tempcounta}{#1}%
\horsize{\tempcountb}{#2}%
\ifnum \tempcounta<\tempcountb \tempcounta=\tempcountb \fi
\divide\tempcounta by20
\horsize{\tempcountb}{#3}%
\ifnum \tempcountb>0 \advance \tempcountb by60 \fi
\ifnum \tempcounta<\tempcountb \tempcounta=\tempcountb \fi
}

\ig{ In this procedure, #1 is the paramater that sticks out all the way,
#2 sticks out the least and #3 is a label sticking out half way.  #4 is
the amount of the offset.}

\def\sladjust[#1`#2`#3]#4{%
\tempcountc=#4
\horsize{\tempcounta}{#1}%
\divide \tempcounta by2
\horsize{\tempcountb}{#2}%
\divide \tempcountb by2
\advance \tempcountb by-\tempcountc
\ifnum \tempcounta<\tempcountb \tempcounta=\tempcountb\fi
\divide \tempcountc by2
\horsize{\tempcountb}{#3}%
\advance \tempcountb by-\tempcountc
\ifnum \tempcountb>0 \advance \tempcountb by80\fi
\ifnum \tempcounta<\tempcountb \tempcounta=\tempcountb\fi
\advance\tempcounta by20
}

\def\putvector(#1,#2)(#3,#4)#5#6{{%
\xpos=#1
\ypos=#2
\run=#3
\rise=#4
\arrowlength=#5
\arrowtype=#6
\ifnum \arrowtype<0
    \ifnum \run=0
        \advance \ypos by-\arrowlength
    \else
        \tempcounta \arrowlength
        \multiply \tempcounta by\rise
        \divide \tempcounta by\run
        \ifnum\run>0
            \advance \xpos by\arrowlength
            \advance \ypos by\tempcounta
        \else
            \advance \xpos by-\arrowlength
            \advance \ypos by-\tempcounta
        \fi
    \fi
    \multiply \arrowtype by-1
    \multiply \rise by-1
    \multiply \run by-1
\fi
\ifnum \arrowtype=1
    \put(\xpos,\ypos){\vector(\run,\rise){\arrowlength}}%
\else\ifnum \arrowtype=2
    \put(\xpos,\ypos){\mvector(\run,\rise)\arrowlength}%
\else\ifnum\arrowtype=3
    \put(\xpos,\ypos){\evector(\run,\rise){\arrowlength}}%
\fi\fi\fi
}}

\def\putsplitvector(#1,#2)#3#4{
\xpos #1
\ypos #2
\arrowtype #4
\halflength #3
\arrowlength #3
\gap 140
\advance \halflength by-\gap
\divide \halflength by2
\ifnum \arrowtype=1
    \put(\xpos,\ypos){\line(0,-1){\halflength}}%
    \advance\ypos by-\halflength
    \advance\ypos by-\gap
    \put(\xpos,\ypos){\vector(0,-1){\halflength}}%
\else\ifnum \arrowtype=2
    \put(\xpos,\ypos){\line(0,-1)\halflength}%
    \put(\xpos,\ypos){\vector(0,-1)3}%
    \advance\ypos by-\halflength
    \advance\ypos by-\gap
    \put(\xpos,\ypos){\vector(0,-1){\halflength}}%
\else\ifnum\arrowtype=3
    \put(\xpos,\ypos){\line(0,-1)\halflength}%
    \advance\ypos by-\halflength
    \advance\ypos by-\gap
    \put(\xpos,\ypos){\evector(0,-1){\halflength}}%
\else\ifnum \arrowtype=-1
    \advance \ypos by-\arrowlength
    \put(\xpos,\ypos){\line(0,1){\halflength}}%
    \advance\ypos by\halflength
    \advance\ypos by\gap
    \put(\xpos,\ypos){\vector(0,1){\halflength}}%
\else\ifnum \arrowtype=-2
    \advance \ypos by-\arrowlength
    \put(\xpos,\ypos){\line(0,1)\halflength}%
    \put(\xpos,\ypos){\vector(0,1)3}%
    \advance\ypos by\halflength
    \advance\ypos by\gap
    \put(\xpos,\ypos){\vector(0,1){\halflength}}%
\else\ifnum\arrowtype=-3
    \advance \ypos by-\arrowlength
    \put(\xpos,\ypos){\line(0,1)\halflength}%
    \advance\ypos by\halflength
    \advance\ypos by\gap
    \put(\xpos,\ypos){\evector(0,1){\halflength}}%
\fi\fi\fi\fi\fi\fi
}

\def\putmorphism(#1)(#2,#3)[#4`#5`#6]#7#8#9{{%
\run #2
\rise #3
\ifnum\rise=0
  \puthmorphism(#1)[#4`#5`#6]{#7}{#8}{#9}%
\else\ifnum\run=0
  \putvmorphism(#1)[#4`#5`#6]{#7}{#8}{#9}%
\else
\setpos(#1)%
\arrowlength #7
\arrowtype #8
\ifnum\run=0
\else\ifnum\rise=0
\else
\ifnum\run>0
    \coefa=1
\else
   \coefa=-1
\fi
\ifnum\arrowtype>0
   \coefb=0
   \coefc=-1
\else
   \coefb=\coefa
   \coefc=1
   \arrowtype=-\arrowtype
\fi
\width=2
\multiply \width by\run
\divide \width by\rise
\ifnum \width<0  \width=-\width\fi
\advance\width by60
\if l#9 \width=-\width\fi
\putbox(\xpos,\ypos){#4}
{\multiply \coefa by\arrowlength
\advance\xpos by\coefa
\multiply \coefa by\rise
\divide \coefa by\run
\advance \ypos by\coefa
\putbox(\xpos,\ypos){#5} }%
{\multiply \coefa by\arrowlength
\divide \coefa by2
\advance \xpos by\coefa
\advance \xpos by\width
\multiply \coefa by\rise
\divide \coefa by\run
\advance \ypos by\coefa
\if l#9%
   \put(\xpos,\ypos){\makebox(0,0)[r]{$#6$}}%
\else\if r#9%
   \put(\xpos,\ypos){\makebox(0,0)[l]{$#6$}}%
\fi\fi }%
{\multiply \rise by-\coefc
\multiply \run by-\coefc
\multiply \coefb by\arrowlength
\advance \xpos by\coefb
\multiply \coefb by\rise
\divide \coefb by\run
\advance \ypos by\coefb
\multiply \coefc by70
\advance \ypos by\coefc
\multiply \coefc by\run
\divide \coefc by\rise
\advance \xpos by\coefc
\multiply \coefa by140
\multiply \coefa by\run
\divide \coefa by\rise
\advance \arrowlength by\coefa
\ifnum \arrowtype=1
   \put(\xpos,\ypos){\vector(\run,\rise){\arrowlength}}%
\else\ifnum\arrowtype=2
   \put(\xpos,\ypos){\mvector(\run,\rise){\arrowlength}}%
\else\ifnum\arrowtype=3
   \put(\xpos,\ypos){\evector(\run,\rise){\arrowlength}}%
\fi\fi\fi}\fi\fi\fi\fi}}

\def\puthmorphism(#1,#2)[#3`#4`#5]#6#7#8{{%
\xpos #1
\ypos #2
\width #6
\arrowlength #6
\putbox(\xpos,\ypos){#3\vphantom{#4}}%
{\advance \xpos by\arrowlength
\putbox(\xpos,\ypos){\vphantom{#3}#4}}%
\horsize{\tempcounta}{#3}%
\horsize{\tempcountb}{#4}%
\divide \tempcounta by2
\divide \tempcountb by2
\advance \tempcounta by30
\advance \tempcountb by30
\advance \xpos by\tempcounta
\advance \arrowlength by-\tempcounta
\advance \arrowlength by-\tempcountb
\putvector(\xpos,\ypos)(1,0){\arrowlength}{#7}%
\divide \arrowlength by2
\advance \xpos by\arrowlength
\vertsize{\tempcounta}{#5}%
\divide\tempcounta by2
\advance \tempcounta by20
\if a#8 %
   \advance \ypos by\tempcounta
   \putbox(\xpos,\ypos){#5}%
\else
   \advance \ypos by-\tempcounta
   \putbox(\xpos,\ypos){#5}%
\fi}}

\def\putvmorphism(#1,#2)[#3`#4`#5]#6#7#8{{%
\xpos #1
\ypos #2
\arrowlength #6
\arrowtype #7
\settowidth{\xlen}{$#5$}%
\putbox(\xpos,\ypos){#3}%
{\advance \ypos by-\arrowlength
\putbox(\xpos,\ypos){#4}}%
{\advance\arrowlength by-140
\advance \ypos by-70
\ifdim\xlen>0pt
   \if m#8%
      \putsplitvector(\xpos,\ypos){\arrowlength}{\arrowtype}%
   \else
      \putvector(\xpos,\ypos)(0,-1){\arrowlength}{\arrowtype}%
   \fi
\else
   \putvector(\xpos,\ypos)(0,-1){\arrowlength}{\arrowtype}%
\fi}%
\ifdim\xlen>0pt
   \divide \arrowlength by2
   \advance\ypos by-\arrowlength
   \if l#8%
      \advance \xpos by-40
      \put(\xpos,\ypos){\makebox(0,0)[r]{$#5$}}%
   \else\if r#8%
      \advance \xpos by40
      \put(\xpos,\ypos){\makebox(0,0)[l]{$#5$}}%
   \else
      \putbox(\xpos,\ypos){#5}%
   \fi\fi
\fi
}}

\def\topadjust[#1`#2`#3]{%
\yoff=10
\vertadjust[#1`#2`{#3}]%
\advance \yext by\tempcounta
\advance \yext by 10
}
\def\botadjust[#1`#2`#3]{%
\vertadjust[#1`#2`{#3}]%
\advance \yext by\tempcounta
\advance \yoff by-\tempcounta
}
\def\leftadjust[#1`#2`#3]{%
\xoff=0
\horadjust[#1`#2`{#3}]%
\advance \xext by\tempcounta
\advance \xoff by-\tempcounta
}
\def\rightadjust[#1`#2`#3]{%
\horadjust[#1`#2`{#3}]%
\advance \xext by\tempcounta
}
\def\rightsladjust[#1`#2`#3]{%
\sladjust[#1`#2`{#3}]{\width}%
\advance \xext by\tempcounta
}
\def\leftsladjust[#1`#2`#3]{%
\xoff=0
\sladjust[#1`#2`{#3}]{\width}%
\advance \xext by\tempcounta
\advance \xoff by-\tempcounta
}
\def\adjust[#1`#2;#3`#4;#5`#6;#7`#8]{%
\topadjust[#1``{#2}]
\leftadjust[#3``{#4}]
\rightadjust[#5``{#6}]
\botadjust[#7``{#8}]}

\def\putsquarep<#1>(#2)[#3;#4`#5`#6`#7]{{%
\setsqparms[#1]%
\setpos(#2)%
\settokens[#3]%
\puthmorphism(\xpos,\ypos)[\tokenc`\tokend`{#7}]{\width}{\arrowtyped}b%
\advance\ypos by \height
\puthmorphism(\xpos,\ypos)[\tokena`\tokenb`{#4}]{\width}{\arrowtypea}a%
\putvmorphism(\xpos,\ypos)[``{#5}]{\height}{\arrowtypeb}l%
\advance\xpos by \width
\putvmorphism(\xpos,\ypos)[``{#6}]{\height}{\arrowtypec}r%
}}

\def\putsquare{\@ifnextchar <{\putsquarep}{\putsquarep%
   <\arrowtypea`\arrowtypeb`\arrowtypec`\arrowtyped;\width`\height>}}
\def\square{\@ifnextchar< {\squarep}{\squarep
   <\arrowtypea`\arrowtypeb`\arrowtypec`\arrowtyped;\width`\height>}}
\def\squarep<#1>[#2`#3`#4`#5;#6`#7`#8`#9]{{
\setsqparms[#1]
\xext=\width                                          
\yext=\height                                         
\topadjust[#2`#3`{#6}]
\botadjust[#4`#5`{#9}]
\leftadjust[#2`#4`{#7}]
\rightadjust[#3`#5`{#8}]
\begin{picture}(\xext,\yext)(\xoff,\yoff)
\putsquarep<\arrowtypea`\arrowtypeb`\arrowtypec`\arrowtyped;\width`\height>%
(0,0)[#2`#3`#4`#5;#6`#7`#8`{#9}]%
\end{picture}%
}}

\def\putptrianglep<#1>(#2,#3)[#4`#5`#6;#7`#8`#9]{{%
\settriparms[#1]%
\xpos=#2 \ypos=#3
\advance\ypos by \height
\puthmorphism(\xpos,\ypos)[#4`#5`{#7}]{\height}{\arrowtypea}a%
\putvmorphism(\xpos,\ypos)[`#6`{#8}]{\height}{\arrowtypeb}l%
\advance\xpos by\height
\putmorphism(\xpos,\ypos)(-1,-1)[``{#9}]{\height}{\arrowtypec}r%
}}

\def\putptriangle{\@ifnextchar <{\putptrianglep}{\putptrianglep
   <\arrowtypea`\arrowtypeb`\arrowtypec;\height>}}
\def\ptriangle{\@ifnextchar <{\ptrianglep}{\ptrianglep
   <\arrowtypea`\arrowtypeb`\arrowtypec;\height>}}

\def\ptrianglep<#1>[#2`#3`#4;#5`#6`#7]{{
\settriparms[#1]%
\width=\height                         
\xext=\width                           
\yext=\width                           
\topadjust[#2`#3`{#5}]
\botadjust[#3``]
\leftadjust[#2`#4`{#6}]
\rightsladjust[#3`#4`{#7}]
\begin{picture}(\xext,\yext)(\xoff,\yoff)
\putptrianglep<\arrowtypea`\arrowtypeb`\arrowtypec;\height>%
(0,0)[#2`#3`#4;#5`#6`{#7}]%
\end{picture}%
}}

\def\putqtrianglep<#1>(#2,#3)[#4`#5`#6;#7`#8`#9]{{%
\settriparms[#1]%
\xpos=#2 \ypos=#3
\advance\ypos by\height
\puthmorphism(\xpos,\ypos)[#4`#5`{#7}]{\height}{\arrowtypea}a%
\putmorphism(\xpos,\ypos)(1,-1)[``{#8}]{\height}{\arrowtypeb}l%
\advance\xpos by\height
\putvmorphism(\xpos,\ypos)[`#6`{#9}]{\height}{\arrowtypec}r%
}}

\def\putqtriangle{\@ifnextchar <{\putqtrianglep}{\putqtrianglep
   <\arrowtypea`\arrowtypeb`\arrowtypec;\height>}}
\def\qtriangle{\@ifnextchar <{\qtrianglep}{\qtrianglep
   <\arrowtypea`\arrowtypeb`\arrowtypec;\height>}}

\def\qtrianglep<#1>[#2`#3`#4;#5`#6`#7]{{
\settriparms[#1]
\width=\height                         
\xext=\width                           
\yext=\height                          
\topadjust[#2`#3`{#5}]
\botadjust[#4``]
\leftsladjust[#2`#4`{#6}]
\rightadjust[#3`#4`{#7}]
\begin{picture}(\xext,\yext)(\xoff,\yoff)
\putqtrianglep<\arrowtypea`\arrowtypeb`\arrowtypec;\height>%
(0,0)[#2`#3`#4;#5`#6`{#7}]%
\end{picture}%
}}

\def\putdtrianglep<#1>(#2,#3)[#4`#5`#6;#7`#8`#9]{{%
\settriparms[#1]%
\xpos=#2 \ypos=#3
\puthmorphism(\xpos,\ypos)[#5`#6`{#9}]{\height}{\arrowtypec}b%
\advance\xpos by \height \advance\ypos by\height
\putmorphism(\xpos,\ypos)(-1,-1)[``{#7}]{\height}{\arrowtypea}l%
\putvmorphism(\xpos,\ypos)[#4``{#8}]{\height}{\arrowtypeb}r%
}}

\def\putdtriangle{\@ifnextchar <{\putdtrianglep}{\putdtrianglep
   <\arrowtypea`\arrowtypeb`\arrowtypec;\height>}}
\def\dtriangle{\@ifnextchar <{\dtrianglep}{\dtrianglep
   <\arrowtypea`\arrowtypeb`\arrowtypec;\height>}}

\def\dtrianglep<#1>[#2`#3`#4;#5`#6`#7]{{
\settriparms[#1]
\width=\height                         
\xext=\width                           
\yext=\height                          
\topadjust[#2``]
\botadjust[#3`#4`{#7}]
\leftsladjust[#3`#2`{#5}]
\rightadjust[#2`#4`{#6}]
\begin{picture}(\xext,\yext)(\xoff,\yoff)
\putdtrianglep<\arrowtypea`\arrowtypeb`\arrowtypec;\height>%
(0,0)[#2`#3`#4;#5`#6`{#7}]%
\end{picture}%
}}

\def\putbtrianglep<#1>(#2,#3)[#4`#5`#6;#7`#8`#9]{{%
\settriparms[#1]%
\xpos=#2 \ypos=#3
\puthmorphism(\xpos,\ypos)[#5`#6`{#9}]{\height}{\arrowtypec}b%
\advance\ypos by\height
\putmorphism(\xpos,\ypos)(1,-1)[``{#8}]{\height}{\arrowtypeb}r%
\putvmorphism(\xpos,\ypos)[#4``{#7}]{\height}{\arrowtypea}l%
}}

\def\putbtriangle{\@ifnextchar <{\putbtrianglep}{\putbtrianglep
   <\arrowtypea`\arrowtypeb`\arrowtypec;\height>}}
\def\btriangle{\@ifnextchar <{\btrianglep}{\btrianglep
   <\arrowtypea`\arrowtypeb`\arrowtypec;\height>}}

\def\btrianglep<#1>[#2`#3`#4;#5`#6`#7]{{
\settriparms[#1]
\width=\height                         
\xext=\width                           
\yext=\height                          
\topadjust[#2``]
\botadjust[#3`#4`{#7}]
\leftadjust[#2`#3`{#5}]
\rightsladjust[#4`#2`{#6}]
\begin{picture}(\xext,\yext)(\xoff,\yoff)
\putbtrianglep<\arrowtypea`\arrowtypeb`\arrowtypec;\height>%
(0,0)[#2`#3`#4;#5`#6`{#7}]%
\end{picture}%
}}

\def\putAtrianglep<#1>(#2,#3)[#4`#5`#6;#7`#8`#9]{{%
\settriparms[#1]%
\xpos=#2 \ypos=#3
{\multiply \height by2
\puthmorphism(\xpos,\ypos)[#5`#6`{#9}]{\height}{\arrowtypec}b}%
\advance\xpos by\height \advance\ypos by\height
\putmorphism(\xpos,\ypos)(-1,-1)[#4``{#7}]{\height}{\arrowtypea}l%
\putmorphism(\xpos,\ypos)(1,-1)[``{#8}]{\height}{\arrowtypeb}r%
}}

\def\putAtriangle{\@ifnextchar <{\putAtrianglep}{\putAtrianglep
   <\arrowtypea`\arrowtypeb`\arrowtypec;\height>}}
\def\Atriangle{\@ifnextchar <{\Atrianglep}{\Atrianglep
   <\arrowtypea`\arrowtypeb`\arrowtypec;\height>}}

\def\Atrianglep<#1>[#2`#3`#4;#5`#6`#7]{{
\settriparms[#1]
\width=\height                         
\xext=\width                           
\yext=\height                          
\topadjust[#2``]
\botadjust[#3`#4`{#7}]
\multiply \xext by2 
\leftsladjust[#3`#2`{#5}]
\rightsladjust[#4`#2`{#6}]
\begin{picture}(\xext,\yext)(\xoff,\yoff)%
\putAtrianglep<\arrowtypea`\arrowtypeb`\arrowtypec;\height>%
(0,0)[#2`#3`#4;#5`#6`{#7}]%
\end{picture}%
}}

\def\putAtrianglepairp<#1>(#2)[#3;#4`#5`#6`#7`#8]{{
\settripairparms[#1]%
\setpos(#2)%
\settokens[#3]%
\puthmorphism(\xpos,\ypos)[\tokenb`\tokenc`{#7}]{\height}{\arrowtyped}b%
\advance\xpos by\height
\advance\ypos by\height
\putmorphism(\xpos,\ypos)(-1,-1)[\tokena``{#4}]{\height}{\arrowtypea}l%
\putvmorphism(\xpos,\ypos)[``{#5}]{\height}{\arrowtypeb}m%
\putmorphism(\xpos,\ypos)(1,-1)[``{#6}]{\height}{\arrowtypec}r%
}}

\def\putAtrianglepair{\@ifnextchar <{\putAtrianglepairp}{\putAtrianglepairp%
   <\arrowtypea`\arrowtypeb`\arrowtypec`\arrowtyped`\arrowtypee;\height>}}
\def\Atrianglepair{\@ifnextchar <{\Atrianglepairp}{\Atrianglepairp%
   <\arrowtypea`\arrowtypeb`\arrowtypec`\arrowtyped`\arrowtypee;\height>}}

\def\Atrianglepairp<#1>[#2;#3`#4`#5`#6`#7]{{%
\settripairparms[#1]%
\settokens[#2]%
\width=\height
\xext=\width
\yext=\height
\topadjust[\tokena``]%
\vertadjust[\tokenb`\tokenc`{#6}]
\tempcountd=\tempcounta                       
\vertadjust[\tokenc`\tokend`{#7}]
\ifnum\tempcounta<\tempcountd                 
\tempcounta=\tempcountd\fi                    
\advance \yext by\tempcounta                  
\advance \yoff by-\tempcounta                 %
\multiply \xext by2 
\leftsladjust[\tokenb`\tokena`{#3}]
\rightsladjust[\tokend`\tokena`{#5}]%
\begin{picture}(\xext,\yext)(\xoff,\yoff)%
\putAtrianglepairp
<\arrowtypea`\arrowtypeb`\arrowtypec`\arrowtyped`\arrowtypee;\height>%
(0,0)[#2;#3`#4`#5`#6`{#7}]%
\end{picture}%
}}

\def\putVtrianglep<#1>(#2,#3)[#4`#5`#6;#7`#8`#9]{{%
\settriparms[#1]%
\xpos=#2 \ypos=#3
\advance\ypos by\height
{\multiply\height by2
\puthmorphism(\xpos,\ypos)[#4`#5`{#7}]{\height}{\arrowtypea}a}%
\putmorphism(\xpos,\ypos)(1,-1)[`#6`{#8}]{\height}{\arrowtypeb}l%
\advance\xpos by\height
\advance\xpos by\height
\putmorphism(\xpos,\ypos)(-1,-1)[``{#9}]{\height}{\arrowtypec}r%
}}

\def\putVtriangle{\@ifnextchar <{\putVtrianglep}{\putVtrianglep
   <\arrowtypea`\arrowtypeb`\arrowtypec;\height>}}
\def\Vtriangle{\@ifnextchar <{\Vtrianglep}{\Vtrianglep
   <\arrowtypea`\arrowtypeb`\arrowtypec;\height>}}

\def\Vtrianglep<#1>[#2`#3`#4;#5`#6`#7]{{
\settriparms[#1]
\width=\height                         
\xext=\width                           
\yext=\height                          
\topadjust[#2`#3`{#5}]
\botadjust[#4``]
\multiply \xext by2 
\leftsladjust[#2`#3`{#6}]
\rightsladjust[#3`#4`{#7}]
\begin{picture}(\xext,\yext)(\xoff,\yoff)%
\putVtrianglep<\arrowtypea`\arrowtypeb`\arrowtypec;\height>%
(0,0)[#2`#3`#4;#5`#6`{#7}]%
\end{picture}%
}}

\def\putVtrianglepairp<#1>(#2)[#3;#4`#5`#6`#7`#8]{{
\settripairparms[#1]%
\setpos(#2)%
\settokens[#3]%
\advance\ypos by\height
\putmorphism(\xpos,\ypos)(1,-1)[`\tokend`{#6}]{\height}{\arrowtypec}l%
\puthmorphism(\xpos,\ypos)[\tokena`\tokenb`{#4}]{\height}{\arrowtypea}a%
\advance\xpos by\height
\putvmorphism(\xpos,\ypos)[``{#7}]{\height}{\arrowtyped}m%
\advance\xpos by\height
\putmorphism(\xpos,\ypos)(-1,-1)[``{#8}]{\height}{\arrowtypee}r%
}}

\def\putVtrianglepair{\@ifnextchar <{\putVtrianglepairp}{\putVtrianglepairp%
    <\arrowtypea`\arrowtypeb`\arrowtypec`\arrowtyped`\arrowtypee;\height>}}
\def\Vtrianglepair{\@ifnextchar <{\Vtrianglepairp}{\Vtrianglepairp%
    <\arrowtypea`\arrowtypeb`\arrowtypec`\arrowtyped`\arrowtypee;\height>}}

\def\Vtrianglepairp<#1>[#2;#3`#4`#5`#6`#7]{{%
\settripairparms[#1]%
\settokens[#2]
\xext=\height                  
\width=\height                 
\yext=\height                  
\vertadjust[\tokena`\tokenb`{#4}]
\tempcountd=\tempcounta        
\vertadjust[\tokenb`\tokenc`{#5}]
\ifnum\tempcounta<\tempcountd%
\tempcounta=\tempcountd\fi
\advance \yext by\tempcounta
\botadjust[\tokend``]%
\multiply \xext by2
\leftsladjust[\tokena`\tokend`{#6}]%
\rightsladjust[\tokenc`\tokend`{#7}]%
\begin{picture}(\xext,\yext)(\xoff,\yoff)%
\putVtrianglepairp
<\arrowtypea`\arrowtypeb`\arrowtypec`\arrowtyped`\arrowtypee;\height>%
(0,0)[#2;#3`#4`#5`#6`{#7}]%
\end{picture}%
}}

\def\putCtrianglep<#1>(#2,#3)[#4`#5`#6;#7`#8`#9]{{%
\settriparms[#1]%
\xpos=#2 \ypos=#3
\advance\ypos by\height
\putmorphism(\xpos,\ypos)(1,-1)[``{#9}]{\height}{\arrowtypec}l%
\advance\xpos by\height
\advance\ypos by\height
\putmorphism(\xpos,\ypos)(-1,-1)[#4`#5`{#7}]{\height}{\arrowtypea}l%
{\multiply\height by 2
\putvmorphism(\xpos,\ypos)[`#6`{#8}]{\height}{\arrowtypeb}r}%
}}

\def\putCtriangle{\@ifnextchar <{\putCtrianglep}{\putCtrianglep
    <\arrowtypea`\arrowtypeb`\arrowtypec;\height>}}
\def\Ctriangle{\@ifnextchar <{\Ctrianglep}{\Ctrianglep
    <\arrowtypea`\arrowtypeb`\arrowtypec;\height>}}

\def\Ctrianglep<#1>[#2`#3`#4;#5`#6`#7]{{
\settriparms[#1]
\width=\height                          
\xext=\width                            
\yext=\height                           
\multiply \yext by2 
\topadjust[#2``]
\botadjust[#4``]
\sladjust[#3`#2`{#5}]{\width}
\tempcountd=\tempcounta                 
\sladjust[#3`#4`{#7}]{\width}
\ifnum \tempcounta<\tempcountd          
\tempcounta=\tempcountd\fi              
\advance \xext by\tempcounta            
\advance \xoff by-\tempcounta           %
\rightadjust[#2`#4`{#6}]
\begin{picture}(\xext,\yext)(\xoff,\yoff)%
\putCtrianglep<\arrowtypea`\arrowtypeb`\arrowtypec;\height>%
(0,0)[#2`#3`#4;#5`#6`{#7}]%
\end{picture}%
}}

\def\putDtrianglep<#1>(#2,#3)[#4`#5`#6;#7`#8`#9]{{%
\settriparms[#1]%
\xpos=#2 \ypos=#3
\advance\xpos by\height \advance\ypos by\height
\putmorphism(\xpos,\ypos)(-1,-1)[``{#9}]{\height}{\arrowtypec}r%
\advance\xpos by-\height \advance\ypos by\height
\putmorphism(\xpos,\ypos)(1,-1)[`#5`{#8}]{\height}{\arrowtypeb}r%
{\multiply\height by 2
\putvmorphism(\xpos,\ypos)[#4`#6`{#7}]{\height}{\arrowtypea}l}%
}}

\def\putDtriangle{\@ifnextchar <{\putDtrianglep}{\putDtrianglep
    <\arrowtypea`\arrowtypeb`\arrowtypec;\height>}}
\def\Dtriangle{\@ifnextchar <{\Dtrianglep}{\Dtrianglep
   <\arrowtypea`\arrowtypeb`\arrowtypec;\height>}}

\def\Dtrianglep<#1>[#2`#3`#4;#5`#6`#7]{{
\settriparms[#1]
\width=\height                         
\xext=\height                          
\yext=\height                          
\multiply \yext by2 
\topadjust[#2``]
\botadjust[#4``]
\leftadjust[#2`#4`{#5}]
\sladjust[#3`#2`{#5}]{\height}
\tempcountd=\tempcountd                
\sladjust[#3`#4`{#7}]{\height}
\ifnum \tempcounta<\tempcountd         
\tempcounta=\tempcountd\fi             
\advance \xext by\tempcounta           %
\begin{picture}(\xext,\yext)(\xoff,\yoff)
\putDtrianglep<\arrowtypea`\arrowtypeb`\arrowtypec;\height>%
(0,0)[#2`#3`#4;#5`#6`{#7}]%
\end{picture}%
}}

\def\setrecparms[#1`#2]{\width=#1 \height=#2}%
\def\recursep<#1`#2>[#3;#4`#5`#6`#7`#8]{{%
\width=#1 \height=#2
\settokens[#3]
\settowidth{\tempdimen}{$\tokena$}
\ifdim\tempdimen=0pt
  \savebox{\tempboxa}{\hbox{$\tokenb$}}%
  \savebox{\tempboxb}{\hbox{$\tokend$}}%
  \savebox{\tempboxc}{\hbox{$#6$}}%
\else
  \savebox{\tempboxa}{\hbox{$\hbox{$\tokena$}\times\hbox{$\tokenb$}$}}%
  \savebox{\tempboxb}{\hbox{$\hbox{$\tokena$}\times\hbox{$\tokend$}$}}%
  \savebox{\tempboxc}{\hbox{$\hbox{$\tokena$}\times\hbox{$#6$}$}}%
\fi
\ypos=\height
\divide\ypos by 2
\xpos=\ypos
\advance\xpos by \width
\xext=\xpos \yext=\height
\topadjust[#3`\usebox{\tempboxa}`{#4}]%
\botadjust[#5`\usebox{\tempboxb}`{#8}]%
\sladjust[\tokenc`\tokenb`{#5}]{\ypos}%
\tempcountd=\tempcounta
\sladjust[\tokenc`\tokend`{#5}]{\ypos}%
\ifnum \tempcounta<\tempcountd
\tempcounta=\tempcountd\fi
\advance \xext by\tempcounta
\advance \xoff by-\tempcounta
\rightadjust[\usebox{\tempboxa}`\usebox{\tempboxb}`\usebox{\tempboxc}]%
\bfig
\putCtrianglep<-1`1`1;\ypos>(0,0)[`\tokenc`;#5`#6`{#7}]%
\puthmorphism(\ypos,0)[\tokend`\usebox{\tempboxb}`{#8}]{\width}{-1}b%
\puthmorphism(\ypos,\height)[\tokenb`\usebox{\tempboxa}`{#4}]{\width}{-1}a%
\advance\ypos by \width
\putvmorphism(\ypos,\height)[``\usebox{\tempboxc}]{\height}1r%
\efig
}}

\def\recurse{\@ifnextchar <{\recursep}{\recursep<\width`\height>}}

\def\puttwohmorphisms(#1,#2)[#3`#4;#5`#6]#7#8#9{{%
\puthmorphism(#1,#2)[#3`#4`]{#7}0a
\ypos=#2
\advance\ypos by 20
\puthmorphism(#1,\ypos)[\phantom{#3}`\phantom{#4}`#5]{#7}{#8}a
\advance\ypos by -40
\puthmorphism(#1,\ypos)[\phantom{#3}`\phantom{#4}`#6]{#7}{#9}b
}}

\def\puttwovmorphisms(#1,#2)[#3`#4;#5`#6]#7#8#9{{%
\putvmorphism(#1,#2)[#3`#4`]{#7}0a
\xpos=#1
\advance\xpos by -20
\putvmorphism(\xpos,#2)[\phantom{#3}`\phantom{#4}`#5]{#7}{#8}l
\advance\xpos by 40
\putvmorphism(\xpos,#2)[\phantom{#3}`\phantom{#4}`#6]{#7}{#9}r
}}

\def\puthcoequalizer(#1)[#2`#3`#4;#5`#6`#7]#8#9{{%
\setpos(#1)%
\puttwohmorphisms(\xpos,\ypos)[#2`#3;#5`#6]{#8}11%
\advance\xpos by #8
\puthmorphism(\xpos,\ypos)[\phantom{#3}`#4`#7]{#8}1{#9}
}}

\def\putvcoequalizer(#1)[#2`#3`#4;#5`#6`#7]#8#9{{%
\setpos(#1)%
\puttwovmorphisms(\xpos,\ypos)[#2`#3;#5`#6]{#8}11%
\advance\ypos by -#8
\putvmorphism(\xpos,\ypos)[\phantom{#3}`#4`#7]{#8}1{#9}
}}

\def\putthreehmorphisms(#1)[#2`#3;#4`#5`#6]#7(#8)#9{{%
\setpos(#1) \settypes(#8)
\if a#9 %
     \vertsize{\tempcounta}{#5}%
     \vertsize{\tempcountb}{#6}%
     \ifnum \tempcounta<\tempcountb \tempcounta=\tempcountb \fi
\else
     \vertsize{\tempcounta}{#4}%
     \vertsize{\tempcountb}{#5}%
     \ifnum \tempcounta<\tempcountb \tempcounta=\tempcountb \fi
\fi
\advance \tempcounta by 60
\puthmorphism(\xpos,\ypos)[#2`#3`#5]{#7}{\arrowtypeb}{#9}
\advance\ypos by \tempcounta
\puthmorphism(\xpos,\ypos)[\phantom{#2}`\phantom{#3}`#4]{#7}{\arrowtypea}{#9}
\advance\ypos by -\tempcounta \advance\ypos by -\tempcounta
\puthmorphism(\xpos,\ypos)[\phantom{#2}`\phantom{#3}`#6]{#7}{\arrowtypec}{#9}
}}

\def\putarc(#1,#2)[#3`#4`#5]#6#7#8{{%
\xpos #1
\ypos #2
\width #6
\arrowlength #6
\putbox(\xpos,\ypos){#3\vphantom{#4}}%
{\advance \xpos by\arrowlength
\putbox(\xpos,\ypos){\vphantom{#3}#4}}%
\horsize{\tempcounta}{#3}%
\horsize{\tempcountb}{#4}%
\divide \tempcounta by2
\divide \tempcountb by2
\advance \tempcounta by30
\advance \tempcountb by30
\advance \xpos by\tempcounta
\advance \arrowlength by-\tempcounta
\advance \arrowlength by-\tempcountb
\halflength=\arrowlength \divide\halflength by 2
\divide\arrowlength by 5
\put(\xpos,\ypos){\bezier{\arrowlength}(0,0)(50,50)(\halflength,50)}
\ifnum #7=-1 \put(\xpos,\ypos){\vector(-3,-2)0} \fi
\advance\xpos by \halflength
\put(\xpos,\ypos){\xpos=\halflength \advance\xpos by -50
   \bezier{\arrowlength}(0,50)(\xpos,50)(\halflength,0)}
\ifnum #7=1 {\advance \xpos by
   \halflength \put(\xpos,\ypos){\vector(3,-2)0}} \fi
\advance\ypos by 50
\vertsize{\tempcounta}{#5}%
\divide\tempcounta by2
\advance \tempcounta by20
\if a#8 %
   \advance \ypos by\tempcounta
   \putbox(\xpos,\ypos){#5}%
\else
   \advance \ypos by-\tempcounta
   \putbox(\xpos,\ypos){#5}%
\fi
}}

\makeatother

\usepackage{amsthm}
\usepackage{amsmath}
\usepackage{dsfont}
\usepackage{stmaryrd}
\usepackage{graphics}
\usepackage[linktocpage]{hyperref}

\newtheorem{theorem}{Theorem}[section]
\newtheorem{lemma}[theorem]{Lemma}
\newtheorem{corollary}[theorem]{Corollary}

\newtheorem{proposition}[theorem]{Proposition}

\makeindex \makeglossary
\begin{document}
\sloppy

\newcommand{\nl}{\hspace{2cm}\\ }

\def\nec{\Box}
\def\pos{\Diamond}
\def\diam{{\tiny\Diamond}}

\def\lc{\lceil}
\def\rc{\rceil}
\def\lf{\lfloor}
\def\rf{\rfloor}
\def\lk{\langle}
\def\rk{\rangle}
\def\blk{\dot{\langle\!\!\langle}}
\def\brk{\dot{\rangle\!\!\rangle}}

\newcommand{\pa}{\parallel}
\newcommand{\lra}{\longrightarrow}
\newcommand{\hra}{\hookrightarrow}
\newcommand{\hla}{\hookleftarrow}
\newcommand{\ra}{\rightarrow}
\newcommand{\la}{\leftarrow}
\newcommand{\lla}{\longleftarrow}
\newcommand{\da}{\downarrow}
\newcommand{\ua}{\uparrow}
\newcommand{\dA}{\downarrow\!\!\!^\bullet}
\newcommand{\uA}{\uparrow\!\!\!_\bullet}
\newcommand{\Da}{\Downarrow}
\newcommand{\DA}{\Downarrow\!\!\!^\bullet}
\newcommand{\UA}{\Uparrow\!\!\!_\bullet}
\newcommand{\Ua}{\Uparrow}
\newcommand{\Lra}{\Longrightarrow}
\newcommand{\Ra}{\Rightarrow}
\newcommand{\Lla}{\Longleftarrow}
\newcommand{\La}{\Leftarrow}
\newcommand{\nperp}{\perp\!\!\!\!\!\setminus\;\;}
\newcommand{\pq}{\preceq}

\newcommand{\lms}{\longmapsto}
\newcommand{\ms}{\mapsto}
\newcommand{\subseteqnot}{\subseteq\hskip-4 mm_\not\hskip3 mm}

\def\o{{\omega}}
\def\sM{{\bf sM}}
\def\bA{{\bf A}}
\def\bEM{{\bf EM}}
\def\bM{{\bf M}}
\def\bN{{\bf N}}
\def\bC{{\bf C}}
\def\bI{{\bf I}}
\def\bK{{\bf K}}
\def\bL{{\bf L}}
\def\bT{{\bf T}}
\def\bS{{\bf S}}
\def\bD{{\bf D}}
\def\bB{{\bf B}}
\def\bW{{\bf W}}
\def\bP{{\bf P}}
\def\bX{{\bf X}}
\def\bU{{\bf U}}
\def\bY{{\bf Y}}
\def\ba{{\bf a}}
\def\bb{{\bf b}}
\def\bc{{\bf c}}
\def\bd{{\bf d}}
\def\bh{{\bf h}}
\def\bi{{\bf i}}
\def\bj{{\bf j}}
\def\bk{{\bf k}}
\def\bm{{\bf m}}
\def\bn{{\bf n}}
\def\bp{{\bf p}}
\def\bq{{\bf q}}
\def\be{{\bf e}}
\def\br{{\bf r}}
\def\bi{{\bf i}}
\def\bs{{\bf s}}
\def\bt{{\bf t}}
\def\bu{{\bf u}}
\def\bv{{\bf v}}
\def\bw{{\bf w}}
\def\bz{{\bf z}}

\def\jeden{{\bf 1}}
\def\dwa{{\bf 2}}
\def\trzy{{\bf 3}}
\def\Lam{{\bf \Lambda}}

\def\cBL{{\cal BL}}
\def\cB{{\cal B}}
\def\cA{{\cal A}}
\def\cC{{\cal C}}
\def\cD{{\cal D}}
\def\cE{{\cal E}}
\def\cEM{{\cal EM}}
\def\cF{{\cal F}}
\def\cG{{\cal G}}
\def\cH{{\cal H}}
\def\cI{{\cal I}}
\def\cJ{{\cal J}}
\def\cK{{\cal K}}
\def\cL{{\cal L}}
\def\cN{{\cal N}}
\def\cM{{\cal M}}
\def\cO{{\cal O}}
\def\cP{{\cal P}}
\def\cQ{{\cal Q}}
\def\cR{{\cal R}}
\def\cS{{\cal S}}
\def\cT{{\cal T}}
\def\cU{{\cal U}}
\def\cV{{\cal V}}
\def\cW{{\cal W}}
\def\cX{{\cal X}}
\def\cY{{\cal Y}}

\def\sgn{{\bf sgn}}
\def\el{{\bf el}}
\def\ll{{\bf ll}}
\def\nextF{{\bf nextF}}
\def\prevF{{\bf prevF}}
\def\nextP{{\bf nextP}}
\def\prevP{{\bf prevP}}
{\def\next{{\bf next}}}
\def\neig{{\bf neig}}
\def\neigl{{\bf neig}_{low}}
\def\neigh{{\bf neig}_{high}}
\def\next{{\bf succ}}
\def\prev{{\bf prev}}
\def\top{{\bf top}}
\def\tr{{\bf tr}}
\def\pu{{\bf pu}}
\def\Fl{{\bf Fl}}
\def\Pl{{\bf Pl}}
\def\pFl{{\bf pFl}}
\def\Flags{{\bf Flags}}
\def\Path{{\bf Path}}
\def\diff{{\bf diff}}
\def\id{{\bf id}}
\def\pFlags{{\bf pFlags}}
\def\MpFlags{{\bf MpFlags}}
\def\MFlags{{\bf MFlags}}
\def\Cyl{{\bf Cyl}}
\def\Cyli{{\bf Cyl}_\iota}
\def\Cylh{{\bf Cyl}_h}
\def\Cylhi{{\bf Cyl}_{h,\iota}}
\def\Cylp{{\bf Cyl}_p}
\def\Ope{{\bf Ope}}
\def\pOpe{{\bf pOpe}}
\def\pOpeCard{{\bf pOpeCard}}
\def\Opei{{\bf Ope}_\iota}
\def\pOpei{{\bf pOpe}_\iota}
\def\Opeo{{\bf Ope}_\omega}
\def\pOpeo{{\bf pOpe}_\omega}
\def\pHg{{\bf pHg}}
\def\pHgi{{\bf pHg_{\iota}}}
\def\pOHgi{{\bf pOHg_{\iota}}}
\def\cHi{{\cal H}_{\iota}}
\def\Fp{{\bf Fp}}
\def\sgn{{\bf sgn}}
\def\lhds{{\lhd\!\!\!\cdot\;\,}}
\def\inn{{\in\!\!\!\rightharpoondown}}

\def\sdiff{\stackrel{.}{-}}

\def\LMF{{\bf LMF}}
\def\Mon{{\bf Mon}}
\def\LAdj{{\bf LAdj}}
\def\Adj{{\bf Adj}}
\def\RAdj{{\bf RAdj}}
\def\Act{{\bf Act}}
\def\Clsd{{\bf Closed}}
\def\Mlcv{{\bf Malcev}}
\def\cMlcv{{\bf coMalcev}}
\def\Mnd{{\bf Mnd}}
\def\bCat{{{\bf Cat}}}
\def\Cat{{{\bf Cat}}}
\def\Gpd{{{\bf Gpd}}}
\def\Gph{{{\bf Gph}}}
\def\Fib{{{\bf Fib}}}
\def\DFib{{{\bf DFib}}}
\def\BiFib{{{\bf BiFib}}}
\def\Catrc{{{\bf Cat}_{rc}}}
\def\Monrc{{{\bf Mon}_{rc}}}
\def\Mod{{\bf Mod}}
\def\Modm{{\bf cMod}}
\def\cMod{{\bf cMod}}
\def\mon{{{\bf mon}}}
\def\act{{{\bf act}}}
\def\bem{{{\bf em}}}
\def\bkem{{{\bf kem}}}
\def\bKl{{{\bf Kl}}}

\def\oG{{{\omega}Gr}}
\def\mts{{MltSet}}
\def\pOpeSet{{\bf pOpeSet}}


\def\oCat{{\bf \o Cat}}
\def\oGph{{\bf \o Gph}}
\def\AMon{{\bf AnMnd}}
\def\An{{\bf An}}
\def\Poly{{\bf Poly}}
\def\San{{\bf San}}
\def\Taut{{\bf Taut}}
\def\PMnd{{\bf PolyMnd}}
\def\SanMnd{{\bf SanMnd}}
\def\RiMnd{{\bf RiMnd}}
\def\End{{\bf End}}

\def\Diag{{\bf Diag}}

\def\ET{\bf ET}
\def\RegET{\bf RegET}
\def\RET{\bf RegET}
\def\LrET{\bf LrET}
\def\RiET{\bf RiET}
\def\SregET{\bf SregET}
\def\Cart{\bf Cart}
\def\wCart{\bf wCart}
\def\CartMnd{\bf CartMnd}
\def\wCartMnd{\bf wCartMnd}

\def\LT{\bf LT}
\def\RegLT{\bf RegLT}
\def\ALT{\bf AnLT}
\def\RiLT{\bf RiLT}

\def\FOp{\bf FOp}
\def\RegOp{\bf RegOp}
\def\SOp{\bf SOp}
\def\RiOp{\bf RiOp}

\def\MonCat{{{\bf MonCat}}}
\def\ActMonCat{{{\bf ActMonCat}}}

\def\F{\mathds{F}}
\def\E{\mathds{E}}
\def\S{\mathds{S}}
\def\I{\mathds{I}}
\def\B{\mathds{B}}

\def\V{\mathds{V}}
\def\W{\mathds{W}}
\def\M{\mathds{M}}
\def\N{\mathds{N}}
\def\R{\mathds{R}}

\pagenumbering{arabic} \setcounter{page}{1}

\title{\bf\Large Positive Opetopes with Contractions form a Test Category}

\author{Marek Zawadowski\\
Instytut Matematyki, Uniwersytet Warszawski\\
}

\maketitle

\begin{abstract} We show that the category of positive opetopes with contraction morphisms, i.e. all face maps and some degeneracies, forms a test category.

The category of positive opetopic sets $\pOpeSet$ can be defined as a full subcategory of the category of polygraphs $\Poly$. An object of $\pOpeSet$ has generators whose codomains are again generators and whose domains are non-identity cells (i.e. non-empty composition of generators).  The category $\pOpeSet$ is a presheaf category with the exponent being called the category of positive opetopes $\pOpe$. Objects of $\pOpe$ are called positive opetopes and morphisms are face maps only. Since $\Poly$ has a full-on-isomorphisms embedding into the category of $\omega$-categories $\oCat$, we can think of morphisms in $\pOpe$ as $\omega$-functors that send generators to generators. The category of positive opetopes with contractions $\pOpe_\iota$ has the same objects and face maps $\pOpe$, but in addition it has some degeneracy maps. A morphism in $\pOpe_\iota$ is an $\omega$-functor that sends generators to either generators or to identities on generators. We show that the category $\pOpe_\iota$ is a test category.
 \end{abstract}

\tableofcontents

\section{Introduction}

The notion of a test category was introduced by Grothendieck in \cite{Gr} to axiomatise small categories so that presheaves on those categories model naturally homotopy types of CW-complexes. The category $\Delta$ is the paradigmatic example of a test category but there are other interesting examples of some cubical categories, Joyal's category $\Theta$ c.f. \cite{Berger} and the dendroidal category c.f. \cite{ACM}, to name some.

The category of opetopic sets was introduced in \cite{BD}, see also \cite{HMP}, \cite{Z3}, \cite{SZ}. Combinatorial definitions were given in \cite{Palm} and \cite{Z2}.  Among several equivalent definitions of opetopic sets the one using polygraphs is the simplest to state, c.f. \cite{Bu}. The category of opetpic sets is the full subcategory of polygraphs (considered here as some special $\o$-categories) whose generators have generators as their codomains.
The category of positive opetopic sets can be described as the full subcategory of polygraphs whose generators have generators as their codomains and non-identity cells as domains. The latter category was described combinatorially in \cite{Z1}. Both categories of opetopic sets and positive opetopic sets are presheaf categories with exponent categories called the category of opetopes $\Ope$ and the category of positive opetopes $\pOpe$, respectively. They have a natural non-full but full on isomorphisms embedding into $\oCat$ so that an $\o$-functor between images of opetopes in $\oCat$ sends generators to generators. Contraction morphisms between opetopes can be identified with $\o$-functors between their images in $\oCat$ that send generators to either generators or to (iterated) identities on generators.

The main result of this paper states that the category $\pOpe_\iota$ of positive opetopes with contraction morphisms is a test category.

The paper is organized as follows. In section \ref{sec-peliminaries} we collect all what is needed to prove the main theorem that $\pOpei$ is a test category.  We recall some notions and facts concerning Grothendieck's test categories. Then we sketch the proof of the announced theorem. Next we present the combinatorial definition of categories of positive opetopic cardinals $\pOpeCard$ and positive opetopes $\pOpe$, c.f. \cite{Z1}. Finally, we introduce the contraction maps called  $\iota$-maps, for short, i.e., we extend previously defined category $\pOpe$ that contains face maps only to $\pOpei$ by introducing some degeneracies. The sections \ref{sec-cylinder-in-pOpe} and \ref{sec-cylinder-in-pOpei} constitute the core of the paper. In section \ref{sec-cylinder-in-pOpe} we describe the cylinder functor $\Cylp$ on the category $\widehat{\pOpe}$ and prove some of its properties. In section \ref{sec-cylinder-in-pOpei} we extend the previously defined functor to the category $\widehat{\pOpei}$ and show that it becomes just a functor $I\times (-):\widehat{\pOpei}\lra \widehat{\pOpei}$ of taking product with interval $I$, the only 1-dimensional opetope.

Positive opetopic cardinals and positive opetopes were introduced in \cite{Z1} under the name of positive face structures and principal positive face structures. Among other things we proved there that the category of presheaves on the category of principal positive face structures, here denoted $\widehat{pOpe}$, is equivalent to the category of positive-to-one polygraphs. We have also shown there the category of strict $\o$-categories is monadic over $\widehat{pOpe}$.

All the technical results used in this paper and not proved here come from \cite{Z1}.


\section{Preliminaries}\label{sec-peliminaries}
\subsection{Test categories}\label{sec-test-cats}

Let $\cN : \Cat \lra \widehat{\Delta}$ be the usual nerve functor. We say that a functor $f: C\ra D$ in $\Cat$ is a weak equivalence iff its nerve $\cN(f)$ is a weak homotopy equivalence of simplicial sets.  We denote by $\cW_\infty$ the class of weak equivalences of categories.

For any small category $\cA$, we have an adjunction ($\int_\cA\dashv \cN_A$)
\begin{center}
\xext=1000 \yext=370
\begin{picture}(\xext,\yext)(\xoff,\yoff)
\putmorphism(0,140)(1,0)[\widehat{\cA}`\Cat`]{1000}{0}a
\putmorphism(0,210)(1,0)[\phantom{\widehat{\cA}}`\phantom{\Cat}`\int_\cA]{1000}{1}a
\putmorphism(0,70)(1,0)[\phantom{\widehat{\cA}}`\phantom{\Cat}`\cN_\cA]{1000}{-1}b
\end{picture}
\end{center}
where $\widehat{\cA}$ is the category of presheaves on $\cA$, $\int_\cA(F)$ is the category of elements of presheaf $F$ in $\widehat{\cA}$, and $\cN_\cA(C)(a)=\widehat{\cA}(A_{/a},C)$ for $C$ in $\Cat$ and $a\in \cA$.

\begin{enumerate}
  \item The category $\cA$ is a {\em weak test category} iff for any category $C$, the counit of adjunction $\varepsilon_C : \int_\cA\, \cN_\cA(C) \ra C$ is in $\cW_\infty$.
  \item The category $\cA$ is a {\em local test category} iff for any object $a$ in $\cA$, the slice category $\cA_{/a}$ is a test category.
  \item The category $\cA$ is a {\em test category} iff it is at the same time a weak test category and a local test category.
  \item A presheaf $F$ in $\widehat{\cA}$ is {\em aspherical} if the morphism from the category of elements to the terminal object $\int_{\cA}(F)\ra 1$ is in $\cW_\infty$.
  \item An {\em interval} $(I,s,t)$ in $\widehat{\cA}$ is a presheaf $I$ together with two morphisms $s$ and $t$ from the terminal object that are disjoint, i.e. the square
  \begin{center}
\xext=400 \yext=450
\begin{picture}(\xext,\yext)(\xoff,\yoff)
\setsqparms[1`1`1`1;400`400]
\putsquare(0,0)[0`1`1`I;``s`t]
\end{picture}
\end{center}
is a pullback.
\item Let $X$ be a family of objects of the category $\cC$. An object $F$ of a category $\cC$ is $X$-{\em straight} iff there is a finite linear order $\lk Q,\leq \rk$ and a family of subobjects $\{ a_q \ra F\}_{q\in Q}$ of $F$ in $\cC$ with domains in $X$ such that
\begin{enumerate}
  \item  for any $q\in Q$ which is not the least element in $Q$, $(\bigcup_{q'<q} a_{q'}) \cap a_q$ is (isomorphic to an object) in $X$,
  \item $F = \bigcup_{q\in Q} a_q$,
\end{enumerate}
\item A {\em presheaf $F$ on $\cA$ is straight} iff it is $Ob(\cA)$-straight in $\widehat{\cA}$. We identify the objects of $\cA$ with representable functors in $\widehat{\cA}$.
\end{enumerate}

If $\cA$ has the terminal object, it is aspherical and then it is a test category iff it is a local test category.  Thus since the $\pOpei$ has the terminal object, it is enough to show that it is a local test category to prove that it is a test category.
\vskip 2mm

\begin{proposition}\label{prop-aspherical-presh-characterization} Let $\cA$ be a small category, $F$ a straight presheaf in $\widehat{\cA}$. Then $F$ is aspherical.
\end{proposition}

{\em Proof.} Representable functors are aspherical, c.f. \cite{Mlt}, Remark 1.2.5. Thus, if presheaf  $F$ is straight, it is an iterated sum of pairs of aspherical functors whose intersection is aspherical. By Prop 1.2.7 of \cite{Mlt}, $F$ is aspherical. $\Box$\vskip 2mm

\begin{proposition}\label{prop-local-test-characterization} Let $\cA$ be a small category, $I$ an interval in $\cA$. If for any object $a$ of $\cA$
the functor $I\times a $ is aspherical, then $A$ is a local test category.
\end{proposition}

{\em Proof.} This is an immediate consequence of  Remark 1.2.5 and  Theorem 1.5.6 of \cite{Mlt}. $\Box$\vskip 2mm

\begin{corollary}\label{prop-test-characterization} Let $\cA$ be a small category, $I$ an interval in $\widehat{\cA}$. If $A$ has a terminal object and for any object $a$ of $A$ the presheaf $I\times a$ is aspherical, then $A$ is a test category.  $\Box$
\end{corollary}

One of the reasons to consider test categories is the following theorem, c.f. \cite{Mlt}, \cite{Th}.
\vskip 2mm
\begin{theorem}
If $A$ is a test category, then there exists a model category structure on $\cA$ with monomorphisms as the class of cofibrations and $\cW_\cA=(\int_{\cA})^{-1}(\cW_\infty)$ as the class of weak equivalences. Moreover, the adjunction $\int_\cA\dashv \cN_\cA$ induces an equivalence of categories between homotopy categories
\begin{center}
\xext=1200 \yext=370
\begin{picture}(\xext,\yext)(\xoff,\yoff)
\putmorphism(0,140)(1,0)[\cW^{-1}_{\widehat{\cA}}\widehat{\cA}`\cW_\infty^{-1}\Cat`]{1200}{0}a
\putmorphism(0,210)(1,0)[\phantom{\cW^{-1}_{\widehat{\cA}}\widehat{\cA}}`\phantom{\cW_\infty^{-1}\Cat}`\bar{\int}_\cA]{1200}{1}a
\putmorphism(0,70)(1,0)[\phantom{\cW^{-1}_{\widehat{\cA}}\widehat{\cA}}`\phantom{\cW_\infty^{-1}\Cat}`\bar{\cN}_\cA]{1200}{-1}b
\end{picture}
\end{center}
In particular, the homotopy category $\cW_\cA^{-1}\widehat{\cA}$ is equivalent to Hot, the category of CW-complexes and homotopy classes  of continuous functions.
\end{theorem}
\vskip 2mm

\subsection{Sketch of the proof}\label{sec-sketch-of-proof}

{\em Sketch of the proof.} Since $\pOpe$ has the terminal object, to use Corollary \ref{prop-test-characterization}, it is enough to verify the assumption of Proposition  \ref{prop-local-test-characterization}, where $I$  is the 1-dimensional opetope. To this end we shall show that, for any opetope $P$, the product $I\times P$ is straight in $\widehat{\pOpei}$ and use Proposition \ref{prop-aspherical-presh-characterization}.

For an arbitrary positive opetope $P$, we shall describe the cylinder hypergraph $\Cylh(P)$ over $P$ that is $\pOpe$-straight object in the category of positive hypergraphs  $\pHg$ but it is not a product $I\times P$ in $\pHg$, in general. Its embedding $\Cylp(P)$ into $\widehat{\pOpe}$ is again straight.  Then we will take a left Kan extension $\Cyli(P)$ of $\Cylp(P)$ along inclusion $\kappa:\pOpe^{op} \ra \pOpei^{op}$. $\Cyli(P)$ is also straight. Finally, we will show that  $\Cyli(P)$ is the product of $I$ and $P$ in $\widehat{\pOpei}$.

$\Cyli(P)$ is built from representable functors that  are identified by the so called {\em flags}, a concept introduced by T. Palm in \cite{Palm} under the name of a {\em chamber}. There is a linear order of flags that will give the order required in  Proposition \ref{prop-aspherical-presh-characterization}. The intersections of a flag with the sum of all flags coming before  will be identified as representable functors corresponding to punctured fags, {\em p-flags} for short.

The proof has two parts,  the first is contained in section \ref{sec-cylinder-in-pOpe} and the second in section \ref{sec-cylinder-in-pOpei}. In section \ref{sec-cylinder-in-pOpe} we describe combinatorial structure of cylinder $\Cylh(P)$ as a positive hypergraph. To this end, we define an order on flags, and then a star operation that lifts the faces from the positive opetope $P$ to their versions in the subopetope $P^{\vec{x}}$ of hypergraph $\Cylh(P)$ corresponding to the flag $\vec{x}$. This operation is used to define $\gamma$ (`codomain') and $\delta$ (`domain') operations on faces of $\Cylh(P)$. We also verify that $P^{\vec{x}}$'s are opetopes and that $\Cylh(P)$ is $\pOpe$-straight in $\pHg$ and its image $\Cylp(P)$ is straight in $\widehat{\pOpe}$.

In section \ref{sec-cylinder-in-pOpei} we show that $\Cyli(P)$ is a product of $I\times P$ in $\widehat{\pOpei}$. To this end we need to construct for a given pair of $\iota$-maps, $\rho :Q\ra I$ and $h: Q\ra P$ in $\pHgi$, a suitable $\iota$-map $H:Q\ra  \Cylhi(P)$. Such a map necessarily factorises through an opetope $P^{\vec{x}}\hra \Cylhi(P)$. The construction of the flag $\vec{x}$ for a given pair of $\iota$-maps $\rho$, $h$ is the main difficulty here.  Then we verify that $H$ has the required properties.

\subsection{Positive hypergraphs}

A {\em positive hypergraph} $S$ is a family $\{ S_k\}_{k\in \o}$ of finite sets of faces, a family of functions $\{ \gamma_k : S_{k+1}\ra S_k \}_{k \in\o }$, and a family of total relations $\{ \delta_k : S_{k+1}\ra S_k\}_{k \in \o}$. Moreover, $\delta_0 : S_1\ra S_0$ is a function and only finitely many among sets $\{ S_k\}_{k\in \o}$ are non-empty. As it is always clear from the context, we shall never use the indices of the functions $\gamma$ and $\delta$. We shall ignore the difference between $\gamma(x)$ and $\{ \gamma(x)\}$ and in consequence we shall consider iterated applications of $\gamma$'s and $\delta$'s as sets of faces, e.g. $\delta\delta(x)=\bigcup_{y\in \delta(x)} \delta(y)$ and  $\gamma\delta(x)=\{\gamma(y)\, |\, y\in \delta(x)\}$.

A {\em morphism of positive hypergraphs} $f:S\lra T$ is a family of functions $f_k : S_k \lra T_k$, for $k \in\o$, such that, for $k>0$ and $a\in S_k$, we have  $\gamma (f(a))=f(\gamma(a))$ and $f_{k-1}$ restricts to a bijection
\[ f_a: \delta(a) \lra \delta(f(a)). \]
The category of positive hypergraphs is denoted by $\bf pHg$.

We define a binary relation of {\em lower order} on $<^{S_k,-}$ for $k>0$ as the transitive closure of  the relation $\lhd^{S_k,-}$ on $S_k$ such that, for $a,b\in S_k$,  $a \lhd^{S_k,-} b$ iff $\gamma (a)\in \delta(b)$.  We write $a\perp^- b$ iff either $a <^- b$ or $b <^- a$, and we write $a \leq^- b$ iff either $a=b$ or $a <^- b$.

We also define a  binary relation of {\em upper order} on $<^{S_k,+}$ for $k\geq0$ as the transitive closure of  the relation $\lhd^{S_k,+}$ on $S_k$ such that, for $a,b\in S_k$,  $a \lhd^{S_k,+} b$ iff there is $\alpha\in S_{k+1}$ so that $a\in \delta(\alpha)$ and $\gamma (\alpha)=b$.  We write $a\perp^+ b$ iff either $a <^+ b$ or $b <^+ a$, and we write $a \leq^+ b$ iff either $a=b$ or $a <^+ b$.

\subsection{Positive opetopic cardinals}

A positive hypergraph $S$ is a {\em positive opetopic cardinal} if it is non-empty,
i.e. $S_0\neq\emptyset$ and it satisfies the following four conditions
\begin{enumerate}
 \item {\em Globularity:}  for  $a\in S_{\geq 2}$:
 \[ \gamma\gamma(a) =\gamma\delta(a)-\delta\delta(a),\hskip 15mm \delta\gamma(a) =\delta\delta(a)-\gamma\delta(a);\]
  \item {\em Strictness:} for $k\in\o$, the relation   $<^{S_k,+}$ is a strict order; $<^{S_0,+}$ is linear;
  \item {\em Disjointness:} for $k>0$,
\[\perp^{S_k,-}\cap \perp^{S_k,+}=\emptyset\]
  \item {\em Pencil linearity:} for any $k>0$   and $x\in S_{k-1}$, the sets
  \[ \{ a\in S_k \; | \; x=\gamma(a) \} \;\;\;\;{\rm   and}\;\;\;\; \{ a\in S_k \; | \; x\in \delta(a) \} \]
  are linearly ordered by $<^{S_k,+}$.
\end{enumerate}

The category of positive opetopic cardinals is the full subcategory of $\bf pHg$ whose objects are positive opetopic cardinals. It is denoted by $\pOpeCard$.

\subsection{Positive opetopes}

The {\em size of positive opetopic cardinal} $S$ is a sequence of natural numbers $size(S)=\{ | S_n - \delta (S_{n+1})| \}_{n\in\o}$, with
all elements above $dim(S)$ being equal $0$. We have an order $<$ on such sequences of natural numbers so that $\{ x_n \}_{n\in\o} < \{ y_n \}_{n\in\o}$ iff there is $k\in\o$ such that $x_k< y_k$ and, for all $l>k$, $x_l = y_l$. This order is well founded and hence facts about positive opetopic cardinals can be proven by induction on their size.

Let $P$ be an positive opetopic cardinal.  We say that $P$ is a {\em positive opetope} iff $size(P)_l\leq 1$, for $l\in \o$. By $\pOpe$ we denote full
subcategory of $\bf pHg$ whose objects are positive opetopes.

\vskip 2mm
{\em Notation.}
\begin{enumerate}
\item Let $S$ be a positive hypergraph $S$, $x$ a face of $S$. Then $S[x]$ denotes the least sub-hypergraph of $S$ containing face $x$. The dimension of $S$ is maximal $k$ such that $S_k\neq\emptyset$. We denote by $dim(S)$ the dimension of $S$.
  \item Let $P$ be an opetope.  If  $dim(P)=n$, then the unique face in $P_n$ is denoted by $\bm_P$.
  \item The function $\gamma^{(k)}: P\ra P_{\leq k}$ is an iterated version of the codomain function $\gamma$ defined as follows. For  any $k,l\in \o$ and $p\in P_l$,
\[ \gamma^{(k)}(p) \;= \; \left\{ \begin{array}{ll}
		\gamma\gamma^{(k+1)}(p)   & \mbox{ if } l>k \\
		p  & \mbox{ if } l\leq k.
                                    \end{array}
			    \right. \]
  \item
 For $p\in P_{>0}$, we define a {\em boundary} of $p$ as
\[ \partial(p)=\gamma(p)\cup \delta(p). \]
  \item If $k\leq m$, $p'\in P_k$ and $p\in P_m$, we write $p'\inn p$ iff either $p'=p$ or $p'\in \partial(p)$ or $p'\in \partial\partial\gamma^{(k+2)}(p)$ and we say that  $p'$ {\em occurs in} $p$.
\end{enumerate}

\subsection{The embedding of $\pOpeCard$ into $\oCat$}

There is an embedding
\[ (-)^*:\pOpeCard \lra \oCat \]
defined as follows, c.f. \cite{Z1}. Let $T$ be an opetopic cardinal. The $\o$-category $T^*$ has as $n$-cells pairs $(S,n)$ where $S$ is a subopetopic cardinal of $T$, $dim(S)\leq n$, and $n\geq 0$. A cell $(S,n)$ is called {\em proper} iff  $n=dim(S)$.

For $k<n$, the domain and codomain operations
\[ \bd^{(k)}, \bc^{(k)}: T^*_{n} \lra T^*_k \]
are given, for  $(S,n)\in T^*_{n}$, by
\[  (\bd^{(k)}(S,n))=(\bd^{(k)}(S),k),\hskip 15mm  (\bc^{(k)}(S,n))=(\bc^{(k)}(S),k) \]
where
\[ (\bd^{(k)}(S)_l \;= \; \left\{ \begin{array}{ll}
		\emptyset   & \mbox{ if } l>k \\
		S_k-\gamma(S_{k+1})  & \mbox{ if } l=k \\
        S_l  & \mbox{ if } l<k
                                    \end{array}
			    \right. \]
and
\[ (\bc^{(k)}(S)_l \;= \; \left\{ \begin{array}{ll}
		\emptyset   & \mbox{ if } l>k \\
		S_k-\delta(S_{k+1})  & \mbox{ if } l=k \\
        S_{k-1}-\iota(S_{k+1})  & \mbox{ if } l=k-1\geq 0 \\
        S_l  & \mbox{ if } l<k-1
                                    \end{array}
			    \right. \]
The identity operation
\[ id : T^*_n \lra T^*_{n+1} \]
is given by
\[  (S,n)\mapsto (S,n+1). \]
The composition operation is defined, for pairs of cells $(S,n),(S',n')\in T^*$ with $k\leq n,n'$ such that $\bd^{(k)}(S,n)=\bc^{(k)}(S',n')$, as the sum
\[ (S,n)\circ (S',n')= (S\cup S', \max(n,n')).\]

Now $T^*$ together with operations of domain, codomain, identity and composition is an $\o$-category. If $f:S\ra T$ is a map of opetopic cardinals and $S'$ is a sub-opetopic cardinal of $S$ then the image $f(S')$ is a sub-opetopic cardinal of $T$. This defined the functor $(-)^*$ on morphisms. We recall from \cite{Z1}

\begin{theorem}\label{thm-embedding-pOpeCard}
The embedding
\[ (-)^*:\pOpeCard \lra \oCat \]
is well defined and full on isomorphisms and it factorises through $\Poly\lra \oCat$ via full and faithful functor, $(-)^*:\pOpeCard \lra \Poly$, into the category of polygraphs.
 \end{theorem}

Since we have a (full) embedding $\pOpe\ra \pHg$, we have a functor
\[ \cH: {\pHg} \lra \widehat{\pOpe}\]
such that
\[ S\; \mapsto\; {[ P \mapsto \pHg(P,S) ]}. \]

\vskip 2mm
{\em Remark.} $\pHg$ does not have all colimits, e.g. two inclusions $f$ and $g$ of 1-dim opetope $\{-\stackrel{a}{\ra} + \}$ into a 2-dim poly with two faces in the domain     as follows
    \begin{center} \xext=800 \yext=520
\begin{picture}(\xext,\yext)(\xoff,\yoff)
 \settriparms[-1`1`1;300]
  \putAtriangle(0,0)[\bullet`\bullet`\bullet;f(a)`g(a)`]
\end{picture}
\end{center}
do not have a coequalizer in $\pHg$. However, sums of subobjects and pushouts along monos exist in $\pHg$ and are preserved by embedding $\cH$.

\vskip 2mm

Since positive opetopes are positive hypergraphs, we can consider $\pOpe$-straight objects of $\pHg$. We have

\vskip 2mm
\begin{corollary}
The functor $\cH$ sends  $ob(\pOpe)$-straight hypergraphs in $\pHg$ to straight presheaves in $\widehat{\pOpe}$.
\end{corollary}

\subsection{$\iota$-maps of positive opetopes}\label{sec-iota-maps}

The embedding $(-)^\ast: \pOpe \ra \oCat$ is not full. The morphisms $P^\ast\ra Q^\ast$ in $\oCat$ between images of opetopes are $\o$-functors that send generators to generators.  The category $\pOpei$ with the same objects as $\pOpe$ will be so defined that the embedding $(-)^\ast: \pOpei \ra \oCat$ (denoted the same way) will be full on $\o$-functors that send generators to either generators or identities on generators of a smaller dimension.

Let $P$ and $Q$ be positive opetopes. A {\em contraction morphism of opetopes} (or {\em $\iota$-map}, for short),  $h : Q \lra P$, is a function $h : |Q| \lra |P|$ between faces of opetopes such that
\begin{enumerate}
  \item  $dim(q)\geq dim(h(q))$, for $q\in  Q$;
  \item (preservation of codomains) $h(\gamma^{(k)} (q))= \gamma^{(k)}(h(q))$, for $k\geq 0$ and $q\in Q_{k+1}$;
  \item (preservation of domains)
  \begin{enumerate}
    \item  if $dim(h(q))= dim(q)$, then $h$ restricts to a bijection
  \begin{center}
\xext=1200 \yext=10 \adjust[`I;I`;I`;`I]
\begin{picture}(\xext,\yext)(\xoff,\yoff)
\putmorphism(0,0)(1,0)[(\delta(q)-\ker(h))`\delta(h(q))`h]{800}{1}a
\end{picture}
\end{center}
for $k\geq 0$ and $q\in Q_{k+1}$, where {\em the kernel of h} is defined as $$\ker(h)=\{ q\in Q| dim(q)> dim(h(q))\};$$
    \item if $dim(h(q))= dim(q)-1$, then $h$ restricts to a bijection
  \begin{center}
\xext=1200 \yext=10 \adjust[`I;I`;I`;`I]
\begin{picture}(\xext,\yext)(\xoff,\yoff)
\putmorphism(0,0)(1,0)[(\delta(q)-\ker(h))`\{ h(q)\}`h]{800}{1}a
\end{picture}
\end{center}
for $k\geq 0$ and $q\in Q_{k+1}$;
    \item if $dim(h(q))< dim(q)-1$, then $\delta^{(k)}(q)\subseteq \ker(h)$.
  \end{enumerate}
\end{enumerate}
We call a face $q\in Q_m$ $k$-{\em collapsing} iff $h(q)\in P_{m-k}$. A $0$-collapsing face is called {\em non-collapsing}.

We have an embedding $\kappa:\pOpe\lra \pOpei$ that induces the usual adjunction $\kappa_!\dashv \kappa^*$

\begin{center}
\xext=1000 \yext=370
\begin{picture}(\xext,\yext)(\xoff,\yoff)
\putmorphism(0,140)(1,0)[\widehat{\pOpe}`\widehat{\pOpei}`]{1000}{0}a
\putmorphism(0,210)(1,0)[\phantom{\widehat{\pOpe}}`\phantom{\widehat{\pOpei}}`\kappa_!]{1000}{1}a
\putmorphism(0,70)(1,0)[\phantom{\widehat{\pOpe}}`\phantom{\widehat{\pOpei}}`\kappa^*]{1000}{-1}b
\end{picture}
\end{center}

\begin{lemma}Let $h: Q\ra P$ be a $\iota$-map, $q_1,q_2\in Q-\ker(h)$ and $l<k\in \o$ such that
\[  \begin{array}{rcl}
		\gamma^{(k+1)}(q_1)  & <^- & \gamma^{(k+1)}(q_2) \\
		\gamma^{(k)}(q_1)   & <^+ & \gamma^{(k)}(q_2)\\
        \ldots &\ldots & \ldots \\
        \gamma^{(l+1)}(q_1)   & <^+ & \gamma^{(l+1)}(q_2)\\
        \gamma^{(l)}(q_1)   & = & \gamma^{(l)}(q_2).\\
                                    \end{array}
			    \]
Then there is $l\leq l'< k$ such that
\[  \begin{array}{rcl}
		h(\gamma^{(k+1)}(q_1))  & <^- & h(\gamma^{(k+1)}(q_2)) \\
		h(\gamma^{(k)}(q_1))   & <^+ & h(\gamma^{(k)}(q_2))\\
        \ldots &\ldots & \ldots \\
        h(\gamma^{(l'+1)}(q_1))   & <^+ & h(\gamma^{(l'+1)}(q_2)) \\
        h(\gamma^{(l')}(q_1))   & = & h(\gamma^{(l')}(q_2))\\
                                    \end{array}
			    \]
\end{lemma}

{\em Proof.} Simple check. $\Box$\vskip 2mm

From the above we get immediately

\begin{corollary}\label{coro-iota-preservation}
Let $h: Q\ra P$ be a $\iota$-map, $q_1,q_2\in Q-\ker(h)$. Then
\begin{enumerate}
  \item $q_1<^- q_2$ iff $h(q_1)<^- h(q_2)$;
  \item if $q_1<^+ q_2$, then $h(q_1)\leq^+ h(q_2)$;
  \item if $h(q_1)<^+ h(q_2)$, then $q_1<^+ q_2$;
  \item if $h(q_1)= h(q_2)$, then $q_1\perp^+ q_2$. $\Box$
\end{enumerate}
\end{corollary}

A set $X$ of $k$-faces in a positive opetope $P$ is a {\em $<^+$-interval} (or {\em interval}, for short)  if it is either empty or there are two k-faces $x_0,x_1\in P_k$ such that $x_0\leq^+ x_1$ and $X=\{ x\in P_k | x_0\leq^+ x \leq^+ x_1 \}$. Any interval in any positive opetope is linearly ordered by $\leq^+$.

\begin{corollary}
Let $h:Q\ra P$ -be a contraction of positive opetopes, $p\in P_k$. Then the fiber of $k$-faces (of non-degenerating faces) $h^{-1}(p)- \ker(h)$ of $h$ over $p$ is an interval.
\end{corollary}

{\em Proof.} From Corollary \ref{coro-iota-preservation}.4 we get that $h^{-1}(p)- \ker(h)$ is linearly ordered. And from Corollary \ref{coro-iota-preservation}.2 that this linear order is an interval. $\Box$\vskip 2mm

\begin{proposition}
$P$ - positive opetope. The onto contraction maps $P\ra I$ are in bijective correspondence with 1-faces in $P_1-\gamma(P_2)$.
\end{proposition}
{\em Proof.} Let $I=\{ -\stackrel{a}{\ra}+\}$ be the 1-dimensional opetope.  For  $p\in P_1-\gamma(P_2)$, we define the contraction map $h_p: P\ra I=\{ -\stackrel{a}{\ra}+\}$ as follows
   \[ h_p(x) \;= \; \left\{ \begin{array}{ll}
		-   & \mbox{ if } \gamma^{(0)}(x)\leq^+\delta(p) \\
		+   & \mbox{ if } \gamma(p)\leq^+\delta\gamma^{(1)}(x) \\
        a   & \mbox{ otherwise, i.e. if } p\leq^+ \gamma^{(1)}(x).\\
                                    \end{array}
			    \right. \]
To see that $h_p$ preserves the domain of a face $q\in Q_2$ such that $p\leq^+ \gamma(q)$, it is enough to notice that in this case there is a unique $q'\in \delta(q)$ such that $h_p(q')=a$. On the other hand, if $h: P\ra I$ is an onto $\iota$-map, then there is exactly one face in $p\in P_1-\gamma(P_2)$ such that $h(p)=a$. Then $h=h_p$. The remaining details are left for the reader. $\Box$\vskip 2mm

\subsection{The embedding of $\pOpei$ into $\oCat$}

We extend the embedding functor $(-)^\ast$ to contractions
\[ (-)^* : \pOpei \lra \oCat. \]
Let $h:Q \ra P$ be a contraction morphism in $ \pOpei$. Then
\[ h^\ast : Q^\ast\ra P^\ast \]
is an $\o$-functor such that
\[ h^\ast(k,A) =  (k,\vec{h}(A)) \]
where $(k,A)\in Q^*_k$, and $\vec{h}(A)$ is the set-theoretic image of the opetopic cardinal $A$ under $h$.
\vskip 2mm

\begin{lemma}\label{lemma-iota-emdedding}
Let $h : Q\ra P$ be a contraction morphism. We have
\begin{enumerate}
  \item $h^\ast(Q[q])=P[h(q)]$, for $q\in Q$;
  \item  if $A$ is an opetopic cardinal contained in $Q$ with $dim(A)\leq l$, $k<l$, and $h^*(A)$ is an opetopic cardinal, then
  \[ h^\ast(\bd^{(k)}A) =\bd^{(k)}h^\ast(A), \;\; h^\ast(\bc^{(k)}A) =\bc^{(k)}h^\ast(A). \]
\end{enumerate}
\end{lemma}
{\em Proof.} Ad 1. Fix $q\in Q_{m+1}$, $m\geq 0$. From the definition of $\iota$-maps we immediately have \mbox{$h^\ast(Q[q])\supseteq P[h(q)]$}. To show the opposite inclusion, we suppose to the contrary that there is a face $q'\in Q[q]_k$ such that $h(q')\not\in P[h(q)]$. We can assume that $q'$ is a face of the least dimension with this property. From condition 2. of the definition of $\iota$-maps it follows that $dim(q')=dim(h(q'))$. Otherwise, we could chose $\gamma(q')$ instead of $q'$. We have two cases
\begin{enumerate}
  \item $k=m$ and $q'\in \partial(q)$;
  \item $k<m$ and $q'\in \gamma^{(k)}(q)\cup \iota\gamma^{(k+2)}(q)\cup \delta\gamma^{(k+1)}(q)$.
\end{enumerate}

In the first case, if $q'=\gamma(q)$, we have $h(q')=h(\gamma(q))=\gamma (h(q)) \in P[h(q)]$.
If $q'\in \delta(q)$, we have
$$ h(q')\in h(\delta(q))\cap P_k = h(\delta(q)-\ker(h)) \subseteq  \delta((h(q)) \subseteq P[h(q){]}.$$

In the second case, since $P[h(\gamma^{(k+1)}(q))]\subseteq P[q]$, we can assume that $q$ is in fact $\gamma^{(k+2)}(q)$. This will simplify the notation a lot.  If $q' \in \gamma\gamma(q)\cup \delta\gamma(q)$, then, by the above argument, we have
\[  h(q')\in P[h(\gamma(q))]\subseteq P[h(q)]. \]
Finally, if $q'\in \iota(q)$, then there are $r,s\in \delta(q)$ such that $\gamma(r)=q'\in\delta(s)$.

If $r\not\in \ker(h)$, then $h(r)\in \delta(h(q))\subseteq P[h(q)]$. Thus
$h(q')=\gamma(h(r))\in P[h(q)]$.

If $s\not\in \ker(h)$, then again $h(s)\in \delta(h(q))\subseteq P[h(q)]$.
Thus $h(q')\in \delta(h(s))\subseteq P[h(q)]$.

If $r,s\in \ker(h)$, let $s_0=r, s_1=s,\ldots, s_n$ be the maximal lower $\delta(q)$-path beginning with $r$  such that $s_i\in \ker(h)$, $i=0,\ldots, n$.
If $\gamma(s_n)=\gamma^{(k)}(q)$, then
$$h(q')=h(\gamma^{(k)}(q))\in P[h(q)].$$
Otherwise, there is $t\in \delta(q)-\ker(h)$ such that $\gamma(s_n)\in \delta(t)$. Then $h(t)\in \delta(h(q))$ and $$ h(q')=h(\gamma(s_n))\in \delta(h(t))\cap P_k \subseteq P[q]$$.

Thus in any case $q'\in P[q]$, as required.

Ad 2. It is enough to verify preservation of domains and codomains for proper faces of $Q^*$. So fix a proper face $A\in Q^*_{k+1}$. Since we don't know yet that $h^*(A)$ is an opetopic cardinal, we assume it is.

To show that the domains are preserved, we need to verify that opetopic cardinals $h^*(\bd A)$ and $\bd (h^*(A))$ that are subsets of $h^*(A)$ are equal. The equality in dimensions $\neq k$ is obvious. Thus we need to  show that \[\bd (h^*(A))_k= (h(A_k)\cap P_k)-\gamma(h (A_{k+1})\cap P_{k+1}) \]
is equal to
 \[    h^*(\bd A)_k = h(A_k-\gamma(A_{k+1}))\cap P_k. \]
We have
\[  (h(A_k)\cap P_k)-\gamma(h (A_{k+1})\cap P_{k+1})\subseteq h(A_k-\gamma(A_{k+1}))\cap P_k\subseteq h(A_k)\cap P_k \]
and hence it is enough to show that the sets
\begin{equation} \label{equation-d-preserv}
 \gamma(h (A_{k+1})\cap P_{k+1}),\hskip 1cm h(A_k-\gamma(A_{k+1}))\cap P_k
\end{equation}
are disjoint.

So suppose to the contrary that there is a face $x\in P_k$  that belongs to both sets, i.e. there is $z\in A_{k+1}$ such that $h(z)\in P_{k+1}$ and $\gamma(h(z))=x$ and, moreover, there is $s\in A_k-\gamma(A_{k+1})$ such that $h(s)=x$. Thus
\[ h(s)= x =\gamma(h(z))= h(\gamma(z))   \]
and hence $s\perp^+\gamma(z)$. Since $s\not\in \gamma(A_{k+1})$, we have $s<^+\gamma(z)$. Let
\[ s,a_1,\ldots, a_l,\gamma(z) \]
 be an upper path in $Q$ witnessing it. Since $h(s)=h(\gamma(z))$, it follows that $a_i\in \ker(h)$, for $i=1,\ldots, l$. By pencil linearity, we have $a_l \perp^+ z$. Since $a_i$ are in the kernel of $h$ and $z$ is not, we have $a_l<^+z$ and, in fact, $a_i<^+z$, for $i=1,\ldots, l$. Hence, as $h(z)\not\in \ker(h)$, $\delta(z)<^+s$ and $s\in \gamma(A_{k+1})$ contrary to the supposition. This shows the preservation of domains by $h^*$.

To show that the codomains are preserved, we need to verify that opetopic cardinals $h^*(\bc A)$ and $\bc (h^*(A))$ that are subsets of $h^*(A)$ are equal. The equality in dimensions $\neq k,k-1$ is again obvious. Thus we need to  show that
\[\bc (h^*(A))_k= (h(A_k)\cap P_k)-\delta(h (A_{k+1})\cap P_{k+1}) \]
is equal to
 \[    h^*(\bc A)_k = h(A_k-\delta(A_{k+1}))\cap P_k \]
and that
\[\bc (h^*(A))_{k-1}= (h(A_{k-1})\cap P_{k-1})-\iota(h (A_{k+1})\cap P_{k+1}) \]
is equal to
 \[    h^*(\bc A)_{k-1} = h(A_{k-1}-\iota(A_{k+1}))\cap P_{k-1}. \]

To show these equalities, it is enough to show that the sets
\begin{equation} \label{equation-c-preserv1}
\delta(h (A_{k+1})\cap P_{k+1}),\;\;\;{\rm and}\;\;\; h(A_k-\delta(A_{k+1}))\cap P_k
\end{equation}
and the sets
\begin{equation} \label{equation-c-preserv2}
 \iota(h (A_{k+1})\cap P_{k+1}) ,\;\;\;{\rm and}\;\;\;h(A_{k-1}-\iota(A_{k+1}))\cap P_{k-1}
\end{equation}
are disjoint.

To see that the sets in \ref{equation-c-preserv1} are disjoint, suppose to the contrary that there are $u\in A_k-\delta(A_{k+1})$ and $z\in A_{k+1}$, with $h(z)\in P_{k+1}$ so that  $h(u)=t\in \delta(h(z))$. Since $h$ is a $\iota$-map, we have $t'\in \delta(z)$ such that $h(t')=t$. Thus $h(u)=t=h(t')$ and, by Corollary \ref{coro-iota-preservation}.4, we have $u\perp^+t'$. Since $u\not\in \delta(A_{k+1})$, we have $t'<^+u$ and even $\gamma(z)\leq^+ u$. Since $h(t')\in \delta(h(z))$ and $h(z)\in P_{k+1}$, we have
$$t=h(t')<^+\gamma(h(z))= h(\gamma(z))\leq^+ h(u)=t.$$
As $<^+$ is a strict order, it is a contradiction. This means that the first pair of sets displayed above is disjoint and that $\bc (h^*(A))_k=  h^*(\bc A)_k$.

To see that the sets in \ref{equation-c-preserv2} are disjoint, we again assume to the contrary that there is $x\in P_{k-1}$ belonging to both of them, i.e. there are $z\in A_{k+1}-\ker(h)$  and $t,t'\in\delta(z)-\ker(h)$ with
$\gamma(h(t))=x\in\delta(h(t'))$ and, moreover, there is $y\in A_{k-1}-\iota(A_{k+1})$ such that $h(y)=x$.

Since $h$ is a $\iota$-map, there is $s'\in\delta(t')$ such that $h(s')=x$. Then $h(\gamma(t))=\gamma(h(t))=x=h(s')=h(y)$ and hence the faces $\gamma(t)$, $s'$ and $y$ are linearly ordered by $\leq^+$. We shall show that each ordering of these faces leads to a contradiction.

If $\gamma(t)\leq^+y\leq s'$, then $y\in \iota(A_{k+1})$ contrary to the supposition.

If $y<^+\gamma(t)$, then, since $y\not\in \iota(A_{k+1})$, we have that $y\not\in \gamma(A_{k+1})$. Thus we have $y\leq^+ s\in \delta(t)$, for some $s$, and then
\[ x=h(y)\leq^+ h(s)<^+h(\gamma(t))=x\]
contradicting strictness of $<^+$.

If $s'<^+ y$, then again, since $y\not\in \iota(A_{k+1})$, we have that $y\not\in \delta(A_{k+1})$. Hence $\gamma(t')\leq^+y$ and we have
\[ x=h(s')<^+h(\gamma(t'))\leq^+h(y)=x \]
contradicting again strictness of $<^+$. This means that the sets \ref{equation-c-preserv2} are disjoint, as well, and that $\bc (h^*(A))_{k-1}=  h^*(\bc A)_{k-1}$ also holds. Thus 2. holds, as well.
$\Box$\vskip 2mm

\begin{proposition}\label{prop-iota-emdedding}
Let $h : Q\ra P$ be a contraction morphism of positive opetopes. Then \mbox{$h^*:Q^*\ra P^*$} is an $\o$-functor.
\end{proposition}
{\em Proof.}
We shall show by simultaneous induction on the size of opetopic cardinals $A$ and $B$ contained in the positive opetope $Q$ that
\begin{enumerate}
  \item $h^*(A)$ is an opetopic cardinal contained in $P$;
  \item $h^*$ preserves domains and codomains of $A$, i.e. for any $k\in \o$
  \[ h^\ast(\bd^{(k)}A) =\bd^{(k)}h^\ast(A), \;\; h^\ast(\bc^{(k)}A) =\bc^{(k)}h^\ast(A); \]
  \item $h^*$ preserves composition $A+_kB$ of $A$ and $B$, if it is well defined, i.e.
  \[ h^\ast(A+_kB) =h^\ast(A)+_k h^\ast(B). \]
\end{enumerate}

If 1. holds for any opetopic cardinal $A$, then, by \ref{lemma-iota-emdedding}.2, condition 2. holds for $A$, as well.

Next note that if 1. and 2. holds for both $A$ and $B$ with $\bc^{(k)}(A)=\bd^{(k)}(B)$, then we have
\[ h^\ast(A+_kB) = h^\ast(A\cap B) = h^\ast(A)\cap h^\ast(B) =  h^\ast(A)+_k h^\ast(B)\]
where the last equation holds by 2. In particular, $h^\ast(A+_kB)$ is an opetopic cardinal contained in $P$.

Thus we are left with 1 to justify. If $A$ is a positive opetope, then 1. holds by \ref{lemma-iota-emdedding}.1.
If $A$ is not an opetope, then there are two positive opetopic cardinals of a smaller size than $A$ such that $A=B+_kC$.
We can assume that for both $B$ and $C$ the conditions 1. and 2. hold. We have
\[ h^\ast(A) = h^\ast(B+_kC) =  h^\ast(B)+_k h^\ast(C).\]
Thus  $h^\ast(A)$ is a composition in $P^*$ of two opetopic cardinals and hence it is an opetopic cardinal, as well.
$\Box$\vskip 2mm

\begin{theorem}\label{theorem-iota-emdedding}The functor
\[ (-)^*: \pOpei \lra \oCat\]
is well defined. The objects of $\pOpei$ are sent under $(-)^*$ to positive-to-one polygraphs.  $(-)^*$ is faithful, conservative and full on those $\o$-functors that send generators to either generators or to (possibly iterated) identities on generators of smaller dimensions. In particular, it is full on isomorphisms.
\end{theorem}
{\em Proof.} All of this is either already proven or obvious except the fullness.

Recall that the generators of $Q^*$ are of form $Q[q]$, for $q\in Q$. To see that $(-)^* $ is indeed full on those $\o$-functors that send generators to either generators or to (possibly iterated) identities on generators of smaller dimensions, note that any $\o$-functor $f:Q^*\ra P^*$ can be restricted to a function $\bar{f}:Q\ra P^*$ such that $\bar{f}(q) =f(Q[q])$, for $q\in Q$. If, moreover, $f$ satisfies the additional hypothesis, then for any $q\in Q$ there is $p\in P$ such that $f(Q[q])=P[p]$. Clearly such a $p$ is unique and we get a function $h:Q\ra P$ such that $f(Q[q])=P[h(q)]$. This shows that $f$ agrees with $h^*$ on generators, and hence $f=h^*$, as required.
$\Box$\vskip 2mm

\begin{proposition}
Let $h: Q\ra P$ be a $\iota$-map, $q\in Q_n$, $l\leq n$. Then the restriction of $h$ to $h_q : Q[q]\ra P[h(q)]$ is an onto $\iota$-map.
If $\gamma^{(l)}(q)\not\in \ker(h)$, then $h(\gamma^{(l)}(q) )$ is $<^+$-maximal in $P[h(q)]_l$.
\end{proposition}
{\em Proof.} This is an easy consequence of Proposition \ref{lemma-iota-emdedding}.1 and Corollary \ref{coro-iota-preservation}.2.
$\Box$\vskip 2mm

We end this section with a Lemma that we will need later on.

\begin{lemma}\label{lemma-2-collaps}
Let $h: Q\ra P$ be a $\iota$-map, $k\geq 0$, $q\in Q_{k+2}$, $q'\in \delta(q)$,  $h(q),h(q')\in P_k$. Then $h(q')=h(\gamma(q))$.
\end{lemma}
{\em Proof.} Let $q'=q_0,\ldots, q_r$ be a lower $\delta(q)$-path such that $\gamma(q_r)=\gamma\gamma(q)$. Since $q$ is 2-collapsing, all faces in this path are at least 1-collapsing.  Since $\gamma(q')=\gamma(q_0)$ is non-collapsing, we have that face in $\gamma(q_i)\in\delta(q_{i+1})$ is the unique non-collapsing face in $q_{i+1}$, for $i=0,\ldots,r-1$. Hence, using induction, we get
\[ h(q')=h(\gamma(q'))=h(\gamma(q_0))=h(\gamma(q_r))=h(\gamma\gamma(q))=h(\gamma(q))\]
as required.
$\Box$\vskip 2mm

\subsection{The embedding of $\pHgi$ into $\widehat{\pOpei}$}

The notion of a $\iota$-map $h:Q\ra P$ was defined as a map between positive opetopes $Q$ and $P$ but it makes sense as a map between opetopic cardinals and even any hypergraph $Q$ and $P$. The category of hypergraphs with $\iota$-maps as morphisms will be denoted by $\pHgi$. Since $\pOpei\hookrightarrow \pHgi$, we have an embedding
\[   \cH_\iota:  \pHgi \lra \widehat{\pOpei} \]
such that
\[ H\mapsto {[ P \mapsto  \pHgi(P,H)]} \]
for any hypergraph $H$.

As positive hypergraphs can be very far from being positive opetopes, e.g. they do not need to satisfy any globularity conditions, it is not true that $\cH_\iota$ is either full or faithful. We call a positive hypergraph $H$ a {\em positive opetopic hypergraphs} iff for each face $x$ in $H$, the hypergraph $H[x]$ generated by $x$ in $H$ is a positive opetope.

\begin{proposition}\label{coro-embedding-cHi}
The functor $\cHi$ induces a bijection between hom-sets
\[ \cHi : \pHgi(P,H) \lra \widehat{\pOpei}(\cHi(P),\cHi(H)) \]
where $P$ is a positive opetope and $H$ is a positive opetopic hypergraph.
\end{proposition}

{\em Proof.}  Any $\iota$-map $h:P\ra H$ between a positive opetope $P$ and a positive opetopic hypergraph restricts to an onto $\iota$-map $h:P\ra H[h(\bm_P)$.  By Yoneda lemma, $\cHi$ is full and faithful on such $\iota$-maps.  $\Box$\vskip 2mm

We have a commuting square of categories and functors

\begin{center} \xext=1600 \yext=600
\begin{picture}(\xext,\yext)(\xoff,\yoff)
\setsqparms[1`1`1`1;800`500]
\putsquare(0,50)[\pOpe`\pHg`\pOpei`\pHgi;`\kappa`\kappa_h`]
\setsqparms[1`1`1`1;800`500]
\putsquare(800,50)[\phantom{\pHg}`\widehat{\pOpe}`\phantom{\pHgi}`\widehat{\pOpei};\cH``\kappa_!`\cHi]
 \end{picture}
\end{center}

\section{The cylinder functor on $\widehat{\pOpe}$}\label{sec-cylinder-in-pOpe}
\subsection{Flags}
Let $P$ be a positive opetope of dimension $n$ fixed for the whole section.

\begin{enumerate}
  \item The faces of $P$ can be described as follows:
  \begin{enumerate}
    \item  $P_n=\{ \bm_P\}$;
    \item $P_{n-1}=\{ \gamma(\bm_P)\}\cup \delta(\bm_P)$;
    \item  and, for each $0\leq k\leq n-2$, we have
  \[ P_k=\gamma^{(k)}(\bm_P)\cup \iota \gamma^{(k+2)}(\bm_P)\cup \delta\gamma^{(k+1)}(\bm_P).\]
  \end{enumerate}

  As usual we ignore braces whenever possible. We call a face $x\in P_k$ of
  \begin{enumerate}
    \item {\em $\gamma$-type} iff $x=\gamma^{(k)}(\bm_P)$;
    \item {\em $\iota$-type} iff $x\in  \iota\gamma^{(k+2)}(\bm_P)$;
    \item {\em $\delta$-type} iff $x\in\delta\gamma^{(k+1)}(\bm_P)$.
  \end{enumerate}

  \item A {\em restricted flag} or an {\em r-flag} in $P$ is a sequence of faces in $P$ $$\left[
                                      \begin{array}{c}
                                        x_m \\
                                        x_{m-1} \\
                                        \ldots \\
                                        x_l \\
                                      \end{array}
                                    \right]$$
  written also as $\vec{x}=[x_m, x_{m-1},\ldots,x_l]$, $0\leq l\leq m\leq n$ where $x_i\in P_i$, for $i=l,\ldots, m$ and $x_i\in \partial(x_{i+1})$, for $i=l,\ldots, m-1$. The subscript indicates the dimension of a face in the flag.   We write $\vec{x}(i)=x_i$.  We also say that such a flag is {\em under} $x_m$ and {\em over} $x_l$.  If $l=0$, we call such an r-flag a {\em flag}. If $k=0$ and $m=n$ (and hence $x_n=\bm_P$), then we call such a flag a maximal flag.

  \item The {\em sign of an (r-)flag} is defined as follows
   \[ \sgn([x_k,x_{k-1},\ldots, x_l]) \;= \; \left\{ \begin{array}{ll}
		1   & \mbox{ if } k=l \\
		\sgn([x_{k-1},\ldots, x_l])   & \mbox{ if }  k>l\mbox{ and }  x_{k-1}=\gamma(x_k) \\
         (-1)\cdot\sgn([x_{k-1},\ldots, x_l])   & \mbox{ if }  k> l\mbox{ and }  x_{k-1}\in \delta(x_k) \\
                                    \end{array}
			    \right. \]
\item We introduce notation for some sets of (r-)flags.  We assume below that $0\leq l \leq k\leq n$, $x_i\in P_i$:
    \begin{enumerate}
      \item $\Flags_P[x_k / x_l]$ is the set of r-flags under $x_k$ and over $x_l$.
      \item  $\Flags_P[x_k / 0)$ is the set of flags under $x_k$.
    \end{enumerate}
\item We define the {\em pencil order} $\prec_x$ on faces in the pencil over $x\in P_k$, i.e. on the set $\Pl_x=\{a\in P_{k+1}| x\in \partial(a) \}$.
If $a,b\in \Pl_x$, then we put $a \prec_x b$ iff
   \begin{enumerate}
    \item either  $\gamma(b)=x\in\delta(a)$,
     \item or $x\in\delta(a)\cap\delta(b)$, and $a<^+ b$,
    \item or  $\gamma(a)=x=\gamma(b)$, and $b<^+ a$.
  \end{enumerate}
  We write $\preceq_x$ for the reflexive closure of $\prec_x$ and $\prec^{op}_x$ for the reverse of the order $\prec_x$. Clearly, both $\prec_x$ and $\prec^{op}_x$ are strict linear orders.  $\next_{\preceq_x}$ is the successor function on pencil $\Pl_x$.

\item We define a {\em flag order} $\lhd$ on sets of flags $\Flags_P[x_m / x_l]$ and $\Flags_P[x_m / 0)$, where $0\leq l<m\leq dim(P)$.  Let $\vec{x}$ and $\vec{y}$ be two different flags in either set. Let $\diff(\vec{x},\vec{y})=\min\{ j | x_j\neq y_j \}=k$. We put
     $\vec{x}\lhd\vec{y}$ iff
    \begin{enumerate}
      \item $k=0$, and $y_0<^+ x_0$,
      \item or $k>0$, $\sgn([x_{k-1},\ldots,x_l])=1$ and $x_k \prec_{x_{k-1}} y_k$
      \item or $k>0$, $\sgn([x_{k-1},\ldots,x_l])=-1$ and $y_k \prec_{x_{k-1}} x_k$
    \end{enumerate}
    (if $\vec{x},\vec{y}\in\Flags_P[x_m / 0)$, then  $l=0$ in the above definition).
    By pencil linearity, the flag order is a strict linear order on $\Flags_P[x_m / x_l]$ and $\Flags_P[x_m / 0)$. Its reflexive closure will be denoted by $\unlhd$.
 \end{enumerate}

  \subsection{Punctured flags}
  \begin{enumerate}
  \item If a flag $\vec{x}=[x_k,\ldots,x_0]$ is neither the first nor the last flag and $k>1$, we can define the {\em low level of a flag} $\vec{x}$ as
  $$\ll(\vec{x}) =\max\{ i<k-1 \,| x_{i+1}\in \delta(x_{i+2}) \}.$$
  \item If $\vec{x}= [ x_k,\ldots,x_0]$ is a flag in $P$, then by $\vec{x}_{(i)}= [ x_k,\ldots, x_{i+1},0, x_{i-1},\ldots,x_0]$ we denote a flag with one face replaced by a dummy node $0$.

  By a {\em punctured flag} or a {\em p-flag} we mean $\vec{z}=\vec{x}_{(i)}$ in $P$  with either $i=k-1$ or $i=\ll(\vec{x})$, for some flag $\vec{x}$. If $i=k-1$, we say that $\vec{z}$ is a {\em high p-flag}, otherwise $\vec{z}$ is a {\em low p-flag}. If $\vec{z}$ is a p-flag with a dummy node at level $i$, then we write $\pu(\vec{z})=i$.

  If $\vec{x}$ is a flag or a high p-flag with top face $x_k$, then it makes sense to define $\vec{x}_{(k)}$. If $\vec{x}$ is a flag, then $\vec{x}_{(k)}= [ x_{k-1},\ldots, x_0]$ and if $\vec{x}$ is a high p-flag, then $\vec{x}_{(k)}= [ x_{k-2},\ldots,x_0]$. In both cases we get a flag again.

  \item By an intersection of two flags we mean a vector in which we put $0$ at each level where flags differ. After Corollary \ref{coro-top-path}, it will be clear that each p-flag is an intersection of two consecutive flags.

    \item $\pFlags_P[x_n,0)$  denote the sets of maximal p-flags in $P$.
    \item The function $$(-)_{high}: \Flags_P[x_n,0)\lra \pFlags_P[x_n,0)$$ is defined, for $\vec{x}\in\Flags_P[x_n,0)$, as
          \[ (\vec{x})_{high} = \vec{x}_{(n-1)}\]
          and the function
          $$(-)_{low}: \Flags_P[x_n,0)\lra \pFlags_P[x_n,0)\cup \{-x_n,+x_n\}$$
           is defined, for $\vec{x}\in\Flags_P[x_n,0)$, as
   \[ (\vec{x})_{low} =  \left\{
                                     \begin{array}{ll}
                                    -x_n & \mbox{if } \vec{x} \mbox{ is the initial flag,} \\
                                       & \\
                                     +x_n & \mbox{if } \vec{x} \mbox{ is the terminal flag,} \\
                                     \vec{x}_{(\ll(\vec{x}))}& \mbox{otherwise.}\\

                                    \end{array}
			    \right. \]
    \item On the set $\pFlags_P[x_k,0)$ we define a relation $\lhd_{x_k,0}$ as follows
    \[ \lhd_{x_k,0} = \{ \lk \vec{x}_{low},\vec{x}_{high}  \rk | \vec{x}\in \Flags_P[x_k,0),\; \sgn(\vec{x})=+1  \}  \cup \]
    \[ \{\lk \vec{x}_{high},\vec{x}_{low}  \rk |  \vec{x}\in \Flags_P[x_k,0),\; \sgn(\vec{x})=-1  \}  \]
    The relation $\lhd_{x_k}$ ($\unlhd_{x_k,0}$) is the transitive (reflexive) closure of $\lhd_{x_k,0}$.
\end{enumerate}

\begin{proposition} Let $P$ be a positive opetope, $x_k\in P_k$, $0\leq k\leq dim(P)$. Then
\begin{enumerate}
  \item The relation $\lhd_{x_k}$ is a strict linear order on $\pFlags_P[x_k,0)$.
  \item The functions $(-)_{high},(-)_{low}: \Flags_P[x_n,0)\lra \pFlags_P[x_n,0)\cup \{-x_n,+x_n\}$ preserve the (non-strict) order on (p-)flags.
\end{enumerate}
\end{proposition}

{\em Proof.} Exercise. $\Box$\vskip 2mm

\subsection{Initial and terminal flags}

The form of the initial and the terminal flags in $\Flags_P[x_n / x_l]$ depend on the type of the face $x_l$. The following Proposition provides a full description of initial and terminal flags in both  $\Flags_P[x_n / x_l]$ and  $\Flags_P[x_n / 0)$.

\begin{proposition}  \label{prop-init-term-flags}
Let $x_n=\bm_P$, $x_l\in P$.
\begin{enumerate}
\item The initial and terminal flags in $\Flags_P[x_n / 0)$ are

$$ \gamma{\rm -flag: }\;\;\; \bI_{P}=\left[
                                       \begin{array}{c}
                                         x_{n} \\
                                        \gamma(x_n)\\
                                        \ldots\\
                                        \gamma^{(1)}(x_n)\\
                                         \gamma^{(0)}(x_n) \\
                                       \end{array}\right]\hskip 5mm {\rm and }\hskip 5mm\delta{\rm -flag: }\;\;\;\bT_{P}=
                                       \left[
                                        \begin{array}{c}
                                         x_{n} \\
                                          \gamma(x_n)\\
                                        \ldots\\
                                         \gamma^{(1)}(x_n)\\
                                         \delta \gamma^{(1)}(x_n)
                                       \end{array}\right]$$
respectively.
\item The initial and terminal flags in $\Flags_P[x_n / x_l]$ are
\begin{enumerate}
  \item if $x_l=\gamma^{(l)}(x_n)$, i.e. if $x_l$ is of $\gamma$-type
  $$ \gamma\gamma{\rm -flag: }\;\; \bI_{x_l}= \left[
                                       \begin{array}{c}
                                         x_{n} \\
                                        \ldots\\
                                        \gamma^{(l+2)}(x_n)\\
                                        \gamma^{(l+1)}(x_n)\\
                                         x_l=\gamma^{(l)}(x_n) \\
                                       \end{array}\right]\hskip 5mm {\rm and }\hskip 5mm
                                       \gamma\delta{\rm -flag: }\;\; \bT_{x_l}=
                                       \left[
                                        \begin{array}{c}
                                         x_{n} \\
                                        \ldots\\
                                         \gamma^{(l+2)}(x_n)\\
                                        x_{l+1}\\
                                         x_l=\gamma^{(l)}(x_n) \\
                                       \end{array}\right]$$
  respectively, where $x_{l+1}\in \delta\gamma^{(l+2)}(x_n)$ is the unique face such that $\gamma(x_{l+1})=x_l$;
  \item if $x_l\in\iota\gamma^{(l+2)}(x_n)$ i.e. if $x_l$ is of $\iota$-type
    $$ \delta\delta{\rm -flag: }\;\;\; \bI_{x_l}= \left[
                                       \begin{array}{c}
                                         x_{n} \\
                                        \ldots\\
                                        \gamma^{(l+2)}(x_n)\\
                                        x_{l+1}\\
                                        x_l \\
                                       \end{array}\right]\hskip 5mm {\rm and }\hskip 5mm
                                       \gamma\delta{\rm -flag: }\;\;\; \bT_{x_l}=
                                       \left[
                                        \begin{array}{c}
                                         x_{n} \\
                                        \ldots\\
                                         \gamma^{(l+2)}(x_n)\\
                                        x'_{l+1}\\
                                        x_l \\
                                       \end{array}\right]$$
  respectively, where $x_{l+1}\in \delta\gamma^{(l+2)}(x_n)$ is the unique face such that $x_l\in \delta(x_{l+1})$
  and $x'_{l+1}\in \delta\gamma^{(l+2)}(x_n)$ is the unique face such that $\gamma(x'_{l+1})=x_l$,
  respectively;
  \item if $x_l\in\delta\gamma^{(l+1)}(x_n)$ i.e. if $x_l$ is of $\delta$-type
      $$ \delta\delta{\rm -flag: }\;\; \bI_{x_l}= \left[
                                       \begin{array}{c}
                                         x_{n} \\
                                        \ldots\\
                                        \gamma^{(l+2)}(x_n)\\
                                        x_{l+1}\\
                                        x_l \\
                                       \end{array}\right]\hskip 5mm {\rm and }\hskip 5mm
                                       \delta\gamma{\rm -flag: }\;\; \bT_{x_l}=
                                       \left[
                                        \begin{array}{c}
                                         x_{n} \\
                                        \ldots\\
                                         \gamma^{(l+2)}(x_n)\\
                                         \gamma^{(l+1)}(x_n)\\
                                        x_l \\
                                       \end{array}\right]$$
  respectively, where $x_{l+1}\in \delta\gamma^{(l+2)}(x_n)$ is the unique face such that $x_l\in \delta(x_{l+1})$
  and $x'_{l+1}\in \delta\gamma^{(l+2)}(x_n)$ is the unique face such that $\gamma(x_{l+1})=x_l$,
  respectively.
\end{enumerate}
\end{enumerate}
\end{proposition}

{\em Proof.} Ad 1. The endpoints in $\Flags_P[x_n / 0)$. Let $\vec{x},\vec{y},\vec{z}$ be three different flags in $\Flags_P[x_n / 0)$ with $\vec{x}$ a $\gamma$-flag and $\vec{z}$ a $\delta$-flag. We shall show that $\vec{x}\lhd\vec{y}\lhd\vec{z}$.

Let $k=\diff(\vec{x},\vec{y})$. Since $\vec{x}$ is a $\gamma$-flag, we have $y_k<^+x_k$. If $k=0$, then we get immediately that $\vec{x}\lhd\vec{y}$. If $k>0$, then $\gamma(x_k)=x_{k-1}=\gamma(y_k)$ and $\sgn(\vec{x}_{\lc k-1})=+1$. So in this case we have $x_k\prec_{x_{k-1}}y_k$ and we also have $\vec{x}\lhd\vec{y}$.

Let $m=\diff(\vec{y},\vec{z})$. If $m=0$ , as $\vec{z}$ is a $\delta$-flag, we have $z_0<^+y_0$. So we get in this case immediately that  $\vec{y}\lhd\vec{z}$. If $m=1$, then, as  $z_0\in\delta(y_1)\cap\delta(z_1)$ and $y_1<^+ z_1=\gamma^{(1)}(x_n)$, we have that $y_1\prec_{z_0}z_1$. Thus $\vec{y}\lhd\vec{z}$, as well. If $m>1$, then $\gamma(z_m )=z_{m-1}=\gamma(y_m)$ and $z_m\prec_{z_{m-1}}y_m$. We have $\sgn(\vec{z}_{\lc m-1})=-1$ and hence again $\vec{y}\lhd\vec{z}$.

Ad 2. The endpoints in $\Flags_P[x_n / x_l]$.  Let $\vec{x},\vec{y},\vec{z}$ be three different flags in $\Flags_P[x_n / x_l)$ with $\vec{x}=\bI_{x_l}$ and $\vec{z}=\bT_{x_l}$. We shall show that $\vec{x}\lhd\vec{y}\lhd\vec{z}$.

Let $k=\diff(\vec{x},\vec{y})$.  Let $m=\diff(\vec{y},\vec{z})$. We shall consider cases concerning the type of face $x_l$ and the relation of $k$ and $m$ to $l$.

Ad 2(a). $x_l$ is of $\gamma$-type, i.e. $x_l=\gamma^{(l)}(x_n)$. This case is similar to 1.

Ad 2(b). $x_l$ is of $\iota$-type, i.e. $x_l=\iota\gamma^{(l+2)}(x_n)$.  In this case $\vec{x}$ is a $\delta\delta$-flag.

If $k=l+1$, then $x_k\prec_{x_{k-1}} y_k$ and hence $\vec{x}\lhd\vec{y}$. If $k=l+2$, then, as $x_{k-1}$ is of $\delta$-type and $x_k=\gamma^{(k)}(x_n)$ is $<^+$-maximal, we have $y_k\prec_{x_{k-1}}x_k$. The sign $\sgn([x_{k-1},x_{k-2}])=-1$ and hence $\vec{x}\lhd\vec{y}$. If $k>l+2$, then $x_{k-1}$ is of $\gamma$-type, $\gamma(x_k)=x_{k-1}=\gamma(y_k)$, $y_k<^+x_k$ and hence $x_k\prec_{x_{k-1}}y_k$. The sign $\sgn(\vec{x}_{\lc k-1})=+1$ and hence $\vec{x}\lhd\vec{y}$ in this case as well. Thus $\vec{x}$ is initial indeed.

To see that the $\gamma\delta$-flag $\vec{z}$ is terminal, we again consider cases. If $m=l+1$, then $\gamma(z_m)=x_{m-1}$ and $z_m$ is $<^+$-minimal. Thus $y_m\prec_{z_{m-1}} z_m$ and hence  $\vec{y}\lhd\vec{z}$. If $m=l+2$, then, as $z_{m-1}$ is of $\delta$-type and $z_m$ is $<^+$-maximal, we have $y_m\prec_{z_{m-1}}z_m$. Since sign $\sgn([z_{m-1},z_{m-2}])=+1$, we have $\vec{y}\lhd\vec{z}$. If $m>l+2$, then $z_{m-1}$ is of $\gamma$-type and $z_m$ is $<^+$-maximal, so $z_m\prec_{z_{m-1}}y_m$. Since sing $\sgn(\vec{z}_{\lc m-1})=-1$, we have $\vec{y}\lhd\vec{z}$, as well.

  Ad 2(c). $x_l$ is of $\delta$-type, i.e. $x_l\in\delta\gamma^{(l+1)}(x_n)$.  In this case $\vec{x}$ is a $\delta\delta$-flag.

If $k=l+1$, then, as  $x_k<^+y_k$, we have $x_k\prec_{x_{k-1}}y_k$ and hence $\vec{x}\lhd\vec{y}$. If $k=l+2$, then, as $x_{k-1}$ is of $\delta$-type i.e. $<^+$-minimal and $x_k=\gamma^{(k)}(x_n)$ is $<^+$-maximal, we have $y_k\prec x_k$. But $\sgn([x_{k-1},x_{k-2}])=-1$ and hence $\vec{x}\lhd\vec{y}$ again. If $k>l+2$, then both $x_k$ and $x_{k-1}$ are of $\gamma$-type and hence $x_k\prec_{x_{k-1}}y_k$. Since $\sgn(\vec{x}_{\lc k-1}=+1$, we have $\vec{x}\lhd\vec{y}$, as well.

Now we assume that $\vec{z}$ is a $\delta\gamma$-flag. If $k=l+1$, then $z_{k-1}\in \delta(y_{k})\cap\delta(z_k)$. Since $y_k<^+z_k=\gamma^{(k)}(x_n)$, we have $\vec{y}\lhd\vec{z}$.

If $k>l+1$, then $\gamma(z_k)=z_{k-1}=\gamma(y_k)$ and $y_k<^+z_k$. Since $\sgn(\vec{z}_{\lc k-1})=-1$, we have $\vec{y}\lhd\vec{z}$, as well.
$\Box$\vskip 2mm

\vskip 2mm
\begin{proposition}\label{prop-succ-fag-order} Let $x_m\in P_m$ and $x_k\in P_k$ be two faces on $P$ and with $n\geq m>k\geq l\geq 0$.
The truncation morphisms
$$\tr : \Flags_P[x_m,x_l]\ra \Flags_P[x_k,x_l]$$
and
$$\tr : \Flags_P[x_m,0)\ra \Flags_P[x_k,0)$$
preserve the flag order and the endpoints.
\end{proposition}

{\em Proof.} Obvious. $\Box$\vskip 2mm

\subsection{Successor flags}

In this section we shall study the successor functions on sets  $\Flags_P[x_n / 0)$ with $n>1$ and on sets $\Flags_P[x_n / x_l]$ with $n>l+1$.

Let $\vec{x}, \vec{y}$ be two different flags in either set. We say that a flag $\vec{y}$ is a {\em $k$-th  neighbour of a flag} $\vec{x}$ (notation $\vec{y}=\neig_k(\vec{x}))$) iff whenever $\vec{y}(i)\neq\vec{x}(i)$, then $i=k$.
The {\em high neighbour of} $\vec{x}$ (notation $\neigh(\vec{x})$) is the $(n-1)$-th neighbour of $\vec{x}$. The {\em low neighbour of} $\vec{x}$ (notation $\neigh(\vec{x})$)  is the $\ll(\vec{x})$-th neighbour of $\vec{x}$.

Let $\vec{x}$ be a (non-terminal) flag in one of those sets. We define
       \[ \next(\vec{x}) \;= \; \left\{ \begin{array}{ll}
		 \neigh(\vec{x})  & \mbox{ if } \sgn(\vec{x})=+1 \\
		\neigl(\vec{x})  & \mbox{ if } \sgn(\vec{x})=-1 \\
                                    \end{array}
			    \right. \]
We may write $\next_{\lhd_{x_l}}(\vec{x})$ to indicate that it is a function on $\Flags_P[x_n / x_l]$.

For $x_l\in \partial(x_{l+1})$, we have a function extending flags by $x_l$
\[ \frac{(-)}{x_l}\cdot_{x_{l+1}} : \Flags[x_n,x_{l+1}]\lra  \Flags[x_n,x_{l}] \]
such that
\[ [x_n,\ldots,x_{l+1}] \mapsto  [x_n,\ldots,x_{l+1},x_l]. \]
Sometimes we write $\frac{(-)}{x_l}$  for $\frac{(-)}{x_l}\cdot_{x_{l+1}}$ if $x_{l+1}$ is understood. We have

\begin{proposition}\label{prop-path-ext}
\begin{enumerate}
  \item For $x_l=\gamma(x_{l+1})$, the function
  \[ {\frac{(-)}{x_l}} \cdot_{x_{l+1}} : \Flags[x_n,x_{l+1}]\lra  \Flags[x_n,x_{l}] \]
preserves the successor and hence the order, as well.

  \item For $x_l\in \delta(x_{l+1})$ the function
  \[ \frac{(-)}{x_l}\cdot_{x_{l+1}} : \Flags^{op}[x_n,x_{l+1}]\lra  \Flags[x_n,x_{l}] \]
preserves the successor and hence the order, as well.

  \item if $x_l\in \partial(x_{l+1})\cap \partial(x'_{l+1})$ and $x_{l+1}\neq x'_{l+1}$, then the images of functions $\frac{(-)}{x_l}\cdot_{x_{l+1}}$ and  $\frac{(-)}{x_l}\cdot_{x'_{l+1}}$ are disjoint.
  \item Let $x_l\in P_l$ and $x_{l+1},x'_{l+1}\in \Pl_{x_l}$ with $\next_{\prec_{x_l}}(x_{l+1})=x'_{l+1}$. Then
  \begin{enumerate}
    \item if $x_l\in\delta(x_{l+1})\cap \delta(x'_{l+1})$, then
    \[ \next_{\lhd_{x_l}}(\frac{\bI_{x_{l+1}}}{x_l}) = \frac{\bT_{\next_{\prec_{x_l}}(x_{l+1})}}{x_l} \]
       \item if $\gamma(x'_{l+1})=x_l\in\delta(x_{l+1})$, then
    \[ \next_{\lhd_{x_l}}(\frac{\bI_{x_{l+1}}}{x_l}) = \frac{\bI_{\next_{\prec_{x_l}}(x_{l+1})}}{x_l} \]
    \item if $\gamma(x'_{l+1})=x_l=\gamma(x_{l+1})$, then
    \[ \next_{\lhd_{x_l}}(\frac{\bT_{x_{l+1}}}{x_l}) = \frac{\bI_{\next_{\prec_{x_l}}(x_{l+1})}}{x_l} \]
  \end{enumerate}
\end{enumerate}
\end{proposition}

{\em Proof.}
Ad 1. Let $x_l=\gamma(x_{l+1})$. To see that the function
\[ {\frac{(-)}{x_l}} \cdot_{x_{l+1}} : \Flags[x_n,x_{l+1}]\lra  \Flags[x_n,x_{l}] \]
preserves the successor, notice that if $\vec{x}\in \Flags[x_n,x_{l+1}]$, then we have
\[ \neigh(\frac{\vec{x}}{x_l}) = \frac{\neigh(\vec{x})}{x_l} \]
and if $\vec{x}$ is not an endpoint flag and $n>l+2$, then
\[ \neigl(\frac{\vec{x}}{x_l}) = \frac{\neigl(\vec{x})}{x_l}. \]

Ad 2.     If $x_l\in \delta(x_{l+1})$, then  $\sgn([x_{l+1},x_l])=-1$ and the function
\[ {\frac{(-)}{x_l}} \cdot_{x_{l+1}} : \Flags^{op}[x_n,x_{l+1}]\lra  \Flags[x_n,x_l] \]
preserves the order. The rest is similar to 1.

Ad 3. If $x_{l+1}\neq x'_{l+1}$, then the values of $\frac{(-)}{x_l}\cdot_{x_{l+1}}$ and  $\frac{(-)}{x_l}\cdot_{x'_{l+1}}$ at level $l+1$ are $x_{l+1}$ and $x'_{l+1}$, respectively, and hence they are fixed for each function and different. Thus the images are disjoint.

Ad 4(a). Let $x_l\in\delta(x_{l+1})\cap \delta(x'_{l+1})$ and $\next_{\prec_{x_l}}(x_{l+1})=x'_{l+1}$.  In the inverse flag order on $\Flags^{op}[x_n,x_{l+1}]$ the flag $\bI_{x_{l+1}}$ is the terminal one, by Lemma \ref{prop-init-term-flags}. Thus, as $x_l\in \delta(x_{l+1})$, in the order $\Flags[x_n,x_l]$ the flag $\frac{\bI_{x_{l+1}}}{x_l}$ is the last flag having $x_{l+1}$ at the level $l+1$. Thus the next flag has to change $x_{l+1}$ to the next face $x'_{l+1}$ in the pencil order over $x_l$. We have
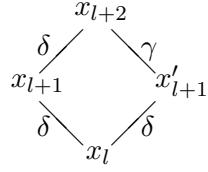
\begin{center}
\xext=500 \yext=500
\begin{picture}(\xext,\yext)(\xoff,\yoff)
   \put(220,500){$x_{l+2}$}
   \put(0,250){$x_{l+1}$}
   \put(500,250){$x'_{l+1}$}
   \put(260,0){$x_{l}$}

\put(350,450){\line(1,-1){150}}
\put(90,360){$\delta$}

\put(250,450){\line(-1,-1){150}}
\put(450,360){$\gamma$}

\put(100,200){\line(1,-1){150}}
\put(90,80){$\delta$}
\put(500,200){\line(-1,-1){150}}
\put(450,80){$\delta$}
\end{picture}
\end{center}
So
$$\vec{x}=[x_n,\ldots, x_{l+2},x_{l+1},x_l]=\frac{\bI_{x_{l+1}}}{x_l}$$
is $\lhd$-smaller than
$$\vec{x}'=[x_n,\ldots, x_{l+2},x'_{l+1},x_l]=\frac{\bT_{\next_{\prec_{x_l}}(x_{l+1})}}{x_l}.$$
We still need to show that there is nothing between these flags in flag order $\lhd_{x_l}$ on $\Flags[x_n,x_l]$. So we need to show that faces in $\vec{x}'$ over the level $l+1$ are the least possible. Over the level $l+2$ it is clear, as
\[ \sgn([x_{l+2},x_{l+1},x_l])=+1,\;\;\;\;{\rm and }\;\;\;\; \sgn([x_{l+2},x'_{l+1},x_l])=-1\]
 and we have a change of sign. Hence what was the largest face in $\vec{x}$ will become the least face in $\vec{x}'$. The element $x_{l+2}\in\delta\gamma^{(l+3)}(x_n)$ is $<^+$-minimal. Since $x_{l+1}\in\delta(x_{l+2})$, it is the least face in $\Pl_{x_{l+1}}$ and since $x'_{l+1}=\gamma(x_{l+2})$, it is the greatest face in $\Pl_{x'_{l+1}}$. We have
\[ \sgn([x_{l+1},x_l])=-1=\sgn([x'_{l+1},x_l]) \]
and the face $x_{l+2}$ is the largest possible over $[x_{l+1},x_l]$ and the least possible over $[x'_{l+1},x_l]$, as required.

Ad 4(b). Let $\gamma(x'_{l+1})=x_l\in\delta(x_{l+1})$ and $\next_{\prec_{x_l}}(x_{l+1})=x'_{l+1}$.  In the inverse flag order on $\Flags^{op}[x_n,x_{l+1}]$ the flag $\bI_{x_{l+1}}$ is the terminal one, by Lemma \ref{prop-init-term-flags}. Thus, as $x_l\in\delta(x_{l+1})$, in the order $\Flags[x_n,x_l]$ the flag $\frac{\bI_{x_{l+1}}}{x_l}$ is the last flag having $x_{l+1}$ at the level $l+1$. Thus the next flag has to change $x_{l+1}$ to the next face $x'_{l+1}$ in the pencil order over $x_l$. By assumption, we have
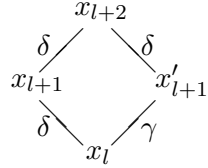
\begin{center}
\xext=500 \yext=500
\begin{picture}(\xext,\yext)(\xoff,\yoff)
   \put(220,500){$x_{l+2}$}
   \put(0,250){$x_{l+1}$}
   \put(500,250){$x'_{l+1}$}
   \put(260,0){$x_{l}$}

\put(350,450){\line(1,-1){150}}
\put(90,360){$\delta$}

\put(250,450){\line(-1,-1){150}}
\put(450,360){$\delta$}

\put(100,200){\line(1,-1){150}}
\put(90,80){$\delta$}
\put(500,200){\line(-1,-1){150}}
\put(450,80){$\gamma$}
\end{picture}
\end{center}
So
$$\vec{x}=[x_n,\ldots, x_{l+2},x_{l+1},x_l]=\frac{\bI_{x_{l+1}}}{x_l}$$
is $\lhd$-smaller than
$$\vec{x}'=[x_n,\ldots, x_{l+2},x'_{l+1},x_l]=\frac{\bI_{\next_{\prec_{x_l}}(x_{l+1})}}{x_l}.$$
We still need to show that there is nothing between these flags in flag order $\lhd_{x_l}$ on $\Flags[x_n,x_l]$. So we need to show that faces in $\vec{x}'$ over the level $l+1$ are the least possible. As before, over the level $l+2$ it is clear, as
\[ \sgn([x_{l+2},x_{l+1},x_l])=+1,\;\;\;\;{\rm and }\;\;\;\; \sgn([x_{l+2},x'_{l+1},x_l])=-1\]
 and we have a change of sign. Hence what was the largest face in $\vec{x}$ will become the least face in $\vec{x}'$. The element $x_{l+2}\in\delta\gamma^{(l+3)}(x_n)$ is $<^+$-minimal. Since $x_{l+1},x'_{l+1}\in\delta(x_{l+2})$, it is the least face in both  $\Pl_{x_{l+1}}$ and $\Pl_{x'_{l+1}}$.
 We have
\[ \sgn([x_{l+1},x_l])=-1\neq +1=\sgn([x'_{l+1},x_l]) \]
and hence the face $x_{l+2}$ is the largest possible over $[x_{l+1},x_l]$ and the least possible over $[x'_{l+1},x_l]$, as required.

Ad 4(c). Let $\gamma(x'_{l+1})=x_l=\gamma(x_{l+1})$ and $\next_{\prec_{x_l}}(x_{l+1})=x'_{l+1}$.  In the flag order on $\Flags[x_n,x_{l+1}]$ the flag $\bT_{x_{l+1}}$ is the terminal one, by Lemma \ref{prop-init-term-flags}. Thus, as $x_l=\gamma(x_{l+1})$, in the order $\Flags[x_n,x_l]$ the flag $\frac{\bT_{x_{l+1}}}{x_l}$ is the last flag having $x_{l+1}$ at the level $l+1$. Thus the next flag has to change $x_{l+1}$ to the next face $x'_{l+1}$ in the pencil order over $x_l$. By assumption, we have
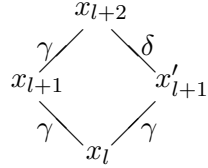
\begin{center}
\xext=500 \yext=500
\begin{picture}(\xext,\yext)(\xoff,\yoff)
   \put(220,500){$x_{l+2}$}
   \put(0,250){$x_{l+1}$}
   \put(500,250){$x'_{l+1}$}
   \put(260,0){$x_{l}$}

\put(350,450){\line(1,-1){150}}
\put(90,360){$\gamma$}

\put(250,450){\line(-1,-1){150}}
\put(450,360){$\delta$}

\put(100,200){\line(1,-1){150}}
\put(90,80){$\gamma$}
\put(500,200){\line(-1,-1){150}}
\put(450,80){$\gamma$}
\end{picture}
\end{center}
So
$$\vec{x}=[x_n,\ldots, x_{l+2},x_{l+1},x_l]=\frac{\bT_{x_{l+1}}}{x_l}$$
is $\lhd$-smaller than
$$\vec{x}'=[x_n,\ldots, x_{l+2},x'_{l+1},x_l]=\frac{\bI_{\next_{\prec_{x_l}}(x_{l+1})}}{x_l}.$$ We still need to show that there is nothing between these flags in flag order $\lhd_{x_l}$ on $\Flags[x_n,x_l]$. So we need to show that faces in $\vec{x}'$ over the level $l+1$ are the least possible. As before, over the level $l+2$ it is clear, as
\[ \sgn([x_{l+2},x_{l+1},x_l])=+1,\;\;\;\;{\rm and }\;\;\;\; \sgn([x_{l+2},x'_{l+1},x_l])=-1\]
 and we have a change of sign. Hence what was the largest face in $\vec{x}$ will become the least face in $\vec{x}'$. The element $x_{l+2}\in\delta\gamma^{(l+3)}(x_n)$ is $<^+$-minimal. Since $\gamma(x_{l+2})=x_{l+1}$ and $x'_{l+1}\in\delta(x_{l+2})$, $x_{l+2}$ is the greatest face in $\Pl_{x_{l+1}}$ and the least face in $\Pl_{x'_{l+1}}$.
 We have
\[ \sgn([x_{l+1},x_l])=-1\neq +1=\sgn([x'_{l+1},x_l]) \]
and hence the face $x_{l+2}$ is the largest possible over $[x_{l+1},x_l]$ and the least possible over $[x'_{l+1},x_l]$, as required.
$\Box$\vskip 2mm

From Proposition \ref{prop-path-ext} we get
\begin{corollary}\label{prop-path-of-flags}
Let $P$ be a positive opetope of dimension $n$, $0\leq l<n-1$, $x_l\in P_l$, $x_n=\bm_P\in P_n$.
The function $\next$ defined on the set $\Flags_P[x_n/x_l]$ is the successor function with respect to the flag order.
\end{corollary}
  {\em Proof.} The proof is by downward induction on $l$.

  Clearly, the Theorem holds for the case $l=n-2$.

Let $l<n-2$.  The set of functions
\[ \{ \frac{(-)}{x_l}\cdot_{x_{l+1}} : \Flags^{op}[x_n,x_{l+1}]\lra  \Flags[x_n,x_{l}] \;|\; x_l\in \delta(x_{l+1})\}  \cup \]
\[ \{ {\frac{(-)}{x_l}} \cdot_{x_{l+1}} : \Flags[x_n,x_{l+1}]\lra  \Flags[x_n,x_{l}]\;|\; x_l=\gamma(x_{l+1})\}\]
is jointly surjective with images disjoint, preserves flag order and successor, by Proposition \ref{prop-path-ext}.1,2,3. Then by Proposition \ref{prop-path-ext}.4, the successor function also agrees with flag order when we pass from the last flag in the image of one such function to the first flag in the image of the next function.
$\Box$\vskip 2mm

\begin{theorem}\label{thm-top-path}
Let $P$ be a positive opetope of dimension $n$, $x_n=\bm_P\in P_n$. The function $\next$ defined on the set $\Flags_P[x_n/0)$ is the successor function with respect to the flag order.
 \end{theorem}

{\em Proof.} We add temporarily two faces to opetope $P$ at dimension $-1$, say $\bs$ for source and $\bt$ for target, and, for all faces $x\in P_0$, we can define $\gamma(x)=\bt$ and $\delta(x)=\bs$. It can be easily verified that it is still an opetope if we start counting dimensions from -1 instead of 0 (alternatively we could shift up one dimension all faces in $P$). Then the statement of Theorem is a consequence of Proposition \ref{prop-path-of-flags} applied to the set of flags $\Flags[x_n,\bt]$.
$\Box$\vskip 2mm

From the above we get
 \begin{corollary}\label{coro-top-path}
The intersection of two consecutive flags is a p-flag. If the first of the two flags is positive, then the intersection is a high p-flag and if the first of the two flags is negative, then the intersection is a low p-flag.

  \end{corollary}
  {\em Proof.} This follows from the description of successor function on maximal flags. $\Box$\vskip 2mm

We shall make the description of the successor function from Proposition \ref{prop-path-ext} more explicit here. Let $\vec{x}$ and $\vec{x}'$ be two maximal flags in $P$ such that $\next(\vec{x})=\vec{x}'$, i.e. $\vec{x}'$ is the successor of $\vec{x}$. Then there is $0\leq k<n$ such that

\[ \vec{x} =          \left[
                                       \begin{array}{c}
                                         x_n\\
                                         \ldots \\
                                         \gamma^{(l+2)}(p)\\
                                         x_{k+1}\\
                                         x_k\\
                                         x_{k-1}\\
                                        \ldots \\
                                         x_0 \\
                                       \end{array}\right]
                                       \lhd\!\!\!\!\cdot
                                       \left[
                                       \begin{array}{c}
                                         x_n\\
                                         \ldots \\
                                         \gamma^{(l+2)}(p)\\
                                         x_{k+1}\\
                                         x'_k\\
                                         x_{k-1}\\
                                        \ldots \\
                                         x_0 \\
                                       \end{array}\right]
                                       = \vec{x}'
                                     \]
If $k=0$, then we assume that $x_{k-1}$ `is' the -1-dimensional face $\bt$ so that $\gamma(x)=\bt$, for $x\in P_0$. There are successors of six kinds.

They are related to the mutual relations of faces $x_{k-1}$, $x_k$, $x'_k$, $x_{k+1}$ and the sign of the flag $\vec{x}_{\lc k-1}=[x_{k-1},\ldots,x_0]$.
For the future reference we will name them all.

\vskip 2mm
{\bf Case} $\sgn(\vec{x}_{\lc k-1})=+1$.

In this case we always have $\sgn([x_{k+1},x_k,x_{k-1}])=+1$. Thus we have three possibilities:

\vskip 1mm
\noindent
 $\vec{x}'$ is a $\delta$-successor of $\vec{x}$ iff

\begin{center}
\xext=500 \yext=500
\begin{picture}(\xext,\yext)(\xoff,\yoff)
   \put(220,500){$x_{k+1}$}
   \put(0,250){$x_{k}$}
   \put(500,250){$x'_{k}$}
   \put(260,0){$x_{k-1}$}

\put(350,450){\line(1,-1){150}}
\put(90,360){$\delta$}

\put(250,450){\line(-1,-1){150}}
\put(450,360){$\gamma$}

\put(100,200){\line(1,-1){150}}
\put(90,80){$\delta$}
\put(500,200){\line(-1,-1){150}}
\put(450,80){$\delta$}
\end{picture}
\end{center}
\vskip 1mm
$\vec{x}'$ is a $\delta\gamma$-successor of $\vec{x}$ iff

\begin{center}
\xext=500 \yext=500
\begin{picture}(\xext,\yext)(\xoff,\yoff)
   \put(220,500){$x_{k+1}$}
   \put(0,250){$x_{k}$}
   \put(500,250){$x'_{k}$}
   \put(260,0){$x_{k-1}$}

\put(350,450){\line(1,-1){150}}
\put(90,360){$\delta$}

\put(250,450){\line(-1,-1){150}}
\put(450,360){$\delta$}

\put(100,200){\line(1,-1){150}}
\put(90,80){$\delta$}
\put(500,200){\line(-1,-1){150}}
\put(450,80){$\gamma$}
\end{picture}
\end{center}
\vskip 1mm
$\vec{x}'$ is a $\gamma$-successor of $\vec{x}$ iff

\begin{center}
\xext=500 \yext=500
\begin{picture}(\xext,\yext)(\xoff,\yoff)
   \put(220,500){$x_{k+1}$}
   \put(0,250){$x_{k}$}
   \put(500,250){$x'_{k}$}
   \put(260,0){$x_{k-1}$}

\put(350,450){\line(1,-1){150}}
\put(90,360){$\gamma$}

\put(250,450){\line(-1,-1){150}}
\put(450,360){$\delta$}

\put(100,200){\line(1,-1){150}}
\put(90,80){$\gamma$}
\put(500,200){\line(-1,-1){150}}
\put(450,80){$\gamma$}
\end{picture}
\end{center}

\vskip 2mm
{\bf Case} $\sgn(\vec{x}_{\lc k-1})=-1$.

In this case we always have $\sgn([x_{k+1},x_k,x_{k-1}])=-1$. Thus we have three possibilities:

\vskip 1mm
\noindent
 $\vec{x}'$ is an inverse $\delta$-successor of $\vec{x}$ iff

\begin{center}
\xext=500 \yext=500
\begin{picture}(\xext,\yext)(\xoff,\yoff)
   \put(220,500){$x_{k+1}$}
   \put(0,250){$x_{k}$}
   \put(500,250){$x'_{k}$}
   \put(260,0){$x_{k-1}$}

\put(350,450){\line(1,-1){150}}
\put(90,360){$\gamma$}

\put(250,450){\line(-1,-1){150}}
\put(450,360){$\delta$}

\put(100,200){\line(1,-1){150}}
\put(90,80){$\delta$}
\put(500,200){\line(-1,-1){150}}
\put(450,80){$\delta$}
\end{picture}
\end{center}
\vskip 1mm
$\vec{x}'$ is an inverse $\gamma\delta$-successor of $\vec{x}$ iff

\begin{center}
\xext=500 \yext=500
\begin{picture}(\xext,\yext)(\xoff,\yoff)
   \put(220,500){$x_{k+1}$}
   \put(0,250){$x_{k}$}
   \put(500,250){$x'_{k}$}
   \put(260,0){$x_{k-1}$}

\put(350,450){\line(1,-1){150}}
\put(90,360){$\delta$}

\put(250,450){\line(-1,-1){150}}
\put(450,360){$\delta$}

\put(100,200){\line(1,-1){150}}
\put(90,80){$\gamma$}
\put(500,200){\line(-1,-1){150}}
\put(450,80){$\delta$}
\end{picture}
\end{center}
\vskip 1mm
$\vec{x}'$ is an inverse $\gamma$-successor of $\vec{x}$ iff

\begin{center}
\xext=500 \yext=500
\begin{picture}(\xext,\yext)(\xoff,\yoff)
   \put(220,500){$x_{k+1}$}
   \put(0,250){$x_{k}$}
   \put(500,250){$x'_{k}$}
   \put(260,0){$x_{k-1}$}

\put(350,450){\line(1,-1){150}}
\put(90,360){$\delta$}

\put(250,450){\line(-1,-1){150}}
\put(450,360){$\gamma$}

\put(100,200){\line(1,-1){150}}
\put(90,80){$\gamma$}
\put(500,200){\line(-1,-1){150}}
\put(450,80){$\gamma$}
\end{picture}
\end{center}

The successors can be still further divided into two classes: {\em high successors} if $k=n-1$ and {\em low successors} if $k<n-1$.

\subsection{The faces of cylinder $\Cyl(P)$ and star operation}
  Faces of {\em the cylinder set} $\Cyl(P)$ are of three kinds:
    \begin{enumerate}
      \item {\em flat faces}:  $\{-\}\times P \cup\{+\}\times P$;
      \item {\em flags}: all flags of all faces of $P$;
      \item {\em p-flags}: all p-flags of all faces in $P$.
    \end{enumerate}
The dimension of a face $-p$ or $+p$ is the dimension of $p$. The dimension of a flag or p-flag is the number of non-zero faces in the sequence.
$\Cyl(P)_k$ is the set of all faces of $\Cyl(P)$ of dimension $k$.

We have a projection function $\pi_P : \Cyl(P) \lra P$ such that, for a face $\varphi\in \Cyl(P)$,
   \[ \pi_P(\varphi) \;= \; \left\{ \begin{array}{ll}
		x   & \mbox{ if } \varphi \in\{  -x, +x\},  \\
		x_k   & \mbox{ if } \varphi = x_k,\ldots,x_0  \\
                                    \end{array}
			    \right. \]

We define a {\em star operation}, an `inverse' operation to  projection on flags
\[    \star: P\times \Flags[\bm_P,0) \ra \Cyl(P)\]
so that, for  $\vec{x}\in  \Flags[\bm_P,0)$ and $p\in P_k$, we have $\pi_P(p\star \vec{x})=p$. We put

\[ p\star\vec{x} =  \left\{
                                    \begin{array}{ll}
                                         \vec{x}_{\lc k}   & \mbox{if }  p=x_k;\\
                                       & \mbox{(flag)}\\
                                       &\\
                                       {[p,0,x_{k-2},\ldots,x_0]} &  \begin{array}{l}
                                       \mbox{otherwise, if }  k\geq 1 \mbox{and } x_{k}<^+  p;\\
                                        \mbox{(high p-flag)}\\
                                          \end{array}\\
                           &   \\

                {[p,\gamma(p),\ldots,\gamma^{(l+2)}(p),t,0,x_{l-1},\ldots,x_0]} &                                       \begin{array}{l}
                                       \mbox{otherwise, if }  k\geq 1,\\
                                       0\leq l=max(\{ l'\leq k-2| \; x_{l'+1}\leq^+  \delta\gamma^{(l'+2)}(p) \} \\
                                       \mbox{is well defined, and } \\
                                       t \mbox{ is such that }   x_{l+1}\leq^+ t \in \delta\gamma^{(l+2)}(p); \\
                                       \mbox{(low p-flag)}\\
                                       \end{array}\\

              &  \\
		 -p   & \mbox{otherwise, if }\gamma^{(0)}(p)\leq^+x_0;  \\
               & \mbox{(bottom flat face)}\\
               & \\
		 +p  & \mbox{otherwise, if } x_0<^+\gamma^{(0)}(p).\\
            & \mbox{(top flat face)}\\
                                    \end{array}
			    \right. \]

We define also an `inverse' operation to  projection on p-flags
\[    \star: P\times \pFlags[\bm_P,0) \ra \Cyl(P)\]
so that, for  $\vec{x}=[x_n,\ldots,\widehat{x_l},\ldots,x_0]\in  \pFlags[\bm_P,0)$, $p\in P_k$, we have

\[ p\star\vec{x} =  \left\{
                                    \begin{array}{lll}
                                       \vec{x}_{\lc k}   & \begin{array}{l}
                                       \mbox{if }   k\neq l,\;  p=x_k;\\
                                       \mbox{((p-)flag)}
                                       \end{array}\\
                            & \\
                             {[p,\ldots,\gamma^{(l+2)}(p),t,0,x_{l-1},\ldots,x_0]}   & \begin{array}{l}
                                       \mbox{otherwise, if } \\
                                       k>l\geq 0 \mbox{ and } x_{l+1}\leq^+ t\in \delta\gamma^{(l+2)}(p);\\
                                       \mbox{(p-flag with puncture at the level l)}
                                       \end{array}\\
               & \\
              {[p,\ldots,\gamma^{(l+1)}(p),t,0,x_{l-2},\ldots,x_0]}   & \begin{array}{l}
                                       \mbox{otherwise, if } \\
                                        k\geq l> 0 \mbox{ and } \gamma(x_{l+1})\leq^+ t\in\delta\gamma^{(l+1)}(p);\\
                                       \mbox{(p-flag with puncture at the level l-1)}
                                       \end{array}\\
             & \\
                  {[p,\ldots,\gamma^{(\bar{l}+2)}(p),t,0,x_{\bar{l}-1},\ldots,x_0]}    & \begin{array}{l}
                                       \mbox{otherwise, if } \\
                                        0\leq \bar{l}=\max\{l'\leq  i-2,k: x_{l'+1}\leq^+\delta\gamma^{(l'+2)}(p)   \} \\
                                       \mbox{ and }  x_{\bar{l}+1}\leq^+ t\in \delta\gamma^{(\bar{l}+2)}(p); \\
                                       \mbox{(p-flag with puncture below the level $\bar{l})$}
                                       \end{array}\\
                           &  \\
                    -p   &       \begin{array}{l}
                                       \mbox{otherwise, if } \\
                                       \mbox{either }l>0  \mbox{ and } \gamma^{(0)}(p)\leq^+x_0\\
                                       \mbox{or } l=0  \mbox{ and }  \gamma^{(0)}(p)<^+\gamma(x_1) ; \\
                                       \mbox{(bottom flat face)}\\
                                       \end{array}  \\
               & \\
	+p   &                                        \begin{array}{l}
                                       \mbox{otherwise, if } \\
                                       \mbox{either }l>0 \mbox{ and }x_0<^+ \gamma^{(0)}(p)\\
                                       \mbox{or } l=0\mbox{ and }  \gamma(x_1)\leq^+\gamma^{(0)}(p). \\
                                       \mbox{(top flat face)}\\
                                       \end{array}  \\
                & \\
                                    \end{array}
			    \right. \]

We can define a set $P\otimes \Cyl(P)$ to be the set of pairs $\lk p,\varphi\rk$ such that $p\in P_k$ occurs in $\pi_P(\varphi)\in P$. Then we can treat $\star$ as an operation on
\[    \star: P\otimes \Cyl(P) \ra \Cyl(P).\]
For flags and p-flags, it is defined as above and, for the flat faces, we extend the definition of $\star$ as follows. If $p\in P$ occurs in $x\in P$, then
\[ p\star -x = -p,\hskip 5mm {\rm and } \hskip 5mm  p\star +x = +p.\]

\begin{lemma} We have
\begin{enumerate}
  \item For $\lk p,\varphi\rk\in P\otimes\Cyl(P)$, we have $$\pi_P(p\star\varphi)=p.$$
  \item If $p\in P_{2}$, $x_1\in P_1$,  $x_0\in \partial(x_1)$ and $x_1\not\leq^+ \delta(p)$, then $x_0<^+\gamma\gamma(p)$ is equivalent to $x_0\leq^+\delta\gamma(p)$
  \item If $p\in P_{\geq 1}$, then in the last case of the above definition the condition $x_0<^+\gamma^{(0)}(p)$ is equivalent to $x_0\leq^+\delta\gamma^{(1)}(p)$.
\end{enumerate}
\end{lemma}

{\em Proof.} 1. is obvious.

Ad 2. Clearly, if $x_0\leq^+ \delta\gamma(p)$, then  $x_0<^+ \gamma\gamma(p)$. We need to show the converse.

So assume $x_0<^+ \gamma\gamma(p)$. Since $x_1\not\leq^+ \delta(p)$, then a fortiori $x_1\not\leq^+ \gamma(p)$. But both $x_1$ and $\gamma(p)$ are one cells, so by Corollary 5.12 of \cite{Z1}, we must have that  $x_1\perp^-\gamma(p)$. If we were to have $\gamma(p)<^-x_1$, then we would have $\gamma\gamma(p)\leq^+\delta(x_1)\leq^+x_0$, contrary to the supposition. On the other hand, if $x_1<^-\gamma(p)$, then $x_0\leq^+\gamma(x_1)\leq^+\delta\gamma(p)$, as required.

Ad 3. Apply 2. to $p$ equal $\gamma^{(2)}(p)$.
 $\Box$\vskip 2mm

  \begin{lemma}
  \label{lemma-flag-tech-J}
Let $P$ be a positive opetope, $k\geq 2$, $z\in P_{k-1}$, $x,y\in k$. Then
\begin{enumerate}
  \item if $x\not\perp^+ y$, $t\in\delta(x)$, $s\in\delta(y)$ and $t\leq^+ s$, then $\gamma(x)\leq^+s$;
  \item  if $z \not\perp^+ \delta(y)$, then, for each $t\in \delta(z)$, if there is $s\in\delta\gamma(y)$ such that $t\leq^+s$, then $\gamma(z)\leq^+s$.
  \end{enumerate}
  \end{lemma}
  {\em Proof.} Ad 1. Assume $x\not\perp^+ y$, $t\in\delta(x)$, $s\in\delta(y)$ and $t\leq^+ s$. If $s=t$, then $x\perp^+y$ and we get a contradiction. Thus $t<^+s$. Let $t,a_1,\ldots,a_l,s$ be a non-empty upper path in $P-\gamma(P)$ from $t$ to $s$.

  If $a_1=x$, then $\gamma(x),a_2,\ldots, s$ is an upper path in $P$ and hence $\gamma(x)\leq^+s$.

  If $a_1\neq x$, then $a_1<^+x$ and by Path Lemma \cite{Z1} either we get $y<^+x$ and we get again a contradiction or $x<^-y$. In the latter case we have
  \[ t\leq^+ \gamma(x)\leq^+ s'  \]
  for some $s' \in \delta(y)$. Since $t\leq^+ s$ and $s,s'\in\delta(y)$, we have $s=s'$ and $\gamma(x)\leq^+s$.

  Ad 2. Assume $z \not\perp^+ \delta(y)$. Fix $t\in \delta(z)$ and $s\in\delta\gamma(y)$ such that $t\leq^+s$. If $s=t$, then $z\perp^+\gamma(y)$.

  If $\gamma(y)\leq^+z$, then $\delta(y)\leq^+z$ and we get a contradiction. If $z<^+\gamma(y)$, then we have a non-empty upper $P-\gamma(P)$-path $z,a_1,\ldots,a_l,\gamma(y)$. If $y=a_l$, then $z\leq \delta(y)$ and we get a contradiction. If $a_l<^+y$, then by Path Lemma \cite{Z1} we get $z\perp^+\delta(y)$ and we get again a contradiction.

  Thus $t<^+s$ and we have a non-empty upper $P-\gamma(P)$-path $t,a_1,\ldots,a_l,s$. If $z=a_1$, then $\gamma(z)\leq^+s$, as required. If $a_1<^+z$, then, by Path Lemma \cite{Z1}, either $\gamma(y)<^+z$ and we get a contradiction, or $z<^-\gamma(y)$. In the latter case $\gamma(z)\leq^+s'\in \delta\gamma(y)$. Since $s,s'\in \delta\gamma(y)$, $t\leq^+ s$ and $t\leq^+ s'$, it follows that $s=s'$ and $\gamma(z)\leq^+s$.
  $\Box$\vskip 2mm

\begin{lemma}[flat values of $\star$]\label{lemma-flat-star}
Let $\vec{x}=[x_n,\ldots, x_0]$ be a flag or p-flag in $P$ with $x_n=\bm_P$.  Let $k<m$, $p'\in P_k-\{x_k\}$,  $p\in P_m-\{x_m\}$ so that $p'\inn p\inn x_n$. Then we have
\begin{enumerate}
  \item if $p\star\vec{x}=-p$, then $p'\star\vec{x}=-p'$;
  \item if $p\star\vec{x}=+p$, then $p'\star\vec{x}=+p'$.
\end{enumerate}
\end{lemma}
 {\em Proof.}
First of all we shall show that if $p'\star\vec{x}$ is a p-flag, then so is $p\star\vec{x}$. So suppose that there is $l>0$ such that $x_l\leq^+ \delta\gamma^{(l+1)}(p')$.

If $x_l\leq^+\delta\gamma^{(l+1)}(p)$, then $p\star\vec{x}$ is a p-flag, as required.

Now assume that $x_l\leq^+\delta\gamma^{(l+1)}(p)$ does not hold. Then there is $t\in\delta\gamma^{(l+1)}(p)$ and an upper $P-\gamma(P)$-path from $t$ through $x_l$, $t'\in\delta\gamma^{(l+1)}(p')$ to $\gamma^{(l)}(p)$
\[ t,a_1,\ldots , a_{s_1},x_l,a_{s_1+1},\ldots,a_{s_2}, t',a_{s_2+1},\ldots,a_{s_3}, \gamma^{(l)}(p) \]
where $0< s_1\leq s_2<s_3$,  $a_i\in  P_{l+1}-\gamma(P_{l+2})$ and $ a_i<^+\gamma^{(l+1)}(p)$, for $i=1,\ldots, s_3$. To see that such a path exists, note that  $\gamma^{(l+1)}(p')\leq^+a_{s_2+1}$ and then use that fact that $x_l\leq t'$ and Path Lemma from \cite{Z1}.

Now consider an upper $P-\gamma(P)$-path from $a_{s_1+1}$ to $\gamma^{(l+1)}(p)$
\[ a_{s_1+1}, \alpha_1,\ldots,\alpha_r,\gamma^{(l+1)}(p)\]
with $r>0$.  Since $x_l\in\delta(a_{s_1+1})$ and $a_{s_1+1}\in\delta(\alpha_1)$, we have that either  $x_l\in\iota(\alpha_1)$ or $x_l\in\delta\gamma(\alpha_1)$. In the latter case $x_l\in\delta\delta(\alpha_2)$ and again  either  $x_l\in\iota(\alpha_2)$ or $x_l\in\delta\gamma(\alpha_2)$.
But we do not have $x_l\in\delta\gamma^{(l+1)}(p)$. Thus there has to be $r_0$ such that $1\leq r_0\leq r$ and $x_l\in \iota(\alpha_{r_0})$.

Since $x_l\in\partial(x_{l+1})$, by Lemma 5.5 \cite{Z1}, we have $x_{l+1}<^+\gamma(\alpha_{r_0})$.

By Path Lemma \cite{Z1}, we have two cases.
\begin{enumerate}
  \item $\alpha_{r_0}\leq^+\gamma^{(l+2)}(p)$,
  \item $\alpha_{r_0}\leq^-\gamma^{(l+2)}(p)$.
\end{enumerate}

In case 1, by Lemma 5.3 \cite{Z1}, we have $x_l\in\iota\gamma^{(l+2)}(p)$. Hence,  by Lemma 5.5 \cite{Z1}, $x_{l+1}\leq^+\delta\gamma^{(l+2)}(p)$.

In case 2, since $x_l\in\partial(x_{l+1})$, again by Lemma 5.5 \cite{Z1}, we have $x_{l+1}<^+\gamma(\alpha_{r_0})$. But in this case $\gamma(\alpha_{r_0})\leq^+ \delta\gamma^{(l+2)}(p)$.

Thus in either case we get that $x_{l+1}\leq^+\delta\gamma^{(l+2)}(p)$ and hence $p\star\vec{x}$ is a p-flag, as required.

It remains to show that the signs of the flat faces agree. But this follows from the fact that if $p'\inn p$, then
$\delta\gamma^{(1)}(p)\leq^+\delta\gamma^{(1)}(p')$ and $\gamma^{(0)}(p')\leq^+\gamma^{(0)}(p)$.

Thus,  if  $\gamma^{(0)}(p)\leq^+x_0$, then $\gamma^{(0)}(p')\leq^+x_0$
and if $x_0\leq^+\delta\gamma^{(1)}(p)$, then $x_0\leq^+\delta\gamma^{(1)}(p')$.
$\Box$\vskip 2mm

\begin{lemma}[Iteration of $\star$]\label{lemma-iteration-of-star}
Let $\vec{x}=[x_n,\ldots, x_0]$ be a flag or p-flag in $P$ with $x_n=\bm_P$.  Let $p'\in P_k$,  $p\in P_m$ so that $p'\inn p\inn x_n$. Then

\begin{enumerate}
  \item  either $p'\neq x_k$ or $p=x_m$ and then
  \[ p'\star (p\star\vec{x}) =p'\star \vec{x}\]
  \item or $p'= x_k$ and $p\neq x_m$ and $p\star\vec{x}$ is flat, and then we have
  \[ p'\star (p\star\vec{x}) = p'\star \vec{x}= \left\{
                                     \begin{array}{ll}
                                     -x_k  & \mbox{if }  p\star\vec{x}=-p;\\
                                      & \\
                                     +x_k  & \mbox{if }  p\star\vec{x}=+p;\\
                                    \end{array}
			    \right. \]
  \item or $p'= x_k$ and $p\neq x_m$ and $p\star\vec{x}$ is a  p-flag, with $l=\pu(p\star\vec{x})$, and we have
  \begin{enumerate}
    \item $$p\star\vec{x}= [ p,\gamma(p),\ldots, \gamma^{(l+2)}(p),t,0, x_{l-1},\ldots, x_0 ]$$
    with $t\in \delta\gamma^{(l+2)}(p)$, and $t=x_{l+1}$, if $l<k$;
    \item \end{enumerate}
\[ p'\star (p\star\vec{x}) =  \left\{
                                     \begin{array}{ll}
                                     \vec{x}_{\lceil k}  & \mbox{if }  k< l;\\
                                      & \\
                                     (\vec{x}_{\lceil k})_{high} & \mbox{if } \mbox{either } k=l \mbox{ and } \gamma(t)=x_k,\\
                                     & \mbox{or } k=l+1;\\
                                       & \\
                                     (\vec{x}_{\lceil k})_{low} & \mbox{if } \mbox{either } k=l \mbox{ and } \gamma(t)\not\leq^+x_k,\\
                                     & \mbox{or } k\geq l+2.\\
                                    \end{array}
			    \right. \]
  \end{enumerate}
\end{lemma}

 {\em Proof.}
 Ad 1. If $p=x_m$, the equality
 \begin{equation}\label{equation-iter-star}
 p'\star(p\star\vec{x})=p'\star\vec{x}
 \end{equation}
 is obvious as both sides of the equation are defined with the reference to the same inequations.

 So assume that $p\neq x_m$ and $p'\neq x_k$. If $p\star \vec{x}$ is flat, then, by Lemma \ref{lemma-flat-star}, we have
 $p'\star(p\star\vec{p})=p'\star\vec{x}$.

So we can now assume that $p\star \vec{x}$ is not flat, and let it be a p-flag as displayed
  \[ p\star\vec{x} =       \begin{array}{l}
                                      \left[
                                       \begin{array}{c}
                                         p\\
                                         \ldots \\
                                         \gamma^{(l+2)}(p)\\
                                         t\\
                                         0\\
                                         x_{l-1}\\
                                        \ldots \\
                                         x_0 \\
                                       \end{array}\right]
                                    \end{array}
                                    \]
 with $x_{l+1}\leq^+t\in \delta\gamma^{(l+2)}(p)$. Since $p'\inn p$, we cannot have $\gamma^{(l')}(p)\leq^+\delta\gamma^{(l'+1)}(p')$ with $l'\geq l+2$.
 Since $x_{l+1}\leq^+t$, if $t\leq^+\delta\gamma^{(l+2)}(p')$, then we have $x_{l+1}\leq^+\delta\gamma^{(l+2)}(p')$. On the other hand, if $x_{l+1}\leq^+t'\in\delta\gamma^{(l+2)}(p')$, we have $t\perp^+t'$ and, since $t,t'\in \delta\gamma^{(l+2)}(p')$, we have $t=t'$. Thus
  \begin{equation}\label{equiv-iter-1}
 t\leq^+\delta\gamma^{(l+2)}(p')\mbox{ iff }  x_{l+1}\leq^+\delta\gamma^{(l+2)}(p')
 \end{equation}
 and if either side of the equivalence (\ref{equiv-iter-1}) holds, we have
\[   p'\star(p\star\vec{x})=   \begin{array}{l}
                                      \left[
                                       \begin{array}{c}
                                         p'\\
                                         \ldots \\
                                         \gamma^{(l+2)}(p')\\
                                         t'\\
                                         0\\
                                         x_{l-1}\\
                                        \ldots \\
                                         x_0 \\
                                       \end{array}\right]
                                    \end{array} =p'\star\vec{x} \]
 where $t'\in  \delta\gamma^{(l+2)}(p')$ is the only element in this set such that $t\leq^+ t'$.

 If neither side of (\ref{equiv-iter-1}) holds, we shall show that
   \begin{equation}\label{equiv-iter-2}
 \gamma(t)\leq^+\delta\gamma^{(l+1)}(p')\mbox{ iff }  x_{l}\leq^+\delta\gamma^{(l+1)}(p').
 \end{equation}
 The implication $\Ra$ is obvious, since $x_l\leq^+\gamma(t)$.

We shall show $\La$. Let $x_{l}\leq^+s\in\delta\gamma^{(l+1)}(p')$. As $x_{l+1}\leq^+ t$ and $x_l\in\partial(x_{l+1})$, there is $u\in \delta(t)$ such that $u\leq^+ x_l$ and
 \[ u\leq ^+x_{l}\leq^+\in\delta\gamma^{(l+1)}(p').\]

 We have that $t$ is $<^+$-minimal in $P[p]$ but not $t\not\leq^+ \delta\gamma^{(l+2)}(p')$, hence $t\not\perp^+ \delta\gamma^{(l+2)}(p')$.
 So we can apply Lemma \ref{lemma-flag-tech-J}.2, and we get
 \[ \gamma(t)\leq^+ \delta\gamma^{(l+1)}(p').\]

Thus the implication $\La$ in (\ref{equiv-iter-2})
also holds. Now, if either of the equivalent statements of (\ref{equiv-iter-2}) holds, we have
\[   p'\star(p\star\vec{x})=   \begin{array}{l}
                                      \left[
                                       \begin{array}{c}
                                         p'\\
                                         \ldots \\
                                         \gamma^{(l+1)}(p')\\
                                         t''\\
                                         0\\
                                         x_{l-2}\\
                                        \ldots \\
                                         x_0 \\
                                       \end{array}\right]
                                    \end{array} =p'\star\vec{x} \]
 where $t''\in  \delta\gamma^{(l+1)}(p')$ is the only element in this set such that $\gamma(t)\leq^+ t''$.

If neither side of equivalences (\ref{equiv-iter-1}) and (\ref{equiv-iter-2}) holds, then the formulas defining both  $p'\star(p\star\vec{x})$ and $p'\star\vec{x}$ become the same and hence we also have
\[   p'\star(p\star\vec{x})= p'\star\vec{x} \]
as required.
\vskip 2mm

 Ad 2. Obvious.
 \vskip 2mm

 Ad 3(a).
 Since $l<k$,  we have $x_{l+1}\leq^+ \delta\gamma^{(l+2)}(p)$ and $x_{l+1}\inn p$. Hence $x_{l+1}\in\delta\gamma^{(l+2)}(p)$ and then $x_{l+1}=t$.
\vskip 2mm

 Ad 3(b).
 The case $l<k$ is obvious.

 Case $l=k$. Thus we have that $x_{k+1}\leq^+ t\in \delta\gamma^{(k+2)}(p)$.

 If $\gamma(t)=x_k$, then we get immediately that $x_k\star(p\star\vec{x})=(\vec{x}_{\lc k})_{high}$.

 If $\gamma(t)\not\leq^+x_k$, then
 \[ x_k\star(p\star\vec{x}) =       \begin{array}{l}
                                      \left[
                                       \begin{array}{c}
                                         x_k\\
                                         \ldots \\
                                         x_{l'+1}\\
                                         0\\
                                         x_{l'-1}\\
                                        \ldots \\
                                         x_0 \\
                                       \end{array}\right] = (\vec{x}_{\lc k})_{low} \\
                                    \end{array}
                                    \]
where $l'=max\{ l'' |\; x_{l''+1}\leq^+\delta\gamma^{(l''+2)}(p)\}$ and $t'$ is the unique element in $\delta\gamma^{(l'+2)}(p)$ such that $x_{l'+1}\leq^+t'$.	It is easy to see that
\[ t'=x_{l'+1}\;\; \mbox{ and }\;\; \gamma^{(l'')}(x_k)=x_{l''} \]
for $l'+2\leq l''<k$. Thus the second equation above holds, as well. 		

 Case $k=l+1$. Then $p\star \vec{x}=[p,\ldots,\gamma^{(k+2)}(p),x_k,0,x_{k-2},\ldots,x_0]$ and hence
 \[ x_k\star(p\star\vec{x}) = (\vec{x}_{\lc k})_{high} \]
 as required.

 Case $k\geq l+2$. We claim that, for $l+1<k'<k$, we have
 \[ x_{k'}=\gamma(x_{k'+1})=\gamma^{(k')}(p). \]
Clearly, we have that
\[ x_{k'}\in \gamma^{(k)}(p)\cup  \iota\gamma^{(k+2)}(p)\cup \delta\gamma^{(k+1)}(p). \]
If we were to have $x_{k'}\in\delta\gamma^{(k+1)}(p)$, we would have $k\leq l+1$, contrary to the supposition.

If $ x_{k'}\in \iota\gamma^{(k+2)}(p)$, then this means that $\partial(x_{k'+1})\in\iota\gamma^{(k+2)}(p)$ and hence, by \cite{Z1}, we have $x_{k'+1}\in\delta\gamma^{(k+2)}(p)$. Hence $k'\leq l$ and we get a contradiction again. Thus we have $ x_{k'}=\gamma^{(k')}(p)$
and hence $ x_{k'}=\gamma(x_{k'+1})$, as well. Therefore

\[ x_k\star(p\star\vec{x}) =          \left[
                                       \begin{array}{c}
                                         x_k\\
                                         \ldots \\
                                         \gamma^{(l+2)}(p)\\
                                         x_{l+1}\\
                                         0\\
                                         x_{l-1}\\
                                        \ldots \\
                                         x_0 \\
                                       \end{array}\right]
                                       =
                                       \left[
                                       \begin{array}{c}
                                         x_k\\
                                         \ldots \\
                                         \gamma^{(l+2)}(p)\\
                                         x_{l+1}\\
                                         x_l\\
                                         x_{l-1}\\
                                        \ldots \\
                                         x_0 \\
                                       \end{array}\right]_{low}
                                       =  \left[
                                       \begin{array}{c}
                                         x_k\\
                                         \ldots \\
                                         x_0 \\
                                       \end{array}\right]_{low}
                                       =  (\vec{x}_{\lc k})_{low}
                                     \]
where the first equation follows from $\gamma^{(k')}(p)=\gamma^{(k')}(x_k)\not\leq^+\delta\gamma^{(k'+1)}(p)$, for $l+1<k'<k$.  $\Box$\vskip 2mm

\subsection{Dual opetope}

Let $P^{op}$ be a positive hypergraph arising from $P$ by switching the domains and codomais of faces of dimension 1, i.e. interchanging the values of functions $\gamma$ and $\delta$ on $P_1$. E. Finster made the following observation

\vskip 2mm
 \begin{lemma}
  $P^{op}$ is an opetope.
 \end{lemma}
  {\em Proof.} Simple check.  $\Box$\vskip 2mm
   We call the opetope $P^{op}$ the {\em dual of opetope} $P$. We also have

\vskip 2mm
 \begin{proposition}[flags in dual opetope]\label{prop-fags-in-dual-ope}
 \begin{enumerate}
   \item The sets of flags, maximal flags, and p-flags of the opetopes $P$ and $P^{op}$ are the same.
   \item If $\vec{x}$ is a flag of dimension $k>1$, then $\sgn^{P^{op}}(\vec{x})=-\sgn^P(\vec{x})$.
   \item The pencils $\Pl_x^P$ and $\Pl_x^{P^{op}}$  over the same face $x$ of the opetopes $P$ and $P^{op}$, respectively, are the same.
   \item If $x\in P_0$, then the linear order of pencil  $\Pl_x^P$ in $P$ is the inverse of the linear order of pencil  $\Pl_x^P$ in $P^{op}$.
   \item If $x\in P_{m}$ and $m>0$, then the linear order of pencil  $\Pl_x^P$ in $P$ is the same as the linear order of pencil  $\Pl_x^P$ in $P^{op}$.
   \item\label{DualMaxFlagOrder} The linear order of flags in  $\Flags^P[k,0)$  is inverse of the linear order of flags in $\Flags^{P^{op}}[k,0)$.
   \item The linear order of flags in  $\pFlags^P[x,0)$  is inverse of the linear order of flags in $\Flags^{P^{op}}[x,0)$.
   \item \label{Star-on-dual} Let $p\in P_m$ and $\vec{x}$ be a maximal flag in $P$. Then
     \[ (p\star^{P^{op}}\vec{x}) =  \left\{
                                     \begin{array}{ll}
                                     +p  & \mbox{if }  p\star^P\vec{x}=-p,\\
                                      & \\
                                     -p  & \mbox{if }  p\star^P\vec{x}=+p,\\
                                     & \\
                                     p\star^P\vec{x} & \mbox{otherwise.}
                                    \end{array}
			    \right. \]
 \end{enumerate}
 \end{proposition}

 \vskip 2mm
 {\em Remark. }
 Note that if we also invert the operation $-$ and $+$ on $P^{op}$, i.e. we consider $-^P$, $+^P$ as an operation on faces of $P$ and  $-^{P^{op}}$ and $+^{P^{op}}$  as an operation on faces of $P^{op}$ so that $-^{P^{op}}p=  +^Pp$ and $+^{P^{op}}p=  -^Pp$, for $p\in P$, then the above formula could be simplified to
 \[ (p\star^{P^{op}}\vec{x}) =  p\star^P\vec{x}.  \]

 \vskip 2mm

  {\em Proof.} All the above observations, except (\ref{DualMaxFlagOrder}), should be easy.

  To see that (\ref{DualMaxFlagOrder}) also holds, note that if we have two flags $\vec{x}$ and $\vec{y}$ in $\Flags^P[k,0]$, with $k>0$,
  then at the level 1, with $x_0=y_0$, the order of faces in pencil $x_0$ in  $\Pl_{x_0}^P$ is inverse of the order $\Pl_{y_0}^{P^{op}}$ and at the level $l>1$ the orders $\Pl_{x_{l-1}}^P$ and  $\Pl_{y_{l-1}}^{P^{op}}$ are the same but the signs $\vec{x}_{\lc l-1}$ and  $\vec{y}_{\lc l-1}$ are opposite. So at any dimension the conditions for being `greater than' on the set $\Flags^P[k,0]=\Flags^{P^{op}}[k,0]$ in $P$ and $P^{op}$ are opposite.
  $\Box$\vskip 2mm

\subsection{Cylinder over a positive opetope as a positive hypergraph}
In this section we equip $\Cyl(P)$ with operations $\gamma$ and $\delta$ so that it becomes a hypergraph $\Cylh(P)$. Eventually, we shall show that $\Cylh(P)$ is $\pOpe$-straight and its image $\Cylp(P)$  in $\widehat{\pOpe}$ is aspherical.

Let $\vec{x}=[x_k,\ldots,x_0]$ be a flag in $\Cyl(P)$. We put
 \[ \gamma(\vec{x})= \vec{x}_{(k-1)}=\vec{x}_{high} \]
and
\[ \delta(\vec{x}) =\{ \vec{x}_{(k)}, \vec{x}_{low} \}=\{ \vec{x}_{side}, \vec{x}_{low} \}\]
in other words
\[ \delta(\vec{x}) =  \left\{
                                    \begin{array}{ll}
                                      \{ \vec{x}_{(k)}, -x_k \} & \mbox{if } \vec{x} \mbox{ is the first $[x_k,0)$-flag}, \\
                                      \{ \vec{x}_{(k)}, +x_k \} & \mbox{if } \vec{x} \mbox{ is the last $[x_k,0)$-flag}, \\
                                      \{ \vec{x}_{(k)}, \vec{x}_{(\ll(\vec{x}))} \} & \mbox{otherwise.}

                                    \end{array}
			    \right. \]

Let $\vec{x}=[x_k,\ldots, \widehat{x_i},\ldots, x_0]$ be a p-falg in $P$. We put

\[ \gamma(\vec{x}) =\gamma(x_k)\star\vec{x}=  \left\{
                                    \begin{array}{ll}
                                      \left[
                                       \begin{array}{c}
                                         \gamma(x_k)\\
                                         0 \\
                                         x_{n-3}\\
                                        \ldots \\
                                         x_0 \\
                                       \end{array}\right] & \mbox{if } i=k-1,k-2, \\
                                       & \\
                                      \left[
                                       \begin{array}{c}
                                          x_{k-1}  \\
                                        x_{k-2} \\
                                        \ldots\\
                                        \widehat{x_{i}} \\
                                        \ldots\\
                                        x_0 \\
                                       \end{array}\right] & \mbox{otherwise.}
                                    \end{array}
			    \right. \]

We put
\[ \delta(\vec{x}) =  \left\{
                                    \begin{array}{ll}
                                      \{ p\star  \vec{x} \; | \; p\in \delta(x_k)  \} & \mbox{if }\vec{x}\mbox{ is a low p-flag,} \\
                                       & \\
                                     \{ p\star  \vec{x} \; | \; p\in \delta(x_k)  \}\cup \{ \vec{x}_{(k)}\} & \mbox{if }\vec{x}\mbox{ is a high p-flag.}
                                    \end{array}
			    \right. \]

For $p\in P_{\geq 1}$, we put
\[ \gamma(-p)=-\gamma(p),\hskip 5mm  \gamma(+p)=+\gamma(p) \]
and
\[ \delta(-p)=\{ -q : q\in \delta(p)\} ,\hskip 5mm  \delta(+p)=\{ +q : q\in \delta(p)\}. \]

\subsection{$P^{\vec{x}}$ as a positive opetope}\label{Faces-in-opetope}

Let $\vec{x}$ be a maximal flag in $P$. By $P^{\vec{x}}$ we denote the least sub-hypergraph of $\Cylh(P)$ containing $\vec{x}$ (and closed under $\gamma$ and $\delta$) together with functions $\gamma$ and $\delta$ restricted to this subset.

 {\bf Faces of flag opetope} $P^{\vec{x}}$.

Let $\vec{x}=[x_n,\ldots,x_0]$. The dimension of $P^{\vec{x}}$ is $n+1$.
Let us describe the faces of $P^{\vec{x}}_l$, for $0\leq l\leq n+1$:

  \[ P^{\vec{x}}_l =  \left\{
                                    \begin{array}{ll}
                                     \{ \vec{x} \} & \mbox{if }l=n+1;\\
                                       & \\
                                     \{ p\star  \vec{x} \; | \; p\in P_l-\{ x_l\}  \}\cup \{ (\vec{x}_{\lceil l})_{low},\; (\vec{x}_{\lceil l})_{high},\; \vec{x}_{\lceil l-1}\}  & \mbox{if }0<l\leq n;\\
                                     & \\
                                      \{ -p \; | \; p\in P_0,\; p\leq^+x_0  \}\cup  \{ +p \; | \; p\in P_0,\; x_0\leq^+p  \} & \mbox{if }l=0. \\
                                    \end{array}
			    \right. \]

It should be clear that all these faces must occur in a subset of $\Cylh(P)$ containing $\vec{x}$ and closed under $\gamma$ and $\delta$.
We shall show that the above set is indeed closed under $\gamma$ and $\delta$.

For the faces of the form $\vec{x}_{\lc l}$ it is clear.

For the face $(\vec{x}_{\lc l})_{high}$, we have
\[ \gamma((\vec{x}_{\lc l})_{high})= \gamma(x_l)\star (\vec{x}_{\lc l})_{high} = \gamma(x_l)\star \vec{x}\in P^{\vec{x}}_{l-1}\]
Moreover
\[   \delta_{side}((\vec{x}_{\lc l})_{high}) = (\vec{x}_{\lc l-2})\in P^{\vec{x}}_{l-1} \]
and, if $x_{l-1}\neq p \in \delta(x_l)$, then
\[ p\star (\vec{x}_{\lc l})_{high} = p\star \vec{x} \in P^{\vec{x}}_{l-1}  \]
if $x_{l-1}= p \in \delta(x_l)$, then
 \[ p\star (\vec{x}_{\lc l})_{high} = (\vec{x}_{\lc l-1})_{low} \in P^{\vec{x}}_{l-1}.  \]
Thus the codomain and all domains of $(\vec{x}_{\lc l})_{high}$ are in $P^{\vec{x}}_{l-1}$.

For the face $(\vec{x}_{\lc l})_{low}$, we have
\[ \gamma((\vec{x}_{\lc l})_{low})= \gamma(x_l)\star (\vec{x}_{\lc l})_{low} = \;= \; \left\{ \begin{array}{ll}
		  (\vec{x}_{\lc l})_{high} & \mbox{ if } x_{l-1}=\gamma(x_l) \\
  &\\
		  \gamma(x_l)\star \vec{x}   & \mbox{ if } x_{l-1}\in\delta(x_l).
                                    \end{array}
			    \right. \]
Moreover, for $p\in \delta(x_l)$, we have
\[ p\star (\vec{x}_{\lc l})_{low}\;= \; \left\{ \begin{array}{ll}
           p\star \vec{x} & p\neq x_{l-1}\\
           &\\
		  (\vec{x}_{\lc l-1})_{high} & \mbox{ if } \pu((\vec{x}_{\lc l})_{low}) =l-2 \\
  &\\
		  (\vec{x}_{\lc l-1})_{low}  & \mbox{ if } \pu((\vec{x}_{\lc l})_{low}) <l-2.
                                    \end{array}
			    \right. \]
Thus the codomain and all domains of $(\vec{x}_{\lc l})_{low}$ are also in $P^{\vec{x}}_{l-1}$.

Finally, for a face of the form $p\star \vec{x}$, the fact that its codomain and domains belong to the set $P^{\vec{x}}$ follows from Lemma \ref{lemma-iteration-of-star}.

\vskip 2mm
{\bf Faces of p-flag opetope}  $P^{\vec{x}}$.

Let $\vec{x}=[x_n,\ldots,\widehat{x_k},\ldots,x_0]$, $0\leq k<n$. The dimension of $P^{\vec{x}}$ is $n>1$. Let us describe the faces of dimension for $0\leq l\leq n$:
  \[ P^{\vec{x}}_l =  \left\{
                                    \begin{array}{ll}
                                    \{ \vec{x} \} & \mbox{if }l=n;\\

                                       & \\
                                        \{ p\star  \vec{x} \; | \; p\in P_l  \} \cup\{ \vec{x}_{\lceil l-1}\;| l=k>0\}  & \mbox{if }k\leq l< n;\\
                                       & \\
                                     \{ p\star  \vec{x} \; | \; p\in P_l-\{ x_l\}  \}\cup \{ (\vec{x}_{\lceil l})_{low},\; (\vec{x}_{\lceil l})_{high}\} \cup\{ \vec{x}_{\lceil l-1}\;| l>0\}  & \mbox{if }0\leq l< k.\\
                                    \end{array}
			    \right. \]

The justification that this set of faces is also closed under $\gamma$ and $\delta$ is left for the reader.

Having defined $\gamma$, $\delta$ on $P^{\vec{x}}$, we have a notion of occurrence in  $P^{\vec{x}}$ and hence we can link it to the operation $\star$. We have

\begin{proposition}\label{prop-unique-proj}
 Let $\vec{x}$ be a maximal flag, $p\in P_k$, $p\not\in \vec{x}$. Then $p\star \vec{x}\in P^{\vec{x}}_k$ and it is the unique face in $P^{\vec{x}}$ that projects on $p$.
\end{proposition}
{\em Proof.} This is an easy consequence of Lemma \ref{lemma-iteration-of-star}.  $\Box$\vskip 2mm

\begin{proposition}\label{prop-opetope-in-prod} Let $P$ be a positive opetope, $\vec{x}$ a maximal  flag in $\Cylh(P)$. Then $P^{\vec{x}}$ with $\gamma$ and $\delta$ described above forms a positive opetope.
\end{proposition}

{\em Proof.} The proof is by checking axioms, one by one. We shall split it into Lemmas \ref{lemma-opetope-in-prod-globularity}, \ref{lemma-opetope-in-prod-strict-disjoint}, \ref{lemma-opetope-in-prod-pencil-lin}  below. $\Box$\vskip 2mm

\vskip 2mm
\begin{lemma}\label{lemma-opetope-in-prod-globularity} Let $P$ be a positive opetope, $\vec{x}$ a  maximal flag in $\Cylh(P)$. Then $P^{\vec{x}}$ with $\gamma$ and $\delta$ described satisfies the globularity axiom.
\end{lemma}

{\em Proof.}
Let us fix a maximal flag $\vec{x}$  in $P$. We show that $P^{\vec{x}}$ satisfies the globularity axiom.

 We need to check the globularity condition of flags, p-flags, and flat faces in $P^{\vec{x}}$.
The globularity conditions for flat faces are directly inherited from the globularity conditions in $P$. We need to check globularity conditions
\begin{equation}\label{equation-glob-gamma}
\gamma\gamma(\vec{x})=\gamma\delta(\vec{x})-\delta\delta(\vec{x})
\end{equation}
\begin{equation}\label{equation-glob-delta}
\delta\gamma(\vec{x})=\delta\delta(\vec{x})-\gamma\delta(\vec{x})
\end{equation}
for flags and p-flags $\vec{x}$ with the top face $x_n=\bm_P$. For other faces in $P$ we can restrict the whole opetope $P$ to make these faces top faces in a smaller opetope. For $n=2$, we can check the equations directly. Thus assume that $n>2$.

\vskip 1mm
{\em Globularity for flags.} We shall make direct verification of the above equations. We consider four cases.
\begin{enumerate}
  \item $\vec{x}$ is the first flag,
  \item $\vec{x}$ is the last flag,
  \item $x_{n-1}=\gamma(x_n)$ (but not the above),
  \item $x_{n-1}\in\delta(x_n)$ (but not the above).
\end{enumerate}
In each case we shall spell what are the four sets of faces involved and then make some comments concerning the equations.

Ad 1. We have $x_i=\gamma^{(i)}(x_n)$ for $i=0,\ldots, n-1$. Then

\[ \gamma\gamma(\vec{x}) =\gamma    \left[
                                       \begin{array}{c}
                                         x_n\\
                                         0 \\
                                         x_{n-2}\\
                                        \vdots \\
                                         x_0 \\
                                       \end{array}\right]
                                      =   \left[
                                       \begin{array}{c}
                                         x_{n-1}\\
                                         0 \\
                                         x_{n-3}\\
                                        \vdots \\
                                         x_0 \\
                                       \end{array}\right]
                                        \]

\[ \delta\gamma(\vec{x}) =\delta   \left[
                                       \begin{array}{c}
                                         x_n\\
                                         0 \\
                                         x_{n-2}\\
                                        \vdots \\
                                         x_0 \\
                                       \end{array}\right]
                                      =\{  \left[
                                       \begin{array}{c}
                                         x_{n-2}\\
                                        \vdots \\
                                         x_0 \\
                                       \end{array}\right], -p : p\in \delta(x_n)\}
                                        \]
Note that $x_0$ is $<^+$-maximal in this case.
\[ \gamma\delta(\vec{x}) =\gamma \{  \left[
                                       \begin{array}{c}
                                         x_{n-1}\\
                                        \vdots \\
                                         x_0 \\
                                       \end{array}\right],-x_n \}
                                      =\{  \left[
                                       \begin{array}{c}
                                         x_{n-1}\\
                                         0\\
                                         x_{n-3}\\
                                        \vdots \\
                                         x_0 \\
                                       \end{array}\right], -x_{n-1} \}
                                        \]

\[ \delta\delta(\vec{x}) =\delta \{  \left[
                                       \begin{array}{c}
                                         x_{n-1}\\
                                        \vdots \\
                                         x_0 \\
                                       \end{array}\right],-x_n \}
                                      =\{  \left[
                                       \begin{array}{c}
                                         x_{n-2}\\
                                        \vdots \\
                                         x_0 \\
                                       \end{array}\right], -x_{n-1}, -p:p\in\delta(x_n) \}
                                        \]
Clearly the equations (\ref{equation-glob-gamma}) and (\ref{equation-glob-delta}) hold in this case.

Ad 2. We have $x_i=\gamma^{(i)}(x_n)$ for $i=1,\ldots, n-1$ and $x_0=\delta(x_1)$. Then

\[ \gamma\gamma(\vec{x}) =\gamma    \left[
                                       \begin{array}{c}
                                         x_n\\
                                         0 \\
                                         x_{n-2}\\
                                        \vdots \\
                                         x_0 \\
                                       \end{array}\right]
                                      =   \left[
                                       \begin{array}{c}
                                         x_{n-1}\\
                                         0 \\
                                         x_{n-3}\\
                                        \vdots \\
                                         x_0 \\
                                       \end{array}\right]
                                        \]

\[ \delta\gamma(\vec{x}) =\delta   \left[
                                       \begin{array}{c}
                                         x_n\\
                                         0 \\
                                         x_{n-2}\\
                                        \vdots \\
                                         x_0 \\
                                       \end{array}\right]
                                      =\{  \left[
                                       \begin{array}{c}
                                         x_{n-2}\\
                                        \vdots \\
                                         x_0 \\
                                       \end{array}\right], +p : p\in \delta(x_n)\}
                                        \]
Note that $x_0$ is $<^+$-minimal in this case.
\[ \gamma\delta(\vec{x}) =\gamma \{  \left[
                                       \begin{array}{c}
                                         x_{n-1}\\
                                        \vdots \\
                                         x_0 \\
                                       \end{array}\right],+x_n \}
                                      =\{  \left[
                                       \begin{array}{c}
                                         x_{n-1}\\
                                         0\\
                                         x_{n-3}\\
                                        \vdots \\
                                         x_0 \\
                                       \end{array}\right], +x_{n-1} \}
                                        \]

\[ \delta\delta(\vec{x}) =\delta \{  \left[
                                       \begin{array}{c}
                                         x_{n-1}\\
                                        \vdots \\
                                         x_0 \\
                                       \end{array}\right],+x_n \}
                                      =\{  \left[
                                       \begin{array}{c}
                                         x_{n-2}\\
                                        \vdots \\
                                         x_0 \\
                                       \end{array}\right], +x_{n-1}, +p:p\in\delta(x_n) \}
                                        \]
Clearly the equations (\ref{equation-glob-gamma}) and (\ref{equation-glob-delta}) hold in this case.

Ad 3. We have $x_{n-1}=\gamma(x_n)$  and $\vec{x}$ is not the first flag. Then $\gamma(\vec{x}_{low})\in\delta(\vec{x}_{side})$ and

\[ \gamma\gamma(\vec{x}) =\gamma    \left[
                                       \begin{array}{c}
                                         x_n\\
                                         0 \\
                                         x_{n-2}\\
                                        \vdots \\
                                         x_0 \\
                                       \end{array}\right]
                                      =   \left[
                                       \begin{array}{c}
                                         x_{n-1}\\
                                         0 \\
                                         x_{n-3}\\
                                        \vdots \\
                                         x_0 \\
                                       \end{array}\right]
                                        \]

\[ \delta\gamma(\vec{x}) =\delta   \left[
                                       \begin{array}{c}
                                         x_n\\
                                         0 \\
                                         x_{n-2}\\
                                        \vdots \\
                                         x_0 \\
                                       \end{array}\right]
                                      =\{  \left[
                                       \begin{array}{c}
                                         x_{n-2}\\
                                        \vdots \\
                                         x_0 \\
                                       \end{array}\right], p\star\left[
                                       \begin{array}{c}
                                         x_{n}\\
                                         0\\
                                         x_{n-2}\\
                                        \vdots \\
                                         x_0 \\
                                       \end{array}\right] : p\in \delta(x_n)\}
                                        \]

\[ \gamma\delta(\vec{x}) =\gamma \{  \left[
                                       \begin{array}{c}
                                         x_{n-1}\\
                                        \vdots \\
                                         x_0 \\
                                       \end{array}\right], \vec{x}_{low}\}
                                      =\{  \left[
                                       \begin{array}{c}
                                         x_{n-1}\\
                                         0\\
                                         x_{n-3}\\
                                        \vdots \\
                                         x_0 \\
                                       \end{array}\right], \gamma(x_n)\star\vec{x}_{low} \}
                                       =\{  \left[
                                       \begin{array}{c}
                                         x_{n-1}\\
                                         0\\
                                         x_{n-3}\\
                                        \vdots \\
                                         x_0 \\
                                       \end{array}\right],
                                       \left[
                                       \begin{array}{c}
                                         x_{n-1}\\
                                        \vdots \\
                                         x_0 \\
                                       \end{array}\right]_{low} \}
                                        \]
For the last equality recall that $\gamma(x_n)=x_{n-1}$ in this case.

\[ \delta\delta(\vec{x}) =\delta \{  \left[
                                       \begin{array}{c}
                                         x_{n-1}\\
                                        \vdots \\
                                         x_0 \\
                                       \end{array}\right],\vec{x}_{low} \}
                                      =\{  \left[
                                       \begin{array}{c}
                                         x_{n-2}\\
                                        \vdots \\
                                         x_0 \\
                                       \end{array}\right],
                                        \left[
                                       \begin{array}{c}
                                         x_{n-1}\\
                                        \vdots \\
                                         x_0 \\
                                       \end{array}\right]_{low},
                                       p\star  \left[
                                       \begin{array}{c}
                                         x_{n}\\
                                        \vdots \\
                                         x_0 \\
                                       \end{array}\right]_{low}
                                       :p\in\delta(x_n) \}
                                        \]
The equations (\ref{equation-glob-gamma}) and (\ref{equation-glob-delta}) hold in this case.

Ad 4. We have $x_{n-1}\in\delta(x_n)$. Then $\gamma(\vec{x}_{side})\in\delta(\vec{x}_{low})$ and

\[ \gamma\gamma(\vec{x}) =\gamma    \left[
                                       \begin{array}{c}
                                         x_n\\
                                         0 \\
                                         x_{n-2}\\
                                        \vdots \\
                                         x_0 \\
                                       \end{array}\right]
                                      =   \gamma(x_n)\star \left[
                                       \begin{array}{c}
                                         x_n\\
                                         0 \\
                                         x_{n-3}\\
                                        \vdots \\
                                         x_0 \\
                                       \end{array}\right]
                                       =   \left[
                                       \begin{array}{c}
                                         \gamma(x_n)\\
                                         0 \\
                                         x_{n-3}\\
                                        \vdots \\
                                         x_0 \\
                                       \end{array}\right]
                                        \]

\[ \delta\gamma(\vec{x}) =\delta   \left[
                                       \begin{array}{c}
                                         x_n\\
                                         0 \\
                                         x_{n-2}\\
                                        \vdots \\
                                         x_0 \\
                                       \end{array}\right]
                                      =\{  \left[
                                       \begin{array}{c}
                                         x_{n-2}\\
                                        \vdots \\
                                         x_0 \\
                                       \end{array}\right], p\star\left[
                                       \begin{array}{c}
                                         x_{n}\\
                                         0\\
                                         x_{n-2}\\
                                        \vdots \\
                                         x_0 \\
                                       \end{array}\right] : p\in \delta(x_n)\}
                                        \]

\[ \gamma\delta(\vec{x}) =\gamma \{  \left[
                                       \begin{array}{c}
                                         x_{n-1}\\
                                        \vdots \\
                                         x_0 \\
                                       \end{array}\right],
                                        \left[
                                       \begin{array}{c}
                                         x_n\\
                                         x_{n-1}\\
                                         0\\
                                         x_{n-3}\\
                                        \vdots \\
                                         x_0 \\
                                       \end{array}\right]\}
                                      =\{  \left[
                                       \begin{array}{c}
                                         x_{n-1}\\
                                         0\\
                                         x_{n-3}\\
                                        \vdots \\
                                         x_0 \\
                                       \end{array}\right],
                                       \gamma(x_n)\star
                                       \left[
                                       \begin{array}{c}
                                         x_n\\
                                         x_{n-1}\\
                                         0\\
                                         x_{n-3}\\
                                        \vdots \\
                                         x_0 \\
                                       \end{array}\right] \}
                                       =\{  \left[
                                       \begin{array}{c}
                                         x_{n-1}\\
                                         0\\
                                         x_{n-3}\\
                                        \vdots \\
                                         x_0 \\
                                       \end{array}\right],
                                       \left[
                                       \begin{array}{c}
                                         \gamma(x_n)\\
                                         0\\
                                         x_{n-3}\\
                                        \vdots \\
                                         x_0 \\
                                       \end{array}\right] \}
                                        \]

\[ \delta\delta(\vec{x}) =\delta \{  \left[
                                       \begin{array}{c}
                                         x_{n-1}\\
                                        \vdots \\
                                         x_0 \\
                                       \end{array}\right],\vec{x}_{low} \}
                                      =\{  \left[
                                       \begin{array}{c}
                                         x_{n-2}\\
                                        \vdots \\
                                         x_0 \\
                                       \end{array}\right],
                                        \left[
                                       \begin{array}{c}
                                         x_{n-1}\\
                                        \vdots \\
                                         x_0 \\
                                       \end{array}\right]_{low},
                                       p\star  \left[
                                       \begin{array}{c}
                                         x_{n}\\
                                         x_{n-1}\\
                                         0\\
                                         x_{n-3}\\
                                        \vdots \\
                                         x_0 \\
                                       \end{array}\right]
                                       :p\in\delta(x_n) \}
                                        \]
To see that the equations (\ref{equation-glob-gamma}) and (\ref{equation-glob-delta}) hold in this case as well, note that
\[ x_{n-1}\star \left[
                                       \begin{array}{c}
                                         x_{n}\\
                                         x_{n-1}\\
                                         0\\
                                         x_{n-3}\\
                                        \vdots \\
                                         x_0 \\
                                       \end{array}\right] = \left[
                                       \begin{array}{c}
                                         x_{n-1}\\
                                         0\\
                                         x_{n-3}\\
                                        \vdots \\
                                         x_0 \\
                                       \end{array}\right].\]
and
\[ x_{n-1}\star \left[
                                       \begin{array}{c}
                                         x_{n}\\
                                         0\\
                                         x_{n-2}\\
                                        \vdots \\
                                         x_0 \\
                                       \end{array}\right] = \left[
                                       \begin{array}{c}
                                         x_{n-1}\\
                                        \vdots \\
                                         x_0 \\
                                       \end{array}\right]_{low}.\]

\vskip 1mm
{\em Globularity for p-flags.} Let $\vec{x}$ be a p-flag.
We consider three cases
\begin{enumerate}
  \item $\pu(\vec{x})=l<n-2$,
  \item $\pu(\vec{x})=l=n-2$,
  \item $\pu(\vec{x})=l=n-1$.
\end{enumerate}

Ad 1. We have
\[ \gamma\gamma(\vec{x}) =  \gamma\left[
                                       \begin{array}{c}
                                         x_{n-1}\\
                                         \vdots\\
                                         x_{l+1}\\
                                         0\\
                                         x_{l-1}\\
                                        \vdots \\
                                         x_0 \\
                                       \end{array}\right]
                                       =  \gamma^2(x_n)\star
                                       \left[
                                       \begin{array}{c}
                                         x_{n-1}\\
                                         \vdots\\
                                         x_{l+1}\\
                                         0\\
                                         x_{l-1}\\
                                        \vdots \\
                                         x_0 \\
                                       \end{array}\right]
                                        \]
Note that $\gamma(x_n)=x_{n-1}$ in this case.

\[ \delta\gamma(\vec{x}) =\delta   \left[
                                       \begin{array}{c}
                                         x_{n-1}\\
                                         \vdots\\
                                         x_{l+1}\\
                                         0\\
                                         x_{l-1}\\
                                        \vdots \\
                                         x_0 \\
                                       \end{array}\right]
                                      =\{  p\star
                                      \left[
                                       \begin{array}{c}
                                         x_{n-1}\\
                                         \vdots\\
                                         x_{l+1}\\
                                         0\\
                                         x_{l-1}\\
                                        \vdots \\
                                         x_0 \\
                                       \end{array}\right]
                                        : p\in \delta(x_{n-1})\}
                                        \]

\[ \gamma\delta(\vec{x}) =\gamma \{  p\star \vec{x} : p\in \delta(x_n)\}
                                      =\{ \gamma(p)\star \vec{x} : p\in \delta(x_n)  \}   \]

\[ \delta\delta(\vec{x}) =\delta \{  p\star \vec{x} : p\in \delta(x_n) \} =\]
\[  =\{    p'\star(p\star \vec{x}) :p'\in \delta(p),\; p\in \delta(x_n) \} = \{    p'\star\vec{x} :p'\in \delta\delta(x_n) \} \]
Note that $x_{n-1}=\gamma(x_n)\not\in\delta(x_n)$.
The equations (\ref{equation-glob-gamma}) and (\ref{equation-glob-delta}) hold in this case.

Ad 2. Now assume that $\pu(\vec{x})=n-2$. We have
\[ \gamma\gamma(\vec{x}) =  \left[
                                       \begin{array}{c}
                                         \gamma^2(x_n)\\
                                         0 \\
                                         x_{n-4}\\
                                        \vdots \\
                                         x_0 \\
                                       \end{array}\right]
                                      =   \gamma^2(x_n)\star \vec{x}
                                       =  \gamma^2(x_n)\star
                                       \left[
                                       \begin{array}{c}
                                         \gamma(x_n)\\
                                         0 \\
                                         x_{n-3}\\
                                        \vdots \\
                                         x_0 \\
                                       \end{array}\right]
                                        \]

\[ \delta\gamma(\vec{x}) =\delta   \left[
                                       \begin{array}{c}
                                         \gamma(x_n)\\
                                         0 \\
                                         x_{n-3}\\
                                        \vdots \\
                                         x_0 \\
                                       \end{array}\right]
                                      =\{  \left[
                                       \begin{array}{c}
                                         x_{n-3}\\
                                        \vdots \\
                                         x_0 \\
                                       \end{array}\right], p\star\vec{x} : p\in \delta\gamma(x_n)\}
                                        \]

The remaining sets are as in the previous case

\[ \gamma\delta(\vec{x}) =\{ \gamma(p)\star \vec{x} : p\in \delta(x_n),  \}   \]

\[ \delta\delta(\vec{x})  = \{    p'\star\vec{x} :p'\in \delta\delta(x_n). \} \]
The equations (\ref{equation-glob-gamma}) and (\ref{equation-glob-delta}) hold in this case.

Ad 3. Finally assume that $\pu(\vec{x})=n-1$, i.e. $\vec{x}$ is a high p-flag. We have
\[ \gamma\gamma(\vec{x}) =  \left[
                                       \begin{array}{c}
                                         \gamma^2(x_n)\\
                                         0 \\
                                         x_{n-4}\\
                                        \vdots \\
                                         x_0 \\
                                       \end{array}\right]
                                        \]

\[ \delta\gamma(\vec{x}) =\delta   \left[
                                       \begin{array}{c}
                                         \gamma(x_n)\\
                                         0 \\
                                         x_{n-3}\\
                                        \vdots \\
                                         x_0 \\
                                       \end{array}\right]
                                      =\{  \left[
                                       \begin{array}{c}
                                         x_{n-3}\\
                                        \vdots \\
                                         x_0 \\
                                       \end{array}\right],p\star\left[
                                       \begin{array}{c}
                                         \gamma(x_n)\\
                                         0 \\
                                         x_{n-3}\\
                                        \vdots \\
                                         x_0 \\
                                       \end{array}\right] : p\in \delta\gamma(x_n)\}
                                        \]

\[ \gamma\delta(\vec{x}) =\gamma\{ (\vec{x}_{\lc n-2}), \gamma(p)\star \vec{x} : p\in \delta(x_n)  \} =    \]

\[ = \{ (\vec{x}_{\lc n-2})_{high}, \gamma(p)\star (p\star\vec{x}) : p\in \delta(x_n)  \} \]

\[ \delta\delta(\vec{x})  = \delta\{ (\vec{x}_{\lc n-2}),
                                          p'\star\vec{x} :p'\in \delta\delta(x_n) \} =
                                          \]

\[ = \{ (\vec{x}_{\lc n-2})_{low},
                                       (\vec{x}_{\lc n-3}),
                                       p'\star(p\star\vec{x}) :p'\in\delta(p),\; p\in\delta(x_n) \} =\]

\[ = \{ (\vec{x}_{\lc n-2})_{low},
                                       (\vec{x}_{\lc n-3}), p'\star\vec{x} :p'\in\delta\delta(x_n)-\{ x_{n-2}\} \} \]
In this case an argument is needed. To see that the equation (\ref{equation-glob-gamma}) holds, we need to show that the face $(\vec{x}_{\lc n-2})_{high}$ does not present a problem. To this end, note that

\begin{enumerate}
  \item either $x_{n-2}=\gamma\gamma(x_n)$ and
  \[ \gamma\gamma(\vec{x})= (\vec{x}_{\lc n-2})_{high} \in\gamma\delta(\vec{x})\]
  \item or there is  $p\in\delta(x_n)$ such that $x_{n-2}\in\delta(p)$, and then
  \[\gamma\delta(\vec{x})\ni (\vec{x}_{\lc n-2})_{high}=  x_{n-2}\star(p\star\vec{x}) \in \delta\delta(\vec{x}).\]
\end{enumerate}

To see that the equation (\ref{equation-glob-delta}) holds, we need to show that the face $(\vec{x}_{\lc n-2})_{low}$ does not present a problem.
Note that
\begin{enumerate}
  \item either there is $p\in\delta(x_n)$ such that $\gamma(p)=x_{n-2}$ and
  \[  \gamma\delta(\vec{x})\ni x_{n-2}\star\vec{x}=(\vec{x}_{\lc n-2})_{low}\in\delta\delta(\vec{x}), \]
  \item or $x_{n-2}\in\delta\gamma(x_n)$ and
  \[ \delta\delta(\vec{x})\ni (\vec{x}_{\lc n-2})_{low}=x_{n-2}\star(\gamma(x_n)\star\vec{x})\in\delta\gamma(\vec{x}).\]
\end{enumerate}

Thus the equations (\ref{equation-glob-gamma}) and (\ref{equation-glob-delta}) hold in this case, as well.
$\Box$\vskip 2mm

\vskip 2mm
\begin{lemma}\label{lemma-opetope-in-prod-strict-disjoint} Let $P$ be a positive opetope, $\vec{x}$ a maximal  flag in $\Cylh(P)$. Then $P^{\vec{x}}$ with $\gamma$ and $\delta$ described satisfies the strictness  and disjointness axioms.
\end{lemma}

{\em Proof.}
Let us fix $\vec{x}$ a maximal flag in $P$. We need to verify that $P^{\vec{x}}$ satisfies the strictness  and disjointness axioms.

We have to verify that $<^+$ on $P^{\vec{x}}_0$ is a linear order and, for $l=1,\ldots,n+1$, the orders $<^+$ and $<^-$ on $P^{\vec{x}}_l$ are strict and disjoint.

We have
\[ P^{\vec{x}}_0 = \{ -p \; | \; p\in P_0,\; p\leq^+x_0  \}\cup  \{ +p \; | \; p\in P_0,\; x_0\leq^+p  \} \]
with each summand linearly ordered. Since $\delta([x_0])=-x_0$ and $\gamma([x_0])=+x_0$, we also have $-x_0<^+ +x_0$. Thus $<^+$ on $P^{\vec{x}}_0$ is a strict linear order.

In dimension $l=n+1$ there is nothing to prove, since $P^{\vec{x}}_{n+1}$ is a singleton.

Let $0<l\leq n$. We have
\[ P^{\vec{x}}= \{ p\star  \vec{x} \; | \; p\in P_l-\{ x_l\}  \}\cup \{ (\vec{x}_{\lceil l})_{low},\; (\vec{x}_{\lceil l})_{high}, \vec{x}_{\lceil l-1}\}.\]
and hence $\star$-operation on
\[ \star : P_l-\{ x_l\} \lra P^{\vec{x}}- \{ (\vec{x}_{\lceil l})_{low},\; (\vec{x}_{\lceil l})_{high}, \vec{x}_{\lceil l-1}\} \]
is an order isomorphism with respect to both $<^+$ and $<^-$ order relations.  Thus it remains to show that if we `substitute' in place of face $x_l$ three faces $\{ (\vec{x}_{\lceil l})_{low},\; (\vec{x}_{\lceil l})_{high}, \vec{x}_{\lceil l-1}\}$, we still get strict orders that are disjoint.

The new upper order  $<^+$ looks as follows. Between new faces we have
\[   (\vec{x}_{\lceil l})_{low} <^+ (\vec{x}_{\lceil l})_{high},\hskip 10mm \vec{x}_{\lceil l-1}<^+ (\vec{x}_{\lceil l})_{high}\]
and, for $p\in P-\{ x_l\}$, we have
\[ p\star \vec{x}<^+ (\vec{x}_{\lceil l})_{low} \hskip 5mm{\rm iff} \hskip 5mm p<^+x_l\]
and
\[(\vec{x}_{\lceil l})_{high} <^+   p\star \vec{x} \hskip 5mm{\rm iff} \hskip 5mm x_l<^+p.\]
The other relations follow from these by transitivity. In particular, $\vec{x}_{\lceil l-1}$ is $<^+$-minimal face in $P^{\vec{x}}_l$. Thus it should be clear that if $<^+$ were strict on $P_l$, it is also strict on $P^{\vec{x}}_l$.

For the new lower order $<^-$, we consider two cases depending whether $\vec{x}_{\lceil l}$ is a $\gamma$-flag, i.e. $\gamma(x_l)=x_{l-1}$ or it is a $\delta$-flag, i.e. $x_{l-1}\in\delta(x_l)$. There is just one relation between new faces

   \[  \left\{ \begin{array}{lll}
        (\vec{x}_{\lceil l})_{low} <^- \vec{x}_{\lceil l-1}& \mbox{ if } \gamma(x_l)=x_{l-1} \\
        &&\\
		\vec{x}_{\lceil l-1}<^- (\vec{x}_{\lceil l})_{low} & \mbox{ otherwise.}
                                    \end{array}
			    \right. \]
The relation $<^-$ between new and old faces is as follows. In either case, for $p\in P_l-\{ x_l\}$, we have
\[ p\star\vec{x}<^- (\vec{x}_{\lceil l})_{high}  \hskip 5mm{\rm iff} \hskip 5mm p<^-x_l \]
and
\[ (\vec{x}_{\lceil l})_{high}<^- p\star\vec{x}  \hskip 5mm{\rm iff} \hskip 5mm x_l<^-p. \]

In case $\vec{x}_{\lceil l}$ is a $\gamma$-flag, we have additionally, for $p\in P_l-\{ x_l\}$,

\[ p\star\vec{x}<^- \vec{x}_{\lceil l-1}  \hskip 5mm{\rm iff} \hskip 5mm \gamma(p)\leq^+x_{l-1} \hskip 2mm{\rm  and } \hskip 2mm x_l\not\leq^+ p\]

\[ \vec{x}_{\lceil l-1}<^- p\star\vec{x}  \hskip 5mm{\rm iff} \hskip 5mm \vec{x}_{low}<^- p\star\vec{x}  \hskip 5mm{\rm iff} \hskip 5mm x_l<^-p \]

\[ p\star\vec{x}<^- \vec{x}_{low}  \hskip 5mm{\rm iff} \hskip 5mm p<^-x_l. \]

In case $\vec{x}_{\lceil l}$ is a $\delta$-flag, we have additionally, for $p\in P_l-\{ x_l\}$,

\[ \vec{x}_{\lceil l-1}<^- p\star\vec{x}  \hskip 5mm{\rm iff} \hskip 5mm \hskip 5mm x_{l-1}\leq^+\delta(p)\hskip 2mm{\rm  and } \hskip 2mm x_l\not\leq^+ p, \]

\[ \vec{x}_{low}<^- p\star\vec{x}  \hskip 5mm{\rm iff} \hskip 5mm x_l<^- p. \]

\[  p\star\vec{x}<^- \vec{x}_{low}  \hskip 5mm{\rm iff} \hskip 5mm p<^-x_l, \]

\[ p\star\vec{x}<^- \vec{x}_{\lceil l-1}   \hskip 5mm{\rm iff} \hskip 5mm \gamma(p)<^+x_{l-1}. \]
From this description it should be clear that the relation $<^-$  on $P^{\vec{x}}$ is a strict order.
$\Box$\vskip 2mm

Before we prove that $P^{\vec{x}}$ satisfies pencil linearity, we need the following technical lemma.
\vskip 2mm
\begin{lemma}\label{lemma-opetope-in-prod-technical} Let $P$ be a positive opetope, $\vec{x}$ a maximal flag in $\Cylh(P)$, $a,\alpha\in P$ and $a,\alpha\not\in \vec{x}$, $k<dim(P)$. Then
\begin{enumerate}
  \item if $a\in\delta(\alpha)$, then $a\star \vec{x}\in\delta(\alpha\star \vec{x})$;
  \item if $a=\gamma(\alpha)$, then $a\star \vec{x}=\gamma(\alpha\star \vec{x})$;
  \item if $x_k\in\partial(a)$ and $a\prec_{x_k}x_{k+1}$, then $(\vec{x}_{\lc k})_{high}\in\partial(a\star \vec{x})$
  (moreover, we can replace $\partial$ by either $\gamma$ or $\delta$);
  \item if $x_k\in\partial(a)$ and $x_{k+1}\prec_{x_k} a$, then $(\vec{x}_{\lc k})_{low}\in\partial(a\star \vec{x})$
  (moreover, we can replace $\partial$ by either $\gamma$ or $\delta$);
  \item if $x_k=\gamma(x_{k+1})$, then
  $$\gamma((\vec{x}_{\lc k+1})_{low})=(\vec{x}_{\lc k})_{low}\in\delta(\vec{x}_{\lc k}),\hskip 10mm  \gamma((\vec{x}_{\lc k+1})_{high})=(\vec{x}_{\lc k})_{high} ;$$

  \item if $x_k\in\delta(x_{k+1})$, then
    $$\gamma(\vec{x}_{\lc k})=(\vec{x}_{\lc k})_{high} \in\delta((\vec{x}_{\lc k+1})_{low}),\hskip 10mm  (\vec{x}_{\lc k})_{low}\in \delta((\vec{x}_{\lc k+1})_{high}) .$$
\end{enumerate}

\end{lemma}

{\em Proof.} 1. and 2. follow immediately from Lemma \ref{lemma-iteration-of-star} and definition of $\gamma$ and $\delta$ in $\Cylh(P)$.
5. and 6. easily follow from the definitions of $\gamma$ and $\delta$.

Ad 3. Assume that $a\prec_{x_k}x_{k+1}$. If $x_k\in\delta(a)$, then $a<^+x_{k+1}$ and we have
\[ x_k\star(a\star\vec{x})=x_k\star[a,(\vec{x}_{\lc k})_{high}] =(\vec{x}_{\lc k})_{high}.\]
Thus $(\vec{x}_{\lc k})_{high}\in\delta(a\star\vec{x})$.

If $x_k=\gamma(a)$, then $x_{k+1}<+a$ and then
\[ \gamma(a\star \vec{x})=\gamma([a,0,\vec{x}_{\lc k-1}])=(\vec{x}_{\lc k})_{high} \]
where the third equality follows from the fact that $\gamma(a)\leq^+x_k$. Thus 3. holds.

Ad 4. Now assume that $x_{k+1}\prec_{x_k}a$. If $x_k\in\delta(a)$, then $x_{k+1}<^+a$ and  we have
\[ x_k\star(a\star \vec{x})=x_k\star [a,0,\vec{x}_{\lc k-1}]= (\vec{x}_{\lc k})_{low}. \]
Thus $ (\vec{x}_{\lc k})_{low}\in \delta(a\star \vec{x})$.

Finally, if $x_k=\gamma(a)$, then
\[ \gamma(a\star \vec{x}) =\gamma([a,(\vec{x}_{\lc k})_{low}])= (\vec{x}_{\lc k})_{low},\]
as required.
 $\Box$\vskip 2mm

\vskip 2mm
\begin{lemma}\label{lemma-opetope-in-prod-pencil-lin} Let $P$ be a positive opetope, $\vec{x}$ a maximal  flag in $\Cylh(P)$. Then $P^{\vec{x}}$ with $\gamma$ and $\delta$ described satisfies pencil linearity axiom.
\end{lemma}

{\em Proof.}
Let us fix a maximal flag  $\vec{x}$ in $P$. We need to verify that $P^{\vec{x}}$ satisfies the pencil linearity axiom.

We shall describe pencils over different kinds of faces in $P^{\vec{x}}$.

{\em Case 1}: pencils over a flag $\vec{x}_{\lc k}$. Since there is no face $a$ in $\Cylh$ with $\gamma(a)$ being a flag, $\gamma$-pencil over $\vec{x}_{\lc k}$ is empty. The $\delta$-pencil over $\vec{x}_{\lc k}$ contains two flags $(\vec{x}_{\lc k})_{high}$ and $(\vec{x}_{\lc k})_{side}=\vec{x}_{\lc k-1}$. Since $(\vec{x}_{\lc k})_{side}\in \epsilon(\vec{x}_{\lc k})$ and  $\gamma(\vec{x}_{\lc k}=(\vec{x}_{\lc k})_{high}$, they are linearly ordered by $<^+$.

{\em Case 2}: pencils over p-flags that project to $x_k\in\vec{x}$. We have two p-flags projecting to $x_k$, i.e. $(\vec{x}_{\lc k})_{low}$ and $(\vec{x}_{\lc k})_{high}$. We shall describe pencils over them together in terms of the pencil over $x_k$ in $P$. The pencils over faces $(\vec{x}_{\lc k})_{low}$ and $(\vec{x}_{\lc k})_{high}$ are formed from lifts of faces in pencils over $x_k$ different from $x_{k+1}$ and additionally the faces in $\partial(\vec{x}_{\lc k+1})$.

We have two cases:
\begin{enumerate}
  \item $x_k=\gamma(x_{k+1})$,
  \item $x_k\in \delta(x_{k+1})$.
\end{enumerate}

Ad 1. In the former case (the fragment of) the pencil over $x_k$ in $P$ is as follows
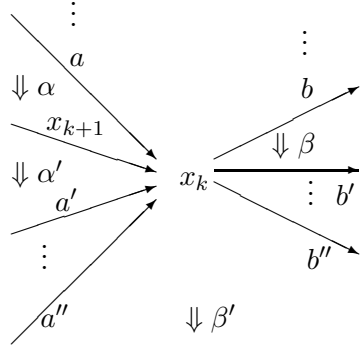
\begin{center}
\xext=2000 \yext=1200
\begin{picture}(\xext,\yext)(\xoff,\yoff)
   \put(1000,1000){$\vdots$}
   \put(200,1100){$\vdots$}

   \put(580,550){$x_k$}

\put(0,1140){\vector(1,-1){500}}
\put(0,760){\vector(3,-1){500}}
\put(0,380){\vector(3,1){500}}
\put(0,0){\vector(1,1){500}}

\put(700,640){\vector(2,1){500}}
\put(700,600){\vector(1,0){500}}
\put(700,560){\vector(2,-1){500}}

\put(200,960){$a$}
\put(0,860){$\Da\alpha$}
\put(120,740){$x_{k+1}$}
\put(0,560){$\Da\alpha'$}
\put(150,460){$a'$}
\put(100,260){$\vdots$}
\put(100,60){$a''$}
\put(600,60){$\Da\beta'$}

\put(1000,850){$b$}
\put(900,660){$\Da\beta$}
\put(1120,480){$b'$}
\put(1020,480){$\vdots$}
\put(1020,260){$b''$}

\end{picture}
\end{center}
with  $\alpha$'s and $\beta$'s in $P-\gamma(P)$ lifts to two pencils
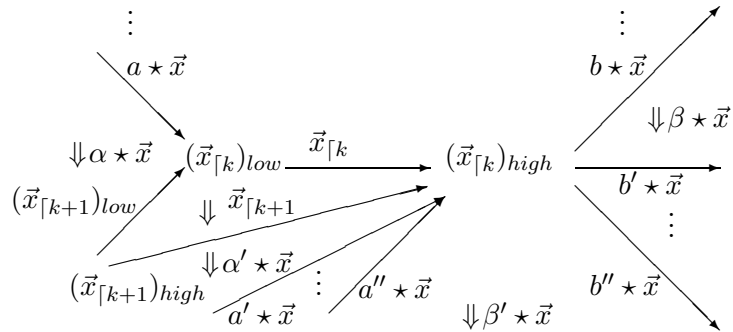
\begin{center}
\xext=2000 \yext=1300
\begin{picture}(\xext,\yext)(\xoff,\yoff)


   \put(300,550){$(\vec{x}_{\lc k})_{low}$}
   \put(730,600){$\vec{x}_{\lc k}$}
   \put(1200,550){$(\vec{x}_{\lc k})_{high}$}
   \put(-100,100){$(\vec{x}_{\lc k+1})_{high}$}
   \put(-300,400){$(\vec{x}_{\lc k+1})_{low}$}
   \put(450,400){$\vec{x}_{\lc k+1}$}
   \put(750,100){$\vdots$}
   \put(1800,1000){$\vdots$}
   \put(100,1000){$\vdots$}

\put(650,540){\vector(1,0){500}}
\put(0,940){\vector(1,-1){300}}
\put(0,240){\vector(1,1){300}}
\put(40,200){\vector(4,1){1100}}
\put(400,40){\vector(2,1){800}}
\put(800,40){\vector(1,1){400}}

\put(-100,560){$\Da$}
\put(-30,560){$\alpha\star \vec{x}$}
\put(100,860){$a\star \vec{x}$}
\put(450,0){$a'\star \vec{x}$}
\put(900,90){$a''\star \vec{x}$}
\put(350,360){$\Da$}
\put(360,180){$\Da$}
\put(430,180){$\alpha'\star \vec{x}$}
\put(1260,0){$\Da$}
\put(1330,0){$\beta'\star \vec{x}$}

\put(1650,540){\vector(1,0){500}}
\put(1650,600){\vector(1,1){500}}

\put(1650,480){\vector(1,-1){500}}

\put(1700,860){$b\star \vec{x}$}
\put(1900,680){$\Da$}
\put(1970,680){$\beta\star \vec{x}$}
\put(1800,450){$b'\star \vec{x}$}
\put(1970,300){$\vdots$}
\put(1700,100){$b''\star \vec{x}$}
\end{picture}
\end{center}
one over $(\vec{x}_{\lc k})_{low}$ and the other over $(\vec{x}_{\lc k})_{high}$  in $P^{\vec{x}}$. This immediately follows from Lemma \ref{lemma-opetope-in-prod-technical}.
If $\alpha$ ($\alpha'$) would be equal $x_{k+2}$, then in the above picture $\alpha\star\vec{x}$ ($\alpha'\star\vec{x}$) would need to be replaced by $(\vec{x}_{k+2})_{low}$.

Ad 2. In the later case (the fragment of) the pencil over $x_k$ in $P$, as follows

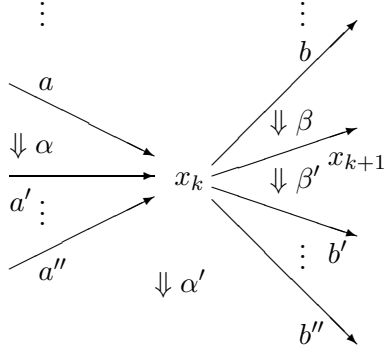
\begin{center}
\xext=2000 \yext=1200
\begin{picture}(\xext,\yext)(\xoff,\yoff)
   \put(1000,1100){$\vdots$}
   \put(100,1100){$\vdots$}

   \put(570,550){$x_k$}
\put(700,620){\vector(1,1){500}}
\put(700,580){\vector(3,1){500}}
\put(700,540){\vector(3,-1){500}}
\put(700,500){\vector(1,-1){500}}

\put(0,260){\vector(2,1){500}}
\put(0,580){\vector(1,0){500}}
\put(0,900){\vector(2,-1){500}}

\put(100,880){$a$}
\put(0,660){$\Da\alpha$}
\put(0,460){$a'$}
\put(100,410){$\vdots$}
\put(100,220){$a''$}
\put(500,180){$\Da\alpha'$}
\put(1000,980){$b$}
\put(900,740){$\Da\beta$}
\put(1100,620){$x_{k+1}$}
\put(1100,280){$b'$}
\put(900,540){$\Da\beta'$}
\put(1000,0){$b''$}

\put(1000,260){$\vdots$}
\end{picture}
\end{center}
with  $\alpha$'s and $\beta$'s in $P-\gamma(P)$ lifts to two pencils
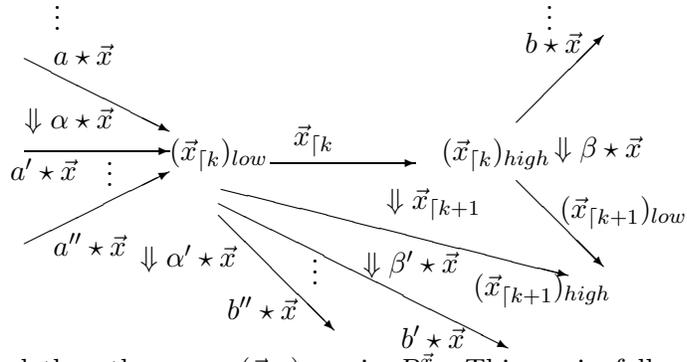
\begin{center}
\xext=2000 \yext=1300
\begin{picture}(\xext,\yext)(\xoff,\yoff)

\put(0,260){\vector(2,1){500}}
\put(0,580){\vector(1,0){500}}
\put(0,900){\vector(2,-1){500}}

\put(100,880){$a\star\vec{x}$}
\put(0,660){$\Da\alpha\star\vec{x}$}
\put(-50,480){$a'\star\vec{x}$}
\put(280,460){$\vdots$}
\put(100,220){$a''\star\vec{x}$}
\put(400,180){$\Da\alpha'\star\vec{x}$}

   \put(500,550){$(\vec{x}_{\lc k})_{low}$}
   \put(930,600){$\vec{x}_{\lc k}$}
   \put(1440,550){$(\vec{x}_{\lc k})_{high}$}

  \put(1700,680){\vector(1,1){300}}
   \put(1700,480){\vector(1,-1){300}}
  \put(680,450){\vector(4,-1){1200}}

  \put(670,400){\vector(2,-1){1000}}
  \put(670,360){\vector(1,-1){400}}
  \put(1730,920){$b\star \vec{x}$}
  \put(1830,560){$\Da\beta\star \vec{x}$}

   \put(1300,-100){$b'\star \vec{x}$}
  \put(1170,160){$\Da\beta'\star \vec{x}$}
  \put(990,120){$\vdots$}
   \put(700,0){$b''\star \vec{x}$}

   \put(1550,80){$(\vec{x}_{\lc k+1})_{high}$}
   \put(1850,350){$(\vec{x}_{\lc k+1})_{low}$}
   \put(1250,380){$\Da\vec{x}_{\lc k+1}$}
  \put(1800,1000){$\vdots$}
   \put(100,1000){$\vdots$}

\put(850,540){\vector(1,0){500}}
\end{picture}
\end{center}
one over $(\vec{x}_{\lc k})_{low}$ and the other over $(\vec{x}_{\lc k})_{high}$  in $P^{\vec{x}}$. This again follows from Lemma \ref{lemma-opetope-in-prod-technical}. Thus in this case it is clear that the pencils are linearly ordered. Similarly as in the previous case, if $\beta$ ($\beta'$) would be equal $x_{k+2}$, then in the above picture $\beta\star\vec{x}$ ($\beta'\star\vec{x}$) would need to be replaced by $(\vec{x}_{k+2})_{low}$.

 {\em Case 3}: pencil over a face $z\star \vec{x}$ in $P^{\vec{x}}$ so that $z\not\in\vec{x}$ and $z\not\in\partial(x_{k+1})$.

In this case the faces in the pencil over $z$ in $P$ lift uniquely to the faces over $z\star\vec{x}\in P^{\vec{x}}$. In particular, we have isomorphism of pencils $(-)\star \vec{x} : (\Pl_z,<^+)\lra (\Pl_{z\star \vec{x}},<^+)$ and hence both $\gamma$- and $\delta$-pencil over $z\star \vec{x}$ are linearly ordered by $<^+$.

 {\em Case 4}: pencil over a face $z\star \vec{x}$ in $P^{\vec{x}}$ so that $z\not\in\vec{x}$.

 As in 2, we consider two cases
 \begin{enumerate}
   \item $z=\gamma(x_{k+1})$;
   \item $z\in \delta(x_{k+1})$.
 \end{enumerate}

Ad 1. In the former case (the fragment of) the pencil over $z$ in $P$ is as follows

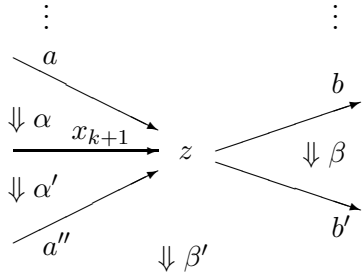
\begin{center}
\xext=1200 \yext=1050
\begin{picture}(\xext,\yext)(\xoff,\yoff)
   \put(1100,1000){$\vdots$}
   \put(100,1000){$\vdots$}

   \put(570,550){$z$}
\put(700,580){\vector(3,1){500}}
\put(700,540){\vector(3,-1){500}}

\put(0,260){\vector(2,1){500}}
\put(0,580){\vector(1,0){500}}
\put(0,900){\vector(2,-1){500}}

\put(100,880){$a$}
\put(-20,660){$\Da\alpha$}
\put(-20,420){$\Da\alpha'$}
\put(200,620){$x_{k+1}$}
\put(100,220){$a''$}
\put(500,180){$\Da\beta'$}
\put(1100,780){$b$}
\put(1100,280){$b'$}
\put(1000,540){$\Da\beta$}
\end{picture}
\end{center}
with  $\alpha$'s and $\beta$'s in $P-\gamma(P)$ lifts to a pencil
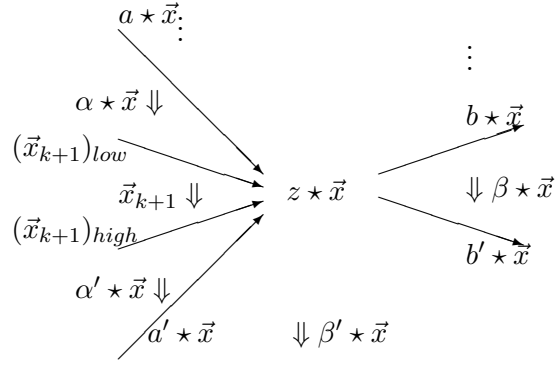
\begin{center}
\xext=2000 \yext=1200
\begin{picture}(\xext,\yext)(\xoff,\yoff)
   \put(1200,1000){$\vdots$}
   \put(200,1100){$\vdots$}

\put(580,550){$z\star \vec{x}$}

\put(0,1140){\vector(1,-1){500}}
\put(0,760){\vector(3,-1){500}}
\put(0,380){\vector(3,1){500}}
\put(0,0){\vector(1,1){500}}

\put(900,640){\vector(3,1){500}}
\put(900,560){\vector(3,-1){500}}

\put(0,1160){$a\star\vec{x}$}
\put(-150,860){$\alpha\star\vec{x}\Da$}
\put(-370,700){$(\vec{x}_{k+1})_{low}$}
\put(-370,420){$(\vec{x}_{k+1})_{high}$}
\put(0,540){$\vec{x}_{k+1}\Da$}
\put(-150,210){$\alpha'\star\vec{x}\Da$}
\put(100,60){$a'\star\vec{x}$}
\put(600,60){$\Da\beta'\star\vec{x}$}

\put(1200,810){$b\star\vec{x}$}
\put(1200,560){$\Da\beta\star\vec{x}$}
\put(1200,330){$b'\star\vec{x}$}
\end{picture}
\end{center}
over $z\star\vec{x}$  in $P^{\vec{x}}$. This again follows from Lemma \ref{lemma-opetope-in-prod-technical}. Thus in this case it is clear that the pencil is linearly ordered. Similarly as in the previous case, if $\alpha$ ($\alpha'$) would be equal $x_{k+2}$, then in the above picture $\alpha\star\vec{x}$ ($\alpha'\star\vec{x}$) would need to be replaced by $(\vec{x}_{k+2})_{low}$.

Ad 2. In the latter case (the fragment of) the pencil over $z$ in $P$ is as follows

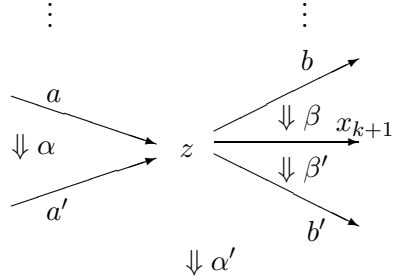
\begin{center}
\xext=2000 \yext=1200
\begin{picture}(\xext,\yext)(\xoff,\yoff)
   \put(1000,1000){$\vdots$}
   \put(120,1000){$\vdots$}

   \put(580,550){$z$}
\put(0,760){\vector(3,-1){500}}
\put(0,380){\vector(3,1){500}}

\put(700,640){\vector(2,1){500}}
\put(700,600){\vector(1,0){500}}
\put(700,560){\vector(2,-1){500}}

\put(120,740){$a$}
\put(0,560){$\Da\alpha$}
\put(120,330){$a'$}
\put(600,160){$\Da\alpha'$}

\put(1000,850){$b$}
\put(920,660){$\Da\beta$}
\put(920,480){$\Da\beta'$}
\put(1120,640){$x_{k+1}$}

\put(1020,260){$b'$}

\end{picture}
\end{center}
with  $\alpha$'s and $\beta$'s in $P-\gamma(P)$ lifts to a pencil
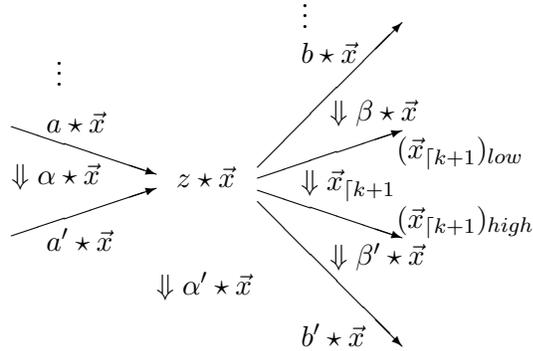
\begin{center}
\xext=2000 \yext=1200
\begin{picture}(\xext,\yext)(\xoff,\yoff)
   \put(1000,1100){$\vdots$}
   \put(150,900){$\vdots$}

\put(120,740){$a\star\vec{x}$}
\put(0,560){$\Da\alpha\star\vec{x}$}
\put(120,330){$a'\star\vec{x}$}
\put(0,760){\vector(3,-1){500}}
\put(0,380){\vector(3,1){500}}

   \put(570,550){$z\star\vec{x}$}
\put(850,620){\vector(1,1){500}}
\put(850,580){\vector(3,1){500}}
\put(850,540){\vector(3,-1){500}}
\put(850,500){\vector(1,-1){500}}


\put(500,180){$\Da\alpha'\star\vec{x}$}
\put(1000,980){$b\star\vec{x}$}
\put(1100,780){$\Da\beta\star\vec{x}$}
\put(1330,650){$(\vec{x}_{\lc k+1})_{low}$}
\put(1000,540){$\Da\vec{x}_{\lc k+1}$}
\put(1330,410){$(\vec{x}_{\lc k+1})_{high}$}
\put(1100,280){$\Da\beta'\star\vec{x}$}
\put(1000,0){$b'\star\vec{x}$}

\end{picture}
\end{center}
over $z\star\vec{x}$  in $P^{\vec{x}}$. This again follows from Lemma \ref{lemma-opetope-in-prod-technical}. Thus in this case it is clear that the pencils are linearly ordered. Similarly as in the previous case, if $\beta$ ($\beta'$) would be equal $x_{k+2}$, then in the above picture $\beta\star\vec{x}$ ($\beta'\star\vec{x}$) would need to be replaced by $(\vec{x}_{k+2})_{low}$.
$\Box$\vskip 2mm

\subsection{Intersections of flags}
In this section we shall describe intersections of flags. We need the following Proposition.

  \begin{proposition}\label{prop-monotone}
  Let $P$ be an opetope of dimension $n$, $k<n$, $p\in P_k$. Then the functions
  \[ p\star(-)_l,  p\star(-)_h : \Flags_P[\bm_P,0)\lra \pFlags_P[p,0) \cup \{ -p,+p\}\]
  are monotone, where
 \[ p\star_l\vec{x} =  \left\{
                                    \begin{array}{ll}
                                     p\star  \vec{x} & \mbox{if } p\neq x_k,\\
                                       & \\
                                     (p\star  \vec{x})_{low} & \mbox{otherwise,}\\
                                    \end{array}
			    \right. \]
and
  \[ p\star_h\vec{x} =  \left\{
                                    \begin{array}{ll}
                                     p\star  \vec{x} & \mbox{if } p\neq x_k,\\
                                       & \\
                                     (p\star  \vec{x})_{high} & \mbox{otherwise}\\
                                    \end{array}
			    \right. \]
for $\vec{x}\in \Flags_P[\bm_P,0)$.

Moreover, if $\vec{x}$ and $\vec{x}'$ are two consecutive flags in $\Flags_P[\bm_P,0)$, then
\[ p\star_l\vec{x}\unlhd  p\star(\vec{x}\cap \vec{x}')\unlhd p\star_l\vec{x}'\]
and
\[ p\star_h\vec{x}\unlhd  p\star(\vec{x}\cap \vec{x}')\unlhd p\star_h\vec{x}'.\]
\end{proposition}

{\em Proof.} Before we start proving the main part of the above Proposition, we explain how we can modify the proof to get the part starting from `Moreover,'. The whole proof is done for just consecutive flags. Differing at level $k$ only. So we need to know that for any face $q$ we have $x_k\leq^+q$ and $x'_k\leq^+q$  iff   $\gamma(x_{k+1})\leq^+q$. This is true, since either one of $x_k$, $x'_k$ is just equal $\gamma(x_{k+1})$ or $\gamma(x_{k+1})$ is the $<^+$ successor of both $x_k$ and $x'_k$.

Now we turn to the proof of the main part of Proposition \ref{prop-monotone}. Let $\vec{x}\unlhd\vec{x}'$, $p\in P_m$. We shall show that
  \begin{equation}\label{equation-star-l}
  p\star_l\vec{x}\unlhd p\star_l\vec{x}'
   \end{equation}
and
   \begin{equation}\label{equation-star-h}
  p\star_h\vec{x}\unlhd p\star_h\vec{x}'.
   \end{equation}

It is enough to assume that $\vec{x}'$ is a successor of $\vec{x}$. Let $\diff(\vec{x},\vec{x}')=k$. Thus we have

\[ \vec{x} =          \left[
                                       \begin{array}{c}
                                         x_n\\
                                         \ldots \\
                                         \gamma^{(l+2)}(x_n)\\
                                         x_{k+1}\\
                                         x_k\\
                                         x_{k-1}\\
                                        \ldots \\
                                         x_0 \\
                                       \end{array}\right]
                                       \lhds
                                       \left[
                                       \begin{array}{c}
                                         x_n\\
                                         \ldots \\
                                         \gamma^{(l+2)}(x_n)\\
                                         x_{k+1}\\
                                         x'_k\\
                                         x_{k-1}\\
                                        \ldots \\
                                         x_0 \\
                                       \end{array}\right]
                                       = \vec{x}'
                                     \]
We shall consider cases concerning relation of $p$ to $x_i$'s and in case $p=x_k$ we shall consider all possible types of successor $\vec{x}'$.

{\em Case 1.} $p\not\in \vec{x}\cup \vec{x}'$.

Let $l=max\{l' \,|\, x_{l'+1}\leq^+\delta \gamma^{(l'+1)}(p)\}$. Then $l\leq k$. Since $x_k$, $x_k'$ are two consecutive faces in pencil $\Pl_{x_{k-1}}$, we have for any $q\in P_{k}-\{ x_k,x_k' \}$ that
\[ x_k\leq^+ q\;\;{\rm iff }\;\; x_k'\leq^+ q. \]
In particular
\[ x_k\leq^+ \delta \gamma^{(k+1)}(p)\;\;{\rm iff }\;\; x_k'\leq^+ \delta \gamma^{(k+1)}(p). \]
and hence $p\star\vec{x}= p\star\vec{x}'$. Since  $p\star\vec{x}$ is a p-flag, we also have that both equations (\ref{equation-star-l}) and (\ref{equation-star-h}) hold.

 {\em Case 2.} $p=x_m$ and $m<k$.

 In this case we have
 \[ p\star_l\vec{x}=(x_m\star\vec{x})_{low} = (\vec{x}_{\lc m})_{low}= \]
 \[ (\vec{x}'_{\lc m})_{low}=(x_m\star\vec{x}')_{low} =p\star_l\vec{x}', \]
 i.e.  equation (\ref{equation-star-l}) holds. Similar calculations show that equation (\ref{equation-star-h}) holds, as well.

  {\em Case 3.} $p=x_m$ and $m>k$.

  We have
  \[ p\star \vec{x} = \vec{x}_{\lc m} \lhds  \vec{x}'_{\lc m} = p\star \vec{x}' \]
  and since the functions $(-)_{low}$ and $(-)_{high}$  are monotone, we get again that both equations (\ref{equation-star-l}) and (\ref{equation-star-h}) hold.

  {\em Case 4.} $p=x_k$ and $m=k$.

  This is the main case where we proceed by inspection of all possible types of successor $\vec{x}'$. Note that we do not need to consider the case $p=x_k'$ separately since it will follow by duality, i.e. it is the case 4 for the dual opetope $P^{op}$.

  There are eight cases to consider.

  \begin{enumerate}
    \item $k<n-1$.
    \begin{enumerate}
      \item $\sgn(\vec{x}_{\lc k-1})=+1$.
      \begin{enumerate}
        \item $\vec{x}'$ is  $\delta$- or $\delta\gamma$-successor of $\vec{x}$. We have
        \[ x_k\star_l\vec{x}=(\vec{x}_{\lc k})_{low}=x_k\star_l\vec{x}'\]
        and
        \[  x_k\star_h\vec{x}=(\vec{x}_{\lc k})_{high}\lhds (\vec{x}_{\lc k})_{low} =x_k\star\vec{x}'=x_k\star_h\vec{x}'.\]
        The last equality follows from the fact that $x_k\star\vec{x}'$ is already a p-flag. This observation we will use often without notice, as in the calculations right below.
        \item $\vec{x}'$ is $\gamma$-successor of $\vec{x}$. We have
         \[  x_k\star_l\vec{x}=(\vec{x}_{\lc k})_{low}\lhds (\vec{x}_{\lc k})_{high} =x_k\star\vec{x}'=x_k\star_l\vec{x}'\]
         and
         \[ x_k\star_h\vec{x}=(\vec{x}_{\lc k})_{high}=x_k\star\vec{x}'=x_k\star_h\vec{x}'.\]
      \end{enumerate}
      \item $\sgn(\vec{x}_{\lc k-1})=-1$.
           \begin{enumerate}
        \item $\vec{x}'$ is inverse $\gamma$- or $\gamma\delta$-successor of $\vec{x}$. We have
        \[ x_k\star_l\vec{x}= (\vec{x}_{\lc k})_{low}=x_k\star\vec{x}'=x_k\star_l\vec{x}'\]
        and
         \[ x_k\star_h\vec{x}= (\vec{x}_{\lc k})_{high}\lhds(\vec{x}_{\lc k})_{low}= x_k\star\vec{x}'=x_k\star_h\vec{x}'.\]
        \item $\vec{x}'$ is inverse $\delta$-successor of $\vec{x}$. We have
        \[ x_k\star_l\vec{x}= (\vec{x}_{\lc k})_{low}\lhds(\vec{x}_{\lc k})_{high}= x_k\star\vec{x}'=x_k\star_l\vec{x}'\]
        and
        \[ x_k\star_h\vec{x}= (\vec{x}_{\lc k})_{high}=x_k\star\vec{x}'=x_k\star_h\vec{x}'.\]
      \end{enumerate}
    \end{enumerate}
    \item $k=n-1$.
        \begin{enumerate}
      \item $\sgn(\vec{x}_{\lc k-1})=+1$.
      \begin{enumerate}
        \item $\vec{x}'$ is $\delta$- or $\delta\gamma$-successor of $\vec{x}$. We have
        \[ x_{n-1}\star_l\vec{x}= (\vec{x}_{\lc n-1})_{low}= x_{n-1}\star\vec{x}'=x_{n-1}\star_l\vec{x}'\]
        and
        \[ x_{n-1}\star_h\vec{x}= (\vec{x}_{\lc n-1})_{high}\lhds (\vec{x}_{\lc n-1})_{low}= x_{n-1}\star\vec{x}'=x_{n-1}\star_h\vec{x}'.\]
        \item $\vec{x}'$ is $\gamma$-successor of $\vec{x}$. We have
         \[ x_{n-1}\star_l\vec{x}= (\vec{x}_{\lc n-1})_{low}\lhds (\vec{x}_{\lc n-1})_{high}= x_{n-1}\star\vec{x}'=x_{n-1}\star_h\vec{x}'\]
         and
         \[ x_{n-1}\star_h\vec{x}= (\vec{x}_{\lc n-1})_{high}= x_{n-1}\star\vec{x}'=x_{n-1}\star_h\vec{x}'.\]
      \end{enumerate}
      \item $\sgn(\vec{x}_{\lc k-1})=-1$.
           \begin{enumerate}
           \item $\vec{x}'$ is inverse $\delta$-successor of $\vec{x}$. We have
             \[ x_{n-1}\star_l\vec{x}= (\vec{x}_{\lc n-1})_{low}\lhds (\vec{x}_{\lc n-1})_{high}= x_{n-1}\star\vec{x}'=x_{n-1}\star_l\vec{x}'\]
        and
        \[ x_{n-1}\star_h\vec{x}= (\vec{x}_{\lc n-1})_{high}= x_{n-1}\star\vec{x}'=x_{n-1}\star_h\vec{x}'.\]
          \item $\vec{x}'$ is  inverse $\gamma$- or $\gamma\delta$-successor. We have
          \[ x_{n-1}\star_l\vec{x}= (\vec{x}_{\lc n-1})_{low}= x_{n-1}\star\vec{x}'=x_{n-1}\star_l\vec{x}'\]
          and
          \[ x_{n-1}\star_h\vec{x}= (\vec{x}_{\lc n-1})_{high}\lhds (\vec{x}_{\lc n-1})_{low}= x_{n-1}\star\vec{x}'=x_{n-1}\star_h\vec{x}'.\]
      \end{enumerate}
    \end{enumerate}
  \end{enumerate}
$\Box$\vskip 2mm

\begin{proposition}\label{lemma-flag-intersection}
  Let $P$ be a positive opetope, $\vec{y}, \vec{x}, \vec{x}'$ maximal flags in $P$ such that $\vec{y}\lhd \vec{x}$ and $\next_{\lhd}(\vec{x})=\vec{x}'$. Then
\begin{enumerate}
  \item $P^{\vec{x}}\cap P^{\vec{x}'}=P^{\vec{x}\cap\vec{x}'}$;
  \item $P^{\vec{y}}\cap P^{\vec{x}'} \subseteq P^{\vec{x}}\cap P^{\vec{x}'} $.
  \end{enumerate}
  \end{proposition}
  {\em Proof.} Ad 1. Let $\vec{x}$, $\vec{x}'$ be maximal flags in $P$ with $\vec{x}'$ being the successor of $\vec{x}$ and $k=\diff(\vec{x},\vec{x}')$.
  Since ${\vec{x}\cap\vec{x}'}$ is a face of both ${\vec{x}}$ and ${\vec{x}'}$, we have   $P^{\vec{x}\cap\vec{x}'}\subseteq P^{\vec{x}}\cap P^{\vec{x}'}$.

  Let $p\inn x_n$. Suppose $p\not\in\vec{x}\cup\vec{x}'$ and
  \[ p\star\vec{x}=         \left[
                                       \begin{array}{c}
                                         p\\
                                         \ldots \\
                                         \gamma^{(l+2)}(p)\\
                                         t\\
                                         0\\
                                         x_{l-1}\\
                                        \ldots \\
                                         x_0 \\
                                       \end{array}\right]
                                       =p\star\vec{x}'
                                     \]
   Then $l\leq k$. So $x_i\not\leq^+\delta\gamma^{(i+1)}(p)$, for $i>l+1$. Moreover, by Lemma \ref{lemma-flag-tech-J}, the conditions $x_k\leq^+ \delta\gamma^{(k+1)}(p)$ and $x'_k\leq^+ \delta\gamma^{(k+1)}(p)$ are equivalent to $\gamma(x_{k+1})\leq^+ \delta\gamma^{(k+1)}(p)$.
   Thus $p\star\vec{x}=p\star(\vec{x}\cap\vec{x}')$ and hence $p\star\vec{x}\in P^{\vec{x}\cap\vec{x}'}$. Thus, for faces of the form $p\star\vec{x}$ with $p\not\in\vec{x}\cup \vec{x}'$, we have $p\star\vec{x}\in P^{\vec{x}\cap\vec{x}'}$ whenever $p\star\vec{x}\in P^{\vec{x}}\cap P^{\vec{x}'}$.

   If $p=x_l$, $l<k$, then
   \[  x_l\star\vec{x}= x_l\star\vec{x}'=x_l\star(\vec{x}\cap\vec{x}') \]
   and $x_l\star\vec{x}\in P^{\vec{x}\cap\vec{x}'}$.

If $p=x_l$, $l>k$, then
\[  (x_l\star\vec{x})_{(k)}= (x_l\star\vec{x}')_{(k)}=x_l\star(\vec{x}\cap\vec{x}') \]
and $(x_l\star\vec{x})_{(k)}\in P^{\vec{x}\cap\vec{x}'}$.

For $p=x_k,x_k'$, the lifts $p\star \vec{x}$ and $p\star \vec{x}'$  are different but they agree exactly on $p\star(\vec{x}\cap\vec{x}')$. To see this, fix $p=x_k$ and assume that $x_k\star\vec{x}'$ is a low p-flag (the cases when $x_k\star\vec{x}'$ is a high p-flag or $p=x'_k$ are similar). Then, using Proposition \ref{prop-monotone}, we have
\[  x_k\star_l\vec{x}\unlhd x_k\star(\vec{x}\cap\vec{x}') \unlhd x_k\star_l\vec{x}'=x_k\star_l\vec{x}\]
This shows that the common face $x_k\star_l\vec{x}= x_k\star_l\vec{x}'$ of both $P^{\vec{x}}$ and $P^{\vec{x}'}$ that projects to $x_k$ belongs to $P^{\vec{x}\cap\vec{x}'}$.

Ad 2. We shall consider the elements of $P^{\vec{y}}\cap P^{\vec{x}'}$  of all possible forms, and we will show that they belong to $P^{\vec{x}}$.

Suppose that the flag $\vec{z}$ is a face in  $P^{\vec{y}}\cap P^{\vec{x}'}$, $\vec{z}=\vec{x}'_{\lc l}$, for some $l$. Then, since $l\leq \diff(\vec{y},\vec{x}')\leq \diff(\vec{x},\vec{x}')$, we have that $\vec{z}$ is a face of $P^{\vec{x}}$, as well.

Suppose now that $p\not\in \vec{x}_{\lc k-1}$, i.e. $p\star\vec{x}'$ is either a p-flag or flat. Then it is the only face in $P^{\vec{x}'}$ that projects to $p$. Suppose $p\star\vec{x}'$ is a low p-flag (the case of $p\star\vec{x}'$ being a high p-flag is similar). If $p\star\vec{x}'\in P^{\vec{y}}$, then
$p\star\vec{x}'=p\star_l\vec{y}$ and,  by Proposition\ref{prop-monotone}, we have
  \[ p\star_l \vec{y}\unlhd p\star_l \vec{x}\unlhd p\star_l \vec{x}'= p\star \vec{x}'=  p\star_l \vec{y}. \]
Thus $p\star \vec{x}'= p\star_l \vec{x}$ is a face of $P^{\vec{x}}$.
$\Box$\vskip 2mm

\subsection{Cylinder $\Cylp(P)$ is aspherical in $\widehat{\pOpe}$}

\begin{proposition}\label{prop-cyl-straight-pHg} Let $P$ be a positive opetope. Then $\Cylh(P)$ is $\pOpe$-straight in $\pHg$.
  \end{proposition}
{\em Proof.} We have a linear order of maximal flags $\Flags_P[\bm_p,0)$. For $\vec{x}\in \Flags_P[\bm_p,0)$, we define
  \[ C^{\vec{x}} =\coprod_{\vec{y}\unlhd\vec{x}} P^{\vec{y}}. \]
 Let $\vec{x}$ and $\vec{x}'$ be two maximal flags with $\next(\vec{x})=\vec{x}'$.  By Proposition \ref{lemma-flag-intersection}, we have that
     $$C^{\vec{x}} \cap P^{\vec{x}'} = P^{\vec{x}\cap \vec{x}'}.$$
Since any face in $\Cylh(P)$ occurs in a maximal flag face, for the terminal maximal flag $\bT_P$, we have $C^{\bT_P}=\Cylh(P)$. Thus $\Cylh(P)$ is $\pOpe$-straight, as required. $\Box$\vskip 2mm

Let the presheaf $\Cylp(P)$ be the image of the hypergraph $\Cylh(P)$ under $\cH : \pHg\ra \widehat{\pOpe}$.

\vskip 2mm
\begin{theorem}\label{thm-cyl-aspherical-pOpe} Let $P$ be a positive opetope. Then $\Cylp(P)$ is aspherical in $\widehat{\pOpe}$.
  \end{theorem}
{\em Proof.} Since $\cH$ sends $\pOpe$-straight objects to straight objects in $\widehat{\pOpe}$, by Proposition \ref{prop-cyl-straight-pHg},
$\Cylp(P)$ is straight. Then, by  Proposition \ref{prop-aspherical-presh-characterization}, $\Cylp(P)$ is aspherical. $\Box$\vskip 2mm

We extend the construction  of the cylinder to a functor
$$\Cylp : \pOpe \lra \widehat{\pOpe}.$$
Let $f:P\ra Q$ be a morphism in $\pOpe$.
For $z\in \Cylp(P)$, we define
 \[ \Cylp(f)(z) =  \left\{
                                    \begin{array}{ll}
                                     -f(y) & \mbox{if } z=-y\in -P;\\
                                       & \\
                                     +f(y) & \mbox{if } z=+y\in +P;\\
                                       & \\
                                     \left[
                                       \begin{array}{c}

                                         f(x_k) \\
                                         \ldots\\
                                         f(x_0) \\
                                       \end{array}\right]   & \mbox{if } z= \left[
                                       \begin{array}{c}

                                         x_k \\
                                         \ldots\\
                                         x_0 \\
                                       \end{array}\right]  \\
                                    \end{array}
			    \right. \]
By convention $f(0)=0$.

 \begin{proposition}
  $\Cylp: \pOpe \lra \widehat{\pOpe}$ is a functor.
 \end{proposition}
  {\em Proof.} Exercise.  $\Box$\vskip 2mm

{\em Remark.} Taking left Kan extension along Yoneda we get a cylinder functor $$\Cylp :  \widehat{\pOpe}\lra \widehat{\pOpe}$$ extended to all positive opetopic sets. It is an elementary homotopic data in the sense of D-C. Cisinski (1.3.6).

\section{Cylinder as a product in $\widehat{\pOpei}$}\label{sec-cylinder-in-pOpei}

We fix two $\iota$-maps $h:Q\ra P$ and $\varrho: Q\ra I$ of positive opetopes for the whole section. We shall show that there is a unique $\iota$-map $H: Q\ra \Cyl_h(P)$ such that $\pi_p\circ H =h$ and $\varrho_I\circ H = \varrho$. We refer to this property that `opetopes thinks that $\Cyl_h(P)$ is a product of $I$ and $P$ in $\pHgi$'. This property will imply that $\Cyli(P)$ which is $\kappa_!(\Cylp(P))$, i.e. the image of $\Cyl_h(P)$ under left Kan extension along $\kappa:\pOpe\ra \pOpei$, is indeed a product of $I$ and $P$ in $\widehat{\pOpei}$.

\subsection{Projection $\iota$-maps}

 \begin{proposition}
  Let $P$ be a positive opetope. We have
  \begin{enumerate}
    \item the restricted projection map $\pi^{\vec{x}}:P^{\vec{x}}\lra P$ is a $\iota$-map.
    \item the projection map $\varrho^{\vec{x}}:P^{\vec{x}}\lra I$ such that, for $q\in P^{\vec{x}}$,
    \[ \varrho(q) =  \left\{
                                    \begin{array}{ll}
                                    -  &  \mbox{ if } q=-p \\
                                    +  &  \mbox{ if } q=+p,\\
                   a  &   \mbox{ othwerwise.} 	
                                    \end{array}
			    \right. \]
is a $\iota$-map.
\end{enumerate}
  \end{proposition}
  {\em Proof.} Simple check. $\Box$\vskip 2mm

\subsection{Splitting faces}

 $I=\{- \stackrel{\ba}{\ra} + \}$. Let $dim (P) =n$. We assume that $h$ is onto.

\begin{enumerate}
  \item We say that the face $q\in Q_1$ $\lk \varrho,h\rk$-{\em splits} face $p\in P_0$ iff $h(q)=1_p$ and $\varrho(q)=\ba$. We say that the face $q\in Q_{k+1}$ $\lk \varrho,h\rk$-{\em splits} face $p\in P_k$ with $k>0$ iff
\begin{enumerate}
  \item $h(q)=p$ and $\rho(q)=\ba$;
  \item There is a face $q'\in \delta(q)$ such that $q'$ $\lk \varrho,h\rk$-splits $h(q')\in P_{k-1}$.
\end{enumerate}
We often say splitting face when we mean $\lk \varrho,h\rk$-splitting. If $q$ is a splitting face, then $p=h(q)$ is the {\em splitted face}.

\vskip 2mm
{\em Example.} Let $\vec{x}$ be a flag in $P$. Putting $Q=P^{\vec{x}}$ and $h=\varrho_P : P^{\vec{x}}\ra P$ and $\varrho:P^{\vec{x}}\ra I$ such that for $f\in P^{\vec{x}}$
     \[ \varrho(f) =  \left\{
                                    \begin{array}{ll}
                                    -  &  \mbox{ if } f=-\varrho_P(f),\\
                                    +  &  \mbox{ if } f=+\varrho_P(f),\\
                   \ba  &   \mbox{ if } f \mbox{ is a flag or p-flag.}	
                                    \end{array}
			    \right. \]
Then the flag $[x_k,\ldots, x_0]$ splits $x_k$.
  \item  We say that $q\in Q_1$ is a {\em threshold face} if $h(q)\in P_1$ and $varrho(q)=\ba$.  We say that $q\in Q_{\geq 2}$ is a {\em threshold face} iff $h(q)\not\in ker(h)$ and   there is a splitting face $s\in \delta(q)$.
  \item Recall that the set $X$ of $k$-faces in an opetope $P$ is called a {\em $<^+$-interval} if it is either empty or there are two k-faces $x_0,x_1\in P_k$
such that $x_0\leq^+ x_1$ and $X=\{ x\in P_k | x_0\leq^+ x \leq^+ x_1 \}$.
  \item  A $<^+$-interval $X$ is {\em initial} iff $x_0\in P_k-\gamma(P_{k+1})$.
  \item Let $X=\{ x_0,\ldots, x_u\}$ be an $<^+$-interval in $P_k$ (listed in $<^+$-order).
The subset $Y=\{ y_1,\ldots, y_{u} \}$ (listed in $<^-$-order) of $P_{k+1}-\gamma(P_{k+2})$ is a {\em witness $<^-$-interval for $X$} iff $x_{i-1}\in \delta(y_i)$ and $\gamma(y_i)=x_{i}$, for $i=1,\ldots u$.
\end{enumerate}

\vskip 2mm
  \begin{proposition}(splitting and threshold faces)\label{prop-splitting-faces}
Let $P$, $\varrho$ and $h$ be as above $h(\bm_Q)=\bm_P=x_n$, $\varrho(\bm_Q)=\ba$. Then there is $0< k\leq n$  and two $<^+$ intervals
of faces of $Q$, for $1\leq i\leq k\leq dim(Q)$,
\[ S^i=\{ s^i_0,\ldots, s^i_{u_i}\},\;\;\; T^i=\{ t^i_0,\ldots, t^i_{v_i} \},\]
and two witness sequences of faces in $Q-\gamma(Q)$
\[ A^{i+1}=\{ a^{i+1}_1,\ldots, a^{i+1}_{u_i}\},\;\;\; B^{i+1}=\{ b^{i+1}_1,\ldots, b^{i+1}_{v_i} \},\]
such that
\begin{enumerate}
\item $u_i$, $v_i$ are numbers $\geq 0$, for $i=1,\ldots, k-1$; $u_k,v_k,u_k+v_k\geq -1$ and either $u_k=-1$ or $v_k=-1$ (i.e. one but not two of the sequences $S^k$, $T^k$ is empty);
  \item $S^i$ is an initial $<^+$-interval in $Q_{i}$ of all splitting faces in $Q$ of dimension $i$;
  \item $A^i$  is the witness $<^-$-interval for $S^i$ of 2-collapsing faces in  $Q_{i+i}-\gamma(Q)$;
  \item $T^i$ is a $<^+$-interval in $Q_{i}$ of all threshold faces in $Q$ of dimension $i$;
  \item $B^i$  is the witness $<^-$-interval for $T^i$ of 1-collapsing and non-collapsing faces in  $Q_{i+i}-\gamma(Q)$;
  \item $S^i\cup T^i$ is a maximal $<^+$-interval in $Q_{i}$; let $\min^i$ $\max^i$ be the least and the largest elements in this set;
  \item $A^i\cup \{s^{i+1}_0\}\cup B^i$ is the witness $<^-$-interval for  $S^i\cup T^i$, for $i<k$;
  \item  for $j=1,\ldots,v_i$,
  \begin{enumerate}
    \item  if $b^i_j$ is non-collapsing, then $h(t^i_{j-1})<^+h(t^i_j)$;
    \item if $b^i_j$ is 1-collapsing, then \mbox{$h(t^i_{j-1})=h(t^i_j)$}.
  \end{enumerate}
  \item for $i=1,\ldots, k-1$, we have
  \[ \min{}^i\in \delta\gamma^{(i+1)}(\max{}^{i+1}), \;\;\;\;\; \gamma^{(i)}(\max{}^{i+1})=\max{}^i.\]
\end{enumerate}
  \end{proposition}

\vskip 2mm

{\em Example.} An example of a sequences as in the above Proposition with $k=4$, $u_4=-1$ is
\[
\begin{array}{cccccccccccccccccc}
  s^1_0 & a^2_1 & \ldots & a^2_{u_1} & s^1_{u_1} & s^2_0 & t^1_0 & b^2_1 & \ldots & b^2_{v_1} & t^1_{v_1} &  &  &  &  &  &  &  \\
   & &  &  &  &  &  &  &  &  & \parallel &  & &  &  &  &  &  \\
   & &  &  &  & a^3_1 &  &  &  &  & \gamma^{(1)}(q) &  &  &  &  &  &  &  \\
   & &  &  &  &  &  &  &  &  &  &  &  &  &  &  &  &  \\
   & &  &  &  & s^2_1 &  &  &  &  &  &  &  &  &  &  &  &  \\
   & &  &  &  &  &  &  &  &  &  &  &  &  &  &  &  &  \\
   & &  &  &  & \ldots &  &  &  &  &  &  &  &  &  &  &  &  \\
   & &  &  &  &  &  &  &  &  &  &  &  &  &  &  &  &  \\
   & &  &  &  & a^3_2 &  &  &  &  &  &  &  &  &  &  &  &  \\
   & &  &  &  &       &  &  &  &  &  &  &  &  &  &  &  &    \\
   & &  &  &  & s^2_{u_2} &  &  &  &  &  &  &  &  &  &  &  &  \\
   & &  &  &  &  &  &  &  &  &  &  &  &  &  &  &  &  \\
   & &  &  &  & s^3_0 & a^4_1 & \ldots & a^4_{u_3} & s^3_{u_3} & t^4_0 & t^3_0 & b^4_1 & \ldots & b^4_{v_3} & t^3_{v_3} &  &  \\
   & &  &  &  &  &  &  &  &  &  &  &  &  &  & \parallel &  &  \\
   & &  &  &  & t^2_0 &  &  &  &  &  b^5_0 &&  &  &  & \gamma^{(2)}(q) &  &  \\
   & &  &  &  &  &  &  &  &  &  &  &  &  &  &  &  &  \\
   & &  &  &  & b^3_1 &  &  &  &  &   t^4_1 & &  &  &  &  &  &  \\
   & &  &  &  &  &  &  &  &  &  &  &  &  &  &  &  &  \\
   & &  &  &  & \ldots &  &  &  & &  \ldots & &  &  &  &  &  &  \\
   & &  &  &  &  &  &  &  &  &  &  &  &  &  &  &  &  \\
   & &  &  &  & b^3_{v_2} &  &  &  &  &  b^5_{v_4} & &  &  &  &  &  & \\
   & &  &  &  &  &  &  &  &  &  &  &  &  &  &  &  &  \\
   & &  &  &  &  t^2_{v_2} &  &  &  &  &   t^4_{v_4} & &  &  &  &  &  & \\
   & &  &  &  &  \parallel &  &  &  &  &  \parallel & &  &  &  &  &  & \\
   & &  &  &  & \gamma^{(2)}(q) &  &  &  &  &   \gamma^{(4)}(q) & &  &  &  &  &  & \\
\end{array}
\]
with
\[ s^1_0\in \delta\gamma^{(2)}(q), \;\; s^2_0\in \delta\gamma^{(3)}(q), \;\; s^3_0\in \delta\gamma^{(3)}(q).\]

{\em Proof.} The proof proceeds by induction. Before we state the inductive hypothesis, we shall make some comments. Note that in order to have any splitting faces whatsoever, we need to have one face $x_1$ such that  $\varrho(x_1)=\ba$ and $h(x_1)=p_0\in P_0$. If $x_1$ is such a splitting 1-face and we have a face $A\in Q_2$ such that $\gamma(A)=x_1$, then there is another splitting face in $\delta(a)$, this $x_1'\in\delta(A)$ for which we have $\varrho(x_1')=\ba$. This shows that the splitting 1-faces form an initial $<^-$-interval. Since $x_1$ is a 1-collapsing face as all splitting faces are, $A$ is a 2-collapsing face.

In particular, we can start counting splitting 1-faces with the one in $Q_1-\gamma(Q_2)$. On the other hand, if $s^1\in\delta(A)$, then we can have that either $A$ is 2-collapsing and $\gamma(A)$ is again a splitting 1-face, or $A$ is 1-collapsing and $A$ is a splitting 2-face, or $A$ is not collapsing, i.e. $h(A)\in P_2$, and there are no more 1-splitting faces and no $\geq 2$ splitting faces. In this case the codomain $\gamma(A)=t^1$ of $A$ is a threshold 1-face.
If we have a threshold 1-face $t$ in domain of a 2-face $A$, then $\gamma(A)$ is again a threshold 1-face.
This situation repeats at the higher dimension. However, to show this, we will need the assumption that the splitting faces at the lower dimensions are already known to form initial $<^+$ intervals. Thus we have to state rather elaborate inductive hypothesis with a list of assumptions in order to be able to obtain all the conclusions.

Thus we assume that the situation as described in Lemma is indeed true at dimensions smaller than $l$, and possibly we already constructed a part of the maximal sequence $S^l\cup T^l$ at the level $l$. Then
\begin{enumerate}
  \item either we already constructed a splitting l-face $s^l_i$ and we have a face $a\in Q_{l+1}-\gamma(Q_{l+2})$ such that $s^l_i\in\delta(a)$;
  \begin{enumerate}
    \item  if $a$ is 2-collapsing, then we put $s^l_{i+1}=\gamma(a)$, the next splitting $l$-face, $a^l_{i+1}=a$;
    \item if  $a$ is 1-collapsing, then we put $t^l_{0}=\gamma(a)$, the first threshold face at dimension $l$, $s^{l+1}_0=a$ is the first splitting face at dimension $l+1$;
    \item if  $a$ is non-collapsing, then we put $t^l_{0}=\gamma(a)$ the first threshold face at dimension $l$, $t^{l+1}_0=a$ the first threshold face at dimension $l+1$; $u_{l+1}$ is $1-1$; in this case the construction ends at dimension $l+1$ where we have threshold faces only;
  \end{enumerate}
    \item or we already constructed all splitting l-faces and some threshold faces till $t^l_i$ and we have a face $a\in Q_{l+1}-\gamma(Q_{l+2})$ such that $t^l_i\in\delta(a)$; then  $a$ is 1-collapsing or non-collapsing and we put $t^l_{i+1}=\gamma(a)$, $b^{l+1}_0=a$,
    \item or we already constructed a splitting l-face $s^l_i$ and we have no face $a\in Q_{l+1}-\gamma(Q_{l+2})$ such that $s^l_i\in\delta(a)$; then we stop the construction;
     \item or we already constructed all splitting l-faces and some threshold faces till $t^l_i$ and we have no face $a\in Q_{l+1}-\gamma(Q_{l+2})$ such that $t^l_i\in\delta(a)$; then we proceed to construct the interval  $S^{l+1}\cup T^{l+1}$ at dimension $l+1$.
\end{enumerate}
$\Box$\vskip 2mm

{\em Notation $\sigma$, $\tau$ and $\xi$.}  Let $q\in Q_{k+1}$, $k> 0$. From the above proposition follows that if there is a splitting face in $\delta(q)$, it is unique and we shall denote it as $\sigma(q)\in \delta(q)$. Similarly, if $k>1$ and there is a threshold face $z\in \delta(q)$,  it is unique and we shall denote it as $\tau(q)\in \delta(q)$. Finally, if $q$ is a 1-collapsing (e.g. splitting) face, then $\xi(q)$ denotes the only face in $\delta(q)$ such that $h(\xi(q))=h(\gamma(q))$.

\vskip 2mm

\begin{lemma}\label{lemma-tau-si-xi-well-def}
Let $q,q'\in Q_{k+1}$. Then
\begin{enumerate}
  \item\label{le-tsx-2}  $\sigma(q)$ and $\tau(q)$, if they exist, are unique;
  \item\label{le-tsx-3} if $q$ and $q'$ are at most 1-collapsing, and $\sigma(q)$ and $\sigma(q')$ are defined, then $q\not\perp^-q'$;
  \item\label{le-tsx-4}  if  $\sigma(q)$ is defined, then $\tau(q)$ is not defined;
  \item\label{le-tsx-5} if $q$ is a splitting face, then $\sigma(q)\perp^-\xi(q)$;
  \item\label{le-tsx-6} if $l,l'\geq 0$, $\tau(\gamma^{{l+2}}(q))$ and  $\tau(\gamma^{{l'+2}}(q))$ are defined, then $l=l'$.
\end{enumerate}
\end{lemma}
{\em Proof.}
\ref{le-tsx-2}. and  \ref{le-tsx-4}. easily follow from Proposition \ref{prop-splitting-faces}.

Ad \ref{le-tsx-3}. Suppose we have both $\sigma(q)$ and $\sigma(q')$ defined.  Assume $q<^-q'$. Then, by Proposition \ref{prop-splitting-faces}, we have $\sigma(q)\leq^+\sigma(q')$. Then $\sigma(q)\leq^+\gamma(q)\leq^+t\in\delta(q')$ with $t$ non-collapsing, and hence $t\neq\sigma(q')$. But the faces $\leq^+$-greater than any face, e.g. $\sigma(q)$, are linearly ordered. Thus $t\perp^+\sigma(q')$. But $t,\sigma(q')\in\delta(q')$ and we get a contradiction.
Changing the roles of $q$ and $q'$ we get that  $q'<^-q$ does not hold either and hence $q\not\perp^-q'$.

Ad \ref{le-tsx-5}. Since $h(\xi(q))=h(\gamma(q))$, it follows  $h(\gamma(\xi(q)))=h(\gamma^2(q))$ and $h(\gamma^2(q))$ is $<^+$-maximal in $P_{k-1}[h(q)]$. Thus $h(\gamma\sigma(q))\leq^+h(\gamma(q))$. Hence, by Corollary \ref{coro-iota-preservation},  we have $\gamma\sigma(q)\perp^+\gamma\xi(q)$, since both faces are not collapsing. We have that $\xi(q),\sigma(q)\in \delta(q)$, so we can't have $\xi(q)\perp^+\sigma(q)$. Therefore   $\xi(q)\perp^-\sigma(q)$.

Ad \ref{le-tsx-6}. Suppose to the contrary that $l>l'$ and  $\tau(\gamma^{{l+2}}(q))$ and  $\tau(\gamma^{{l'+2}}(q))$ are defined. Then
$\sigma^{(l'+2)}\tau\gamma^{(l+2)}(q)\leq^+ \gamma^{(l'+2)}(q)$.  But since splitting faces of the same dimension form an initial interval, there is a splitting face $\sigma(\gamma^{(l'+2)}(q)) \in \delta(\gamma^{(l'+2)}(q))$. Thus both $\tau$ and $\sigma$ are defined on the face $\gamma^{(l'+2)}(q)$. This contradicts  \ref{le-tsx-4}.
 $\Box$
 \vskip 2mm

\begin{corollary}
Let $q\in Q_{k+1}$, $k>0$, be 2-collapsing face so that $\sigma(q)$ is defined. Then $\gamma(q)$ is again a splitting face.
\end{corollary}
 {\em Proof.} For $k=1$, Corollary is obvious. Assume $k>1$. Since $q$ is $2$-collapsing, $\gamma(q)$ is 1-collapsing.  Moreover, $\delta\gamma(q)\leq^+ \sigma^2(q)$ and hence the initial interval of splitting faces of dimension $k-1$ intersects $\delta\gamma(q)$, see Proposition \ref{prop-splitting-faces}. $\Box$

\subsection{The $\iota$-map $H$ into product in $\pHgi$}

Given $\iota$-maps $h$ and $\varrho$ as above, we define a $\iota$-map
\[ H: Q\lra \Cylh(P)\]
as follows. Let $q\in Q$.
If $\varrho(q)\in \{ -,+\}$, then we put
 \[ H(q) =  \left\{  \begin{array}{ll}
                      -h(q) & \mbox{ if } \varrho(q)=-\\
                      +h(q) & \mbox{ if } \varrho(q)=+
   \end{array}
  \right. \]

If $\varrho(q)=\ba$, we put

 \[ H(q) =  \left\{
                                    \begin{array}{lll}
                               (H1,S1)  &   [h(q)] &  \begin{array}{l}
                                                   \mbox{ if } q\in Q_{1} \mbox{ is a splitting face, i.e. } h(q)\in P_0\\
                                                   (H(q) \mbox{ is a flag of length 1});
                                              \end{array} \\
                & & \\
                               (H2,N1)  &   {[h(q),0]} &  \begin{array}{l}
                                                   \mbox{ if } q\in Q_{1} \mbox{ is not a splitting face, i.e. } h(q)\in P_1\\
                                                   (H(q) \mbox{ is a p-flag of length 2});
                                              \end{array} \\
                & & \\
                                (H3,A)&       {[h(q),H(\sigma(q))]}  &  \begin{array}{l}
                                                 \mbox{ if } q\in Q_{\geq 2} \mbox{ is a splitting face} \\
                                                 (H(q) \mbox{ is a flag });
                                              \end{array} \\
                & & \\
                   (H4,B)& {[h(q),0,H(\sigma(q))]}  & \begin{array}{l}  \mbox{ if } q\in Q_{\geq 2}-\ker(h) \mbox{ and } \sigma(q) \mbox{ is defined} \\
                                              (H(q) \mbox{ is a high p-flag});
                                              \end{array} \\
                & & \\
                   (H4.1,B)& {[h(q),0,H(\sigma(q))]}  & \begin{array}{l}  \mbox{ if } q\in Q_{\geq 2} \mbox{ is a threshold face} \\
                                              (H(q) \mbox{ is a high p-flag});
                                              \end{array} \\
                & & \\
                              (H5,C)    &  {[h(q),H(\tau(q))]}  &  \begin{array}{l}
                                                 \mbox{ if } q\in Q_{\geq 2}-\ker(h) \mbox{ and } \tau(q) \mbox{ is defined}\\
                                                 (H(q) \mbox{ is a low p-flag of codim 2});
                                              \end{array} \\
                & & \\

        (H6,D) & {[h(q),H(\gamma(q))]}   &
                                       \begin{array}{l}
                                       \mbox{otherwise, if } q\in Q_{\geq 2}-\ker(h) \\
                                       (H(q) \mbox{ is a low p-flag of codim } > 2);
                                       \end{array}
                                       \\
           &  &  \\
	   (H7,E)& H(\gamma(q))  & \mbox{ otherwise } (H(q)  \mbox{ is an  identity on a face}).		
                                    \end{array}
			    \right. \]

\begin{lemma}\label{lemma-H-unique}
Let $\vec{x},\vec{x}'\in \Cylh(P)$, with $\pi_P(\vec{x})=\pi_P(\vec{x}')=x_{k+1}\in P_{k+1}$, $\delta(\vec{x})=\delta(\vec{x}')$ and $\gamma(\vec{x})=\gamma(\vec{x}')$. Then $\vec{x}=\vec{x}'$.
\end{lemma}

{\em Proof.} For $k\leq 1$ or for flat faces $\vec{x},\vec{x}'\in \Cyli(P)$, the Lemma is obvious.

So assume that $k>1$ and $\rho_I(\vec{x})=\rho_I(\vec{x}')=\ba$.

Since $\gamma(\vec{x})=\gamma(\vec{x}')$, it follows that $\vec{x}$ and $\vec{x}'$ have the same dimension. As $\pi_P(\vec{x})=\pi_P(\vec{x}')$, either both  $\vec{x}$ and $\vec{x}'$ are flags or p-flags.

Let $\vec{x}=[x_{k+1},\ldots,x_0]$ and $\vec{x}'=[x_{k+1},x'_{k-1}\ldots,x'_0]$ be flags.  Then side flags agree, i.e.
\[ [x_{k},\ldots,x_0]=\vec{x}_{side}=\vec{x}'_{side}=[x'_{k},\ldots,x'_0].\]
As we have $x_{k+1}=\pi_P(\vec{x})=\pi_P(\vec{x}')=x'_{k+1}$, so $\vec{x}=\vec{x}'$, as required.

Now let $\vec{x}=[x_{k+1},\ldots,\widehat{x_l},\ldots,x_0]$ and $\vec{x}'=[x'_{k+1},\ldots,\widehat{x'_j},\ldots,x'_0]$ be p-flags.

The largest flag occurring in faces $\partial(\vec{x})$ is  $[x_{l-1},\ldots,x_0]$. The largest flag occurring in faces $\partial(\vec{x}')$ is  $[x'_{j-1},\ldots,x'_0]$. Since $\partial(\vec{x})=\partial(\vec{x}')$, we have that $l=j<k$ and
$$[x_{l-1},\ldots,x_0]=[x'_{l-1},\ldots,x'_0].$$
Since $x_{k+1}=x'_{k+1}$, we have
\[ \vec{x} =[x_{k+1},\gamma(x_{k+1}),\ldots,\gamma^{(l+2)}(x_{k+1}), x_{l+1},0,x_{l-1},\ldots,x_0] \]
and
\[ \vec{x}' =[x_{k+1},\gamma(x_{k+1}),\ldots,\gamma^{(l+2)}(x_{k+1}), x'_{l+1},0,x_{l-1},\ldots,x_0] \]
i.e. $\vec{x}$ may differ from $\vec{x}'$ only at level $l+1$. Suppose that $x'_{l+1}\neq x_{l+1}$.

If $l+1<k$, then the face $\vec{x}_{\lc k}\in\partial(\vec{x})$, the  only face in $\partial(\vec{x})$ that projects to $\gamma(x_{k+1})$, is different than the face $\vec{x}'_{\lc k}\in\partial(\vec{x}')$, the only face $\partial(\vec{x}')$ that projects to $\gamma(x_{k+1})$. Thus $\partial(\vec{x})\neq\partial(\vec{x}')$, contrary to the assumption.

If $l+1=k$, then
\[ \vec{x} =[x_{k+1},x_k,0,x_{k-2},\ldots,x_0] \]
and
\[ \vec{x}' =[x_{k+1},x'_k,0,x_{k-2},\ldots,x_0] \]
with $x_k \neq x'_k$. Then $x_k, x'_k\in \delta(x_{k+1})$ and hence $x_k\not\perp^+ x'_k$. Thus $x_k\star \vec{x}'\in \delta(\vec{x}')$ is a low p-flag and
it is the only face in $\delta(\vec{x}')$ that projects to $x_k$. On the other hand, the only face in $\delta(\vec{x})$ that projects to $x_k$ is $x_k\star \vec{x}=\vec{x}_{\lc k}$ and it is a high p-flag. Thus the sets $\delta(\vec{x})$ and $\delta(\vec{x}')$ cannot be equal in this case either.
This shows that $\vec{x}=\vec{x}'$, as required.
$\Box$\vskip 2mm

\begin{proposition}\label{prop-univ-prop-cylh-pHgi}
The map $H$ defined above is the unique $\iota$-map $H: Q\ra \Cylh(P)$ such that $\varrho=\varrho_I\circ H$ and $h=\pi_P\circ H$.
\end{proposition}

{\em Proof.} The equalities $\varrho=\varrho_I\circ H$ and $h=\pi_P\circ H$ are obvious. We need to check that $H$ is a $\iota$-map, i.e. it preserves codomains and domains and that it is unique.

That $H$ preserves codomains and domains is the content of Propositions \ref{prop-H-pres-codomains} and \ref{prop-H-pres-domains} below.
We show below that $H$ is unique.

The proof goes by induction of the level of $q\in Q_{k+1}$. We proceed by cases that define the map $H$. For $k=0,1$ or $\varrho(q)\neq\ba$, the proof is easy. Thus we assume that $k>0$ and $\varrho(q)=\ba$.

Now assume that, for $q\in Q_{k+1}$, $k>1$ and for $\bar{q}\inn q$, $\bar{q}\neq q$, $H(\bar{q})$ is the unique value making $\varrho(\bar{q})=\varrho_I\circ H(\bar{q})$ and $h(\bar{q})=\pi_P\circ H(\bar{q})$. Let $H':Q\ra \Cyli(P)$ be a $\iota$-map such that $H'(\bar{q})=H(\bar{q})$, for $\bar{q}\inn q$, $\bar{q}\neq q$, and $\pi_P(H'(q))=h(q)$. We need to show that $H(q)=H'(q)$, as well. But this is the content of Lemma \ref{lemma-H-unique}.
$\Box$\vskip 2mm

\subsection{Splitting sequences}

We define below, for some faces $q\in Q$, a sequence for faces in $Q$ that we call a splitting sequence for $q$. We shall show that for any face $q\in Q$ such that $\rho(q)=\ba$, i.e. those $q$ for which the values $H(q)$ are not flat, such a sequence exists and is unique. Moreover, the value of $H(q)$ is simply calculated from  the splitting sequence for $q$.

Let $q\in Q$. A {\em splitting sequence for} $q$ is defined as follows.
\begin{enumerate}
  \item if $q\in Q_0$, then there is no splitting sequence for $q$;
  \item if $q\in Q_1$ and $\rho(q)=\ba$, then the splitting sequence for $q$ is
  \[ [q]; \]
  \item if $q\in Q_{\geq 2}-ker(h)$ and there is $l\geq 0$ such that $\tau\gamma^{(l+2)}(q)$ is defined, then the splitting sequence for $q$ is
  \[ [q,\gamma(q),\ldots,\tau\gamma^{(l+2)}(q),\sigma\tau\gamma^{(l+2)}(q),\ldots,\sigma^{(1)}\tau\gamma^{(l+2)}(q) ]; \]
  \item if $q\in Q_{\geq 2}$ is 1-collapsing and $\sigma(q)$ is defined, then the splitting sequence for $q$ is
  \[ [q,\sigma(q),\ldots,\sigma^{(1)}(q) ]; \]
  \item otherwise, if the splitting sequence for $\gamma(q)$ is defined, then the splitting sequence for $q$ is the splitting sequence for $\gamma(q)$, ;
  \item no other splitting sequences are defined.
\end{enumerate}

\begin{lemma}\label{lemma-split-seq}
Let $q\in Q$. Then
\begin{enumerate}
  \item\label{le-ss-1} there is a splitting sequence for $q$ iff $\varrho(q)=\ba$;
  \item\label{le-ss-2} if a splitting sequence exists for $q$, then it is unique;
  \item\label{le-ss-3} if $q\not\in ker(h)$ and a splitting sequence for $q$ exists, then
  \[ \hskip -10mm H(q)= [h(q),h(\gamma(q)),\ldots,h(\tau\gamma^{(l+2)}(q)),0,h(\sigma\tau\gamma^{(l+2)}(q)),\ldots,h(\sigma^{(1)}\tau\gamma^{(l+2)}(q)) ]; \]
  \item\label{le-ss-4} if $q$ is 1-collapsing and a splitting sequence for $q$ exists, then
  \[ \hskip -10mm H(q)= [h(q),h(\sigma(q)),\ldots,h(\sigma^{(1)}(q)) ]; \]
  \item\label{le-ss-5} otherwise, if a splitting sequence for $\gamma(q)$ exists, then $H(q)=H(\gamma(q))$;
  \item\label{le-ss-6} Let $x\in Q_1$ such that $\varrho(x)=\ba$. If there is no splitting sequence for $q$, then
  \begin{enumerate}
    \item either $\gamma^{(1)}(q)<^-x$ and then $H(q)=-h(q)$,
    \item or $x<^-\gamma^{(1)}(q)$ and then  $H(q)=+h(q)$.
  \end{enumerate}
\end{enumerate}
\end{lemma}

{\em Proof.} Ad 1, 2, 3, 4 and 5. Follows immediately from the definitions of a splitting sequence and morphism $H$.

Ad \ref{le-ss-6}. Since we have for any faces $q_1,q_2\in Q_1$ that either $q_1\perp^+ q_2$ or $q_1\perp^- q_2$, we can assume that $x\in  Q_1-\gamma(Q_2)$.  Since there is no splitting sequence for $q$, $x\not\leq^+q$ and hence $x\not\perp^+q$. Thus $x\perp^-q$. The rest is easy.
$\Box$\vskip 2mm

\subsection{Some technical lemmas}

Before we prove that $H$ preserves codomains and domains, we need to prove some Lemmas.

\vskip 2mm
\begin{lemma}\label{lemma-H-collaps-ga}
Let $k\geq 1$, $q\in Q_{k+1}$ a 1-collapsing but not splitting. Then
\begin{equation}\label{equation-H-lemma}
H(\xi(q))=H(\gamma(q)).
\end{equation}
\end{lemma}

{\em Proof.} If $\rho(q)\neq\ba$, the Lemma is obvious. In particular, if $q\in Q_1$, it holds.  If $q\in Q_2$ and $\rho(q)=\ba$, then both $\xi(q)$ and $\gamma(q)$ are threshold faces and again Lemma holds. So we assume from now on that $q\in Q_{\geq 3}$.

We shall show that we have three excluding cases
\begin{enumerate}
  \item $\sigma(\xi(q))\da$ iff $\sigma\gamma(q)\da$,
  \item $\tau(\xi(q))\da$ iff $\tau\gamma(q)\da$,
  \item $\sigma(\xi(q))\ua$, $\tau(\xi(q))\ua$, $\sigma\gamma(q)\ua$, $\tau\gamma(q)\ua$,

\end{enumerate}
  and that in either case

\begin{equation} \label{equation-H(t)=H(g(q))}
  H(\xi(q))=H(\gamma(q).
\end{equation}

\vskip 2mm

Ad 1. Assume that $\sigma(t)\da$. Then, as $\xi(q)\leq^+\gamma(q)$, by Proposition \ref{prop-splitting-faces}, we have $\sigma\gamma(q)\da$ and hence
\[ H(\xi(q))=[h(\xi(q)),H(\sigma(\xi(q)))]=[h(\gamma(q)),H(\sigma\gamma(q))]= H(\gamma(q)).\]

Suppose now that  $\sigma\gamma(q)\da$. If we have a splitting $\sigma(\xi(q))$, then we proceed as above. So assume that $\sigma(\xi(q))$ is not defined. If there is a face $z\in\delta(\xi(q))$ such that $\sigma\gamma(q)\leq^+z$, then $z$ must be a threshold face and we have an upper $\delta(q)$-path of 1- and 2-collapsing faces
\[ \sigma\gamma(q),a_1,\ldots, a_r,z.\]
Then $a_j$ with $j=\max\{i: \sigma(a_i)\da \}$ is a splitting face and, as $q$ is 1-collapsing, it is also a splitting face contrary to the assumption. Thus we cannot have $\sigma\gamma(q)\leq^+\delta(\xi(q))$. Thus we must have a face $t'\in\delta(q)$ such that there are $z,z'\in\delta(t')$ such that $\gamma(\xi(q))\leq^+z$  and $\sigma\gamma(q)\leq^+z'$ with $z\neq z'$.

As $\xi(q)$ is the only non-collapsing face in $\delta(q)$, $t'$ is collapsing. $z$ is non-collapsing as $\gamma(\xi(q))$ is non-collapsing.  Then $z'$ must be collapsing. But if $z'$ is collapsing and it is $<^+$-greater than a splitting face, $z'$ must be a splitting face. Then $t'$ is a splitting face and hence $q$ is a splitting face, contrary to the supposition.  Thus $\sigma(\xi(q))\da$ iff $\sigma\gamma(q)\da$ and if either condition is true, (\ref{equation-H(t)=H(g(q))}) holds.

\vskip 2mm

Ad 2. If $\tau(\xi(q))\da$, then, by Proposition \ref{prop-splitting-faces}, we must have either that  $\sigma\gamma(q)\da$ or that $\tau\gamma(q)\da$. The former is impossible as we saw above. In the latter case we have a $\delta(q)$-path of 1-collapsing faces from $\tau(\xi(q))$ to  $\tau\gamma(q)$ and hence
\[ H(\xi(q))=[h(\xi(q)),H(\tau(\xi(q)))]=[h(\gamma(q)),H(\tau\gamma(q))]= H(\gamma(q)).\]

We shall show that if $\tau\gamma(q)\da$, then $\tau(\xi(q))\da$. So assume to the contrary that $\tau\gamma(q)\da$ and $\tau(\xi(q))\ua$. Thus we must have a face $t'\in\delta(q)$ such that there are $z,z'\in\delta(t')$ such that $\tau\gamma(q)\leq^+z$ and $\gamma(\xi(q))\leq^+z'$ with $z\neq z'$. As $\xi(q)$ is the only non-collapsing face in $\delta(q)$, $t'$ is collapsing and $z$, $z'$ are not collapsing. This is impossible.

Ad 3. Now assume that  $\sigma(\xi(q))\ua$, $\tau(\xi(q))\ua$, $\sigma\gamma(q)\ua$, and $\tau\gamma(q)\ua$. Then, as $\xi(q)$ and $\gamma(\xi(q))$ are non-collapsing, we have an upper $\delta(q)$-path of 1-collapsing faces
\[ \gamma(\xi(q)),b_1,\ldots, b_s,\gamma^2(q).\]
If any of $b_i$'s would be splitting face, $q$ would be splitting, as well. Thus using inductive assumption i.e. this Lemma for $q$'s of lower dimension, we get that $H(\gamma(\xi(q)))=H(\gamma^2(q))$. Then we have
\[ H(\xi(q))=[h(\xi(q)),H(\gamma(\xi(q)))]=[h(\gamma(q)),H(\gamma^2(q))]= H(\gamma(q)),\]
as required.
$\Box$\vskip 2mm

\vskip 2mm
\begin{lemma}\label{lemma-H-xi-sigma} Let $k\geq 1$, $q\in Q_{k+1}$ be a splitting face.
\begin{enumerate}
    \item If $\sigma(q)<^-\xi(q)$, then
  \begin{enumerate}
    \item $\gamma(\sigma(q))=\tau(\xi(q))\in \delta(\xi(q))$ is defined
    \item $H(\xi(q))=[h(\xi(q)),H(\tau(\xi(q)))]=H(q)_{low}$.
\end{enumerate}
  \item If $\xi(q)<^-\sigma(q)$, then
  \begin{enumerate}
    \item $h(\sigma(q))=h(\gamma\sigma(q))=h(\gamma^2(q))=h(\gamma(\xi(q)))$,
    \item $H(\xi(q))=[h(\xi(q)),H(\sigma(q))_{low}]=H(q)_{low}$.
  \end{enumerate}

\end{enumerate}
\end{lemma}

{\em Proof.} Ad 1. We have $h(\sigma(q))\in\delta(\xi(q))$.
Moreover
\[ \sigma^2(q)<^+\gamma\sigma(q)\leq^+ \tau\xi(q)\in \delta(\xi(q))\]
where $\tau\xi(q)$ is the unique face in $\delta(\xi(q))$ $<^+$-greater than $\sigma^2(q)$. Since it is greater than $\gamma\sigma(q)$, it is non-collapsing and hence a threshold face. This justifies (a).

We also have
\[ H(\xi(q)) = [h(\xi(q)),H(\tau(\xi(q)))= \]
(Lemma \ref{lemma-H-collaps-ga})
\[ =[h(q),H(\gamma\sigma(q) ]= \]
($\gamma\sigma(q)$ is a threshold face)
\[ =[h(q),h(\gamma\sigma(q)),0,H(\sigma\gamma\sigma(q)) ] = \]
($H$ on splitting faces of the same dimension is equal)
\[ =[h(q),h(\sigma(q)),0,H(\sigma^3(q)) ] =  \]
(first remark above)
\[ = H(q)_{low} \]
as required.
\vskip 2mm

Ad 2. Part (a) is clear. We have $h(\sigma(q))=\gamma(h(q))$. Moreover, as we have an upper $\delta(q)$-path from $\gamma\xi(q)$ to $\xi\sigma(q)$, by Lemma \ref{lemma-H-collaps-ga} we have that $H(\gamma\xi(q))=H(\xi\sigma(q))$. Since $\gamma\xi(q)\leq^+\xi\sigma(q)\perp^-\sigma^2(q)$, it follows that neither $\tau\xi(q)$ nor $\sigma\xi(q)$ is defined. Thus we have

\[ H(\xi(q)= [h(\xi(q)),H(\gamma\xi(q))]=\]
(Lemma \ref{lemma-H-collaps-ga})
\[ =  [h(q),H(\xi\sigma(q))]=\]
(induction)
\[ = [h(q),H(\sigma(q))_{low}]= \]
(second remark above)
\[ = [h(q),H(\sigma(q))]_{low}= H(q)_{low} \]
as required. $\Box$\vskip 2mm

\vskip 2mm
\begin{lemma}\label{lemma-H-collaps}
Let $k\geq 1$, $q\in Q_{k+1}$ be a splitting face. Then
    \begin{equation}\label{equation-H-lemma-S-low}
H(\xi(q))=H(q)_{low}
\end{equation}
and
    \begin{equation}\label{equation-H-lemma-S-high}
H(\gamma(q))=H(q)_{high}
\end{equation}
\end{lemma}

{\em Proof.}  The equation (\ref{equation-H-lemma-S-low}) is an immediate consequence of Lemma \ref{lemma-H-xi-sigma} below.
We shall show that the equation (\ref{equation-H-lemma-S-high}) holds. Note that as $q$ is a splitting face, $\gamma(q)$ is a non-collapsing face with $\sigma(\gamma(q))$ defined. Thus
 \[  H(\gamma(q)) = \left[
                                       \begin{array}{c}
                                         h(\gamma(q))\\
                                         0\\
                                         h(\sigma\gamma(q))\\
                                         \vdots \\
                                         h(\sigma^{(1)}(q))\\
                                       \end{array}\right]= \left[
                                       \begin{array}{c}
                                         h(q)\\
                                         0\\
                                         h(\sigma^2(q))\\
                                           \vdots \\
                                         h(\sigma^{(1)}(q))\\
                                       \end{array}\right] = \left[
                                       \begin{array}{c}
                                         h(q)\\
                                          h(\sigma(q))\\
                                         h(\sigma^2(q))\\
                                           \vdots \\
                                         h(\sigma^{(1)}(q))\\
                                       \end{array}\right]_{high} =H(q)_{high}\]
$\Box$\vskip 2mm

\subsection{Preservation of codomains}

\begin{proposition}\label{prop-H-pres-codomains}
The map $H$ defined above preserves codomains.
\end{proposition}

{\em Proof.} Let $q\in Q_{k+1}$. For $k=0$ or if $\varrho(q)\neq\ba$, the Proposition is easy.

{\em Case (A).} $q$ is a splitting face and hence $H(q)$ is a flag.

In this case, by Lemma \ref{lemma-H-collaps} (\ref{equation-H-lemma-S-high}) we have
\[ H(\gamma(q))= H(q)_{high}=\gamma(H(q)). \]

{\em Case (B).} $h(q)\in P_{k+1}$ and $\sigma(q)$ is defined.

 In this case we must have $k>1$.

Then there are $\sigma(q)\in \delta(q)$ and $\sigma^2(q)\in \delta(\sigma(q))$. Moreover, by $(B)$ and $(A)$, we have
$$H(q)=[h(q),0,H(\sigma(q))]$$
and
$$H(\sigma(q))=[h(\sigma(q)),H(\sigma^2(q))].$$
As above, there is $\bar{q}\in\delta\gamma(q)$ with $\bar{q}\leq^+\sigma^2(q)$ and $H( \sigma^2(q))=H(\bar{q})$. We have
\[ H(\gamma(q))=[h(\gamma(q)),0,H(\bar{q})]= [\gamma(h(q)),0,H(\sigma^2(q))]=\]
\[ \gamma([h(q),0,h(\sigma(q)),H(\sigma^2(q))])=\gamma([h(q),0,H(\sigma(q))])= \gamma(H(q)) \]
as required.

{\em Case (C).} $h(q)\in P_{k+1}$ and $\tau(q)$ is defined.

In this case we must have $k>1$ and we also have $\sigma\tau(q)$ defined. As $\tau(q)<^+\gamma(q)$, by Proposition \ref{prop-splitting-faces},  $\sigma(\gamma(q))\in\delta\gamma(q)$ is also defined
and we have
\[ H(\gamma(q))=[h(\gamma(q)),0,H(\sigma\gamma(q))]= [\gamma(h(q)),0,H(\sigma^2(q))]=\]
\[ \gamma([h(q),h(\tau(q)),0,H(\sigma^2(q))])=\gamma(H(q)) \]
as required.

{\em Cases (D) and (E).} Easy exercise.
$\Box$\vskip 2mm

\subsection{Preservation of domains}

\begin{lemma}\label{lemma-H-star-pres}
Let $k\geq 1$, $q\in Q_{k+1}-ker(h)$. Then, for $\bar{q}\in\delta(q)-ker(h)$, we have
\begin{equation}\label{equation-H-star-pres}
H(\bar{q})=h(\bar{q})\star H(q).
\end{equation}
\end{lemma}

{\em Proof.}
Let $k\geq 1$, $q\in Q_{k+1}-ker(h)$, $\bar{q}\in\delta(q)-ker(h)$. We shall show that (\ref{equation-H-star-pres}) holds.

\vskip 2mm
Case $\rho(q)=\ddagger\in \{ +,-\}$. Then  $\rho(q)=\rho(\bar{q})=\ddagger$. Thus
\[ H(\bar{q})=\ddagger h(\bar{q})=h(\bar{q})\star \ddagger h(q)=h(\bar{q})\star H(q). \]
\vskip 2mm

So now we assume that $\rho(q)=\ba$ for the rest of the proof. Thus, by Lemma \ref{lemma-split-seq}, we have $l\geq 0$ and a splitting sequence for $q$
\[ [q,\gamma(q),\ldots,\tau\gamma^{(l+2)}(q),\sigma\tau\gamma^{(l+2)}(q),\ldots,\sigma^{(1)}\tau\gamma^{(l+2)}(q) ] \]
with
\[ H(q)=[h(q),h(\gamma(q)),\ldots,h(\tau\gamma^{(l+2)}(q)),0,h(\sigma\tau\gamma^{(l+2)}(q)),\ldots,h(\sigma^{(1)}\tau\gamma^{(l+2)}(q)) ]. \]

\vskip 2mm {\em The idea of the proof.} We look at the maximal $\bar{l}$ and a face $t$ such that $H(q)_{\bar{l}+1}\leq^+ t\in\delta\gamma^{(\bar{l}+2)}(\bar{q})$ to compute $h(\bar{q})\star H(q)$. Then we show that this $t$ is a threshold face, as well, and hence it takes part in defining p-flag $H(\bar{q})$.  If such an $\bar{l}$ does not exist, then depending on whether $h(\sigma^{(1)}\tau\gamma^{(l+2)}(q)$ or $\gamma h(\tau\gamma^{(2)}(q)$ we get equality of flat values, either $-h(\bar{q})$ or $+h(\bar{q})$.
\vskip 2mm

We define $\bar{l}$ to be the largest $l'$ such that
\begin{enumerate}
  \item either $l'=l$ and $h(\tau\gamma^{(l+2)}(q))\leq^+ t\in \delta\gamma^{(l+2)}(h(\bar{q}))$
  \item or $l'=l-1$ and $h(\gamma\tau\gamma^{(l+2)}(q))\leq^+ t\in \delta\gamma^{(l+1)}(h(\bar{q}))$
  \item or $l'<l-1$ and $h(\gamma\sigma^{(l'+2)}\tau\gamma^{(l+2)}(q))\leq^+ t\in \delta\gamma^{(l'+2)}(h(\bar{q}))$
\end{enumerate}
In either case $t\in\delta\gamma^{(\bar{l}+2)}(h(\bar{q}))-ker(h)$ and, as $h$ is a $\iota$-map, there is a unique  $\bar{t}\in\delta\gamma^{(\bar{l}+2)}(\bar{q})$ such that $h(\bar{t})=t$.

If $\bar{l}=l$, then, as $\tau\gamma{(l+2)}(q)$ is $<^+$-minimal, we have $\tau\gamma^{(l+2)}(q)\leq^+\bar{t}$. In particular, it is a threshold face, as it is $<^+$-larger than a threshold face. So the splitting sequence for $\bar{q}$ is
\[ Sp(\bar{q})=[\bar{q},\gamma(\bar{q}),\ldots, \gamma^{(l+2)}(\bar{q}),\bar{t},\sigma\tau\gamma^{(l+2)}(\bar{q}), \ldots,\sigma^{(1)}\tau\gamma^{(l+2)}(\bar{q}) ] \]
and we have
\[  H(\bar{q}) = \left[
                                       \begin{array}{c}
                                         h(\bar{q})\\
                                         h(\gamma(\bar{q}))\\
                                         \vdots \\
                                         h(\gamma^{(l+2)}(\bar{q}))\\
                                         h(\bar{t})\\
                                         0\\
                                         h(\sigma(\bar{t}))\\
                                        \vdots \\
                                         h(\sigma^{(1)}(\bar{t})) \\
                                       \end{array}\right]= \left[
                                       \begin{array}{c}
                                         h(\bar{q})\\
                                         \gamma(h(\bar{q}))\\
                                         \vdots \\
                                         \gamma^{(l+2)}(h(\bar{q}))\\
                                         t\\
                                         0\\
                                         h(\sigma\tau\gamma^{(l+2)}(q))\\
                                        \vdots \\
                                         h(\sigma^{(1)}\tau\gamma^{(l+2)}(q)) \\
                                       \end{array}\right] =h(\bar{q})\star H(q)\]

If $0\leq \bar{l}=l-1$, then we claim
\vskip 2mm

{\em Claim 1.} $\gamma\tau\gamma^{(l+2)}(q)\leq^+ \bar{t}$.

{\em Proof of Claim 1.} From the above we have $\gamma\tau\gamma^{(l+2)}(q)\perp^+ \bar{t}$. If $h(\gamma\tau\gamma^{(l+2)}(q))<^+ h(\bar{t})$, then Claim 1 is obvious. Thus we assume that $h(\gamma\tau\gamma^{(l+2)}(q))=h(\bar{t})$. Suppose to the contrary that
\[ \bar{t}<^+ \gamma\tau\gamma^{(l+2)}(q).\]
Then we have an upper $Q-\gamma(Q)$-path of 1-collapsing faces
\[ \bar{t},a_1,\ldots, a_r,\gamma\tau\gamma^{(l+2)}(q)\]
with $r\geq 1$. By pencil linearity we have $a_1\leq^+ \gamma^{(l+1)}(\bar{q})$. By Path Lemma, since $\gamma^{(l+1)}(\bar{q})$ is a non-collapsing face and all $a_i$'s are 1-collapsing, we have that
\[ \delta\gamma^{(l+1)}(\bar{q})\ni \bar{t}<^+\gamma\tau\gamma^{(l+2)}(q)<^+\gamma^{(l)}(\bar{q})  \]
and hence
\begin{equation}\label{equation-H1}
 \tau\gamma^{(l+2)}(q) <^+ \gamma^{(l+1)}(\bar{q}).
\end{equation}

Since $\tau\gamma^{(l+2)}(q)$ is $<^+$-minimal in $Q[q]$, the inequality (\ref{equation-H1}) implies that
$$\tau\gamma^{(l+2)}(q) \leq^+ \delta\gamma^{(l+2)}(\bar{q})$$
i.e $l=\bar{l}$ contrary to the assumption. This ends the proof of Claim 1.
\vskip 2mm

In particular, $\bar{t}$ is a threshold face, as it is $<^+$-larger than a threshold face. So the splitting sequence for $\bar{q}$ is
\[ Sp(\bar{q})=[\bar{q},\ldots, \gamma^{(l+1)}(\bar{q}),\bar{t},\sigma\tau\gamma^{(l+1)}(\bar{q}), \ldots,\sigma^{(1)}\tau\gamma^{(l+1)}(\bar{q}) ]. \]
Thus
\[  H(\bar{q}) = \left[
                                       \begin{array}{c}
                                         h(\bar{q})\\
                                         \vdots \\
                                         h(\gamma^{(l+1)}(\bar{q}))\\
                                         h(\bar{t})\\
                                         0\\
                                         h(\sigma(\bar{t}))\\
                                        \vdots \\
                                         h(\sigma^{(1)}(\bar{t})) \\
                                       \end{array}\right]= \left[
                                       \begin{array}{c}
                                         h(\bar{q})\\
                                         \vdots \\
                                         \gamma^{(l+1)}(h(\bar{q}))\\
                                         t\\
                                         0\\
                                         h(\sigma^{(l)}\tau\gamma^{(l+2)}(q))\\
                                        \vdots \\
                                         h(\sigma^{(1)}\tau\gamma^{(l+2)}(q)) \\
                                       \end{array}\right] =h(\bar{q})\star H(q)\]

If $0\leq \bar{l}<l-1$ then we claim
\vskip 2mm

{\em Claim 2.} $\gamma\sigma^{(\bar{l}+2)}\tau\gamma^{(l+2)}(q)\leq^+ \bar{t}$.

{\em Proof of Claim 2.} This proof is an iterated version of the proof of the  previous Claim.

From the above we have $\gamma\sigma^{(\bar{l}+2)}\tau\gamma^{(l+2)}(q))\perp^+ \bar{t}$. If $h(\gamma\sigma^{(\bar{l}+2)}\tau\gamma^{(l+2)}(q))<^+ h(\bar{t})$, then Claim 2 holds. Thus we assume that $h(\gamma\sigma^{(\bar{l}+2)}\tau\gamma^{(l+2)}(q))=h(\bar{t})$. Suppose to the contrary that
\[ \bar{t}<^+ \gamma\sigma^{(\bar{l}+2)}\tau\gamma^{(l+2)}(q).\]
Then we have an upper $Q-\gamma(Q)$-path of 1-collapsing faces
\[ \bar{t},a_1,\ldots, a_r,\gamma\sigma^{(\bar{l}+2)}\tau\gamma^{(l+2)}(q)\]
with $r\geq 1$. By pencil linearity we have $a_1\leq \gamma^{(\bar{l}+2)}(\bar{q})$. By Path Lemma, since $\gamma^{(\bar{l}+2)}(\bar{q})$ is a non-collapsing face and all $a_i$'s are 1-collapsing, we have that
\[ \delta\gamma^{(l+1)}(\bar{q})\ni \bar{t}<^+\gamma\sigma^{(\bar{l}+2)}\tau\gamma^{(l+2)}(q)<^+\gamma^{(\bar{l}+1)}(\bar{q})  \]
and hence
\begin{equation}\label{equation-H2}
 \sigma^{(\bar{l}+2)}\tau\gamma^{(l+2)}(q) <^+ \gamma^{(\bar{l}+2)}(\bar{q}).
\end{equation}

The inequality implies that
\begin{enumerate}
  \item either $\sigma^{(\bar{l}+2)}\tau\gamma^{(l+2)}(q) \leq^+ \delta\gamma^{(\bar{l}+3)}(\bar{q})$ and there is either splitting of threshold face in $\delta\gamma^{(\bar{l}+3)}(\bar{q})$ contradicting the maximality of $\bar{l}$;
  \item or else  $\delta\gamma^{(\bar{l}+3)}(\bar{q}) <^+\sigma^{(\bar{l}+2)}\tau\gamma^{(l+2)}(q)<^+\gamma^{(\bar{l}+2)}(\bar{q})$.
\end{enumerate}

In the latter case we obtain
\begin{equation}\label{equation-H3}
\sigma^{(\bar{l}+3)}\tau\gamma^{(l+2)}(q)<^+\gamma^{(\bar{l}+3)}(\bar{q})
\end{equation}
lifting the inequality (\ref{equation-H2}) one dimension up. Iterating this argument we get
\begin{equation}\label{equation-H4}
\tau\gamma^{(l+2)}(q)<^+\gamma^{(l+1)}(\bar{q})
\end{equation}
and, as $\tau\gamma^{(l+2)}(q)$ is $<^+$-minimal, we also have
\[ \tau\gamma^{(l+2)}(q)<^+\delta\gamma^{(l+2)}(\bar{q})\]
contradicting maximality of $\bar{l}$.  This ends the proof of Claim 2.
\vskip 2mm

Again $\bar{t}$ is a threshold face, as it is $<^+$-larger than a threshold face. So the splitting sequence for $\bar{q}$ is
\[ Sp(\bar{q})=[\bar{q},\ldots, \gamma^{(\bar{l}+2)}(\bar{q}),\bar{t},\sigma\tau\gamma^{(\bar{l}+2)}(\bar{q}), \ldots,\sigma^{(1)}\tau\gamma^{(\bar{l}+2)}(\bar{q}) ]. \]
Thus
\[  H(\bar{q}) = \left[
                                       \begin{array}{c}
                                         h(\bar{q})\\
                                         \vdots \\
                                         h(\gamma^{(\bar{l}+2)}(\bar{q}))\\
                                         h(\bar{t})\\
                                         0\\
                                         h(\sigma(\bar{t}))\\
                                        \vdots \\
                                         h(\sigma^{(1)}(\bar{t})) \\
                                       \end{array}\right]= \left[
                                       \begin{array}{c}
                                         h(\bar{q})\\
                                         \vdots \\
                                         \gamma^{(\bar{l}+2)}(h(\bar{q}))\\
                                         t\\
                                         0\\
                                         h(\sigma^{(l)}\tau\gamma^{(l+2)}(q))\\
                                        \vdots \\
                                         h(\sigma^{(1)}\tau\gamma^{(l+2)}(q)) \\
                                       \end{array}\right] =h(\bar{q})\star H(q)\]

If $\bar{l}$ is not defined, then we have two cases
\begin{enumerate}
  \item $l=0$ and $\gamma^{(1)}(\bar{q})\perp^-\tau\gamma^{(2)}(q)$,
  \item $l>0$ and $\gamma^{(1)}(\bar{q})\perp^-\gamma\sigma^{(2)}\tau\gamma^{(l+2)}(q)$.
\end{enumerate}

Ad 1. If  $\gamma^{(1)}(\bar{q})<^-\tau\gamma^{(2)}(q)$, then
\[ \gamma^{(0)}(h(\bar{q}))<^+\delta h(\tau\gamma^{(2)}(q))<^+\gamma h(\tau\gamma^{(2)}(q))=\gamma(H(q)_1)\]
and by definition of $\star$ on p-flags we have
\begin{equation}\label{equation-star-minus}
h(\bar{q})\star H(q)=-h(\bar{q})
\end{equation}
On the other hand, we have
\[ \varrho(\gamma^{(0)}(\bar{q}))\leq^+\varrho(\delta\tau\gamma^{(2)}(q))=\delta\varrho(\tau\gamma^{(2)}(q))=\delta(\ba)=- \]
and hence $\varrho(\bar{q})=-$ and $H(\bar{q})=-h(\bar{q})$.

If  $\tau\gamma^{(2)}(q)<^-\gamma^{(1)}(\bar{q})$, then

\[ \gamma(H(q)_1)=\gamma h(\tau\gamma^{(2)}(q))<^+\delta h(\gamma^{(1)}(\bar{q}))<^+\gamma^{(0)}(h(\bar{q}))\]
and by definition of $\star$ on p-flags we have
\begin{equation}\label{equation-star-plus}
h(\bar{q})\star H(q)=+h(\bar{q})
\end{equation}
On the other hand, we have
\[ +=\gamma(\ba)=\gamma\varrho(\tau\gamma^{(2)}(q))= \varrho(\gamma\tau\gamma^{(2)}(q))\leq^+ \varrho(\delta\gamma^{(1)}(\bar{q})) \]
and hence $\varrho(\bar{q})=+$ and $H(\bar{q})=+h(\bar{q})$.

Thus in either case
\[ h(\bar{q})\star H(q)=H(\bar{q}).\]

Ad 2. If $\gamma^{(1)}(\bar{q})<^-\gamma\sigma^{(2)}\tau\gamma^{(l+2)}(q)$, then, as it is equality of 1-faces, we also have
\[ \gamma^{(1)}(\bar{q})<^-\sigma^{(1)}\tau\gamma^{(l+2)}(q). \]
Thus we have
\[  \gamma^{(0}(\bar{q})\leq^+\delta\sigma^{(1)}\tau\gamma^{(l+2)}(q)<^+\gamma\sigma^{(1)}\tau\gamma^{(l+2)}(q). \]
Thus
\[ \gamma^{(0)}(h(\bar{q}))=h(\gamma^{(0)}(\bar{q}))\leq^+h(\delta\sigma^{(1)}\tau\gamma^{(l+2)}(q))=h(\sigma^{(1)}\tau\gamma^{(l+2)}(q))=H(q)_0 \]
i.e. $h(\bar{q})\star H(q)=-h(\bar{q})$.
On the other hand
\[ \varrho(\gamma^{(0)(\bar{q})})\leq^+ \varrho(\delta\sigma^{(1)}\tau\gamma^{(l+2)}(q))= \delta\varrho(\sigma^{(1)}\tau\gamma^{(l+2)}(q))=\delta(\ba)=-  \]
i.e. $H(\bar{q})=-h(\bar{q})$.

If $\gamma\sigma^{(2)}\tau\gamma^{(l+2)}(q)<^-\gamma^{(1)}(\bar{q})$, then
\[ \gamma\sigma^{(1)}\tau\gamma^{(l+2)}(q)\leq^+\delta\gamma^{(1)}(\bar{q})<^+\gamma^{(0)}(\bar{q})\]
and, as $\gamma^{(1)}(\bar{q})$ is non-collapsing, we have a strict inequality
\[ H(q)_0= h(\sigma^{(1)}\tau\gamma^{(l+2)}(q))=h(\gamma\sigma^{(1)}\tau\gamma^{(l+2)}(q))\leq^+h(\delta\sigma^{(1)}\tau\gamma^{(l+2)}(q))<^+h(\gamma^{(0)}(\bar{q}))= \gamma^{(0)}(h(\bar{q})). \]
Thus $h(\bar{q})\star H(q)=+h(\bar{q})$.
Moreover
\[ +=\gamma(\ba)=\gamma\varrho(\sigma^{(1)}\tau\gamma^{(l+2)}(q))=\varrho(\gamma\sigma^{(1)}\tau\gamma^{(l+2)}(q))\leq^+\varrho(\delta\gamma^{(1)}(\bar{q})). \]
Hence $\varrho(\bar{q})=+$ and $H(\bar{q})=+h(q)$. Thus in either case
\[ h(\bar{q})\star H(q)=H(\bar{q}).\]
$\Box$\vskip 2mm

\begin{proposition}\label{prop-H-pres-domains}
The map $H$ defined above preserves domains.
\end{proposition}

 {\em Proof.}  For $q\in Q_1$, preservation of domains is easy. We assume $q\in Q_{k+1}$ with $k>0$ and we shall check that in all cases (A), (B), (C), (D), and (E), $H$ preserves the domain of $q$ i.e. we have
 a bijection
 \begin{equation} \label{equation-H-pres-dom}
 \delta(q)-ker(H)\lra  \delta(H(q))
 \end{equation}
 induced by $H$.

\vskip 2mm
{\em Case (A).} $q$ is a splitting face.

Since $h(q)\in P_k$, by definition of $\iota$-maps, we have a bijection
\[ \{\xi(q)\} =\delta(q)-ker(h)\lra \delta^{(k)}(h(q))=\{ h(q)\} \]
induced by $h$. Moreover, $H(\sigma(q))$ is a flag, and hence it is of dimension $k$. The other face $q'\in \delta(q)-\{ \xi(q),\sigma(q)\}$ (if exists) is not a splitting face and $h(q')\in P_{<k}$. So such $q'\in ker(H)$. Thus $$\delta(q)-ker(H)=\{ \sigma(q), \xi(q)\}.$$

On the other hand, $H(q)$ is a flag and
 $$\delta(H(q))=\{ H(q)_{side}, H(q)_{low}\}.$$

 By case (A) of the definition of $H$, we have that $H(\sigma(q))=H(q)_{side}$. Thus to prove that $H$ induces a bijection (\ref{equation-H-pres-dom}), we need to verify that
\begin{equation}\label{equation-H-q-xi-low}
H(\xi(q))=H(q)_{low}.
\end{equation}
But this is the content of Lemma \ref{lemma-H-collaps}(\ref{equation-H-lemma-S-low}). This proves case (A) of preservation of domains.

\vskip 2mm

{\em Case (B).} $h(q)\in P_{k+1}$ and $\sigma(q)$ is defined.

In this case $H(q)$ is a high p-flag and
\[ \delta(H(q)) = \{ h(\bar{q})\star H(q) \}_{\bar{q}\in \delta(q)-ker(h)}\cup \{ H(q)_{(k+1)} \}. \]

From the definition of $H$ we have that if a face $q'\in ker(h)$, then either $q'$ is a splitting face or $q'\in ker(H)$. Thus
\[ \delta(q)-ker(H) = ( \delta(q)-ker(h))\cup \{ \sigma(q) \}. \]
As $H(q)$ is a high p-flag, from the definitions $\delta$ on $\Cyl(P)$ and case (A) of $H$ we get that $H(\sigma(q))=H(q)_{(k+1)}$. From Lemma \ref{lemma-H-star-pres} we get that, for $\bar{q}\in\delta(q)-ker(h))$, we have
\[ H(\bar{q})= \bar{q})\star H(q). \]
Thus in this case $H$ induces bijection (\ref{equation-H-pres-dom}), as required.

\vskip 2mm

{\em Case (C) and (D).} $h(q)\in P_{k+1}$ and $\sigma(q)$ is not defined.

In these cases there are no faces in $\delta(q)$ that are sent to flag faces in $\Cyl(P)$. Hence
 $$\delta(q)-ker(H)=\delta(q)-ker(h)$$
and
\[ \delta(H(q))= \{ h(\bar{q})\star H(q) \}_{\bar{q}\in \delta(q)-ker(h)}.\]
Thus again using Lemma \ref{lemma-H-star-pres} we get that $H$ induces bijection (\ref{equation-H-pres-dom}) in this case, as well.

\vskip 2mm

{\em Case (E).} $q\in ker(h)$ and $q$ is not splitting.

The condition is non-trivial if $h(q)\in P_{k}$. Then we have that
\[ \delta(q)-ker(H)=\delta(q)-ker(h)=\{ \xi(q)\} \]
and
\[ \delta(H(q))=\{ H(\gamma(q))\}.\]
Thus the bijection (\ref{equation-H-pres-dom}) is the consequence of Lemma \ref{lemma-H-collaps-ga} in this case.
$\Box$\vskip 2mm

\subsection{$\pOpei$ is a test category}

Thus we have proved that `opetopes think' that $\Cylh$ is a product of $I$ and $P$ in $\pHgi$.
Recall that $\Cyli(P)$ is $\kappa_!(\Cylp(P))$, i.e. left Kan extension of $\Cyli(P)$ along $\kappa : \pOpe\ra \pOpei$.

\begin{corollary}\label{coro-product-IxP-pOpei}
$\Cyli(P)$ is the product of $I$ and $P$ in $\widehat{\pOpei}$.
\end{corollary}

{\em Proof.}
$\cHi$ is full embedding on hom-sets $\pHgi(Q,H)$ where $Q$ is a positive opetope and $H$ is a opetopic hypergraph. Therefore
the universal property of  $\Cylhi(P)$ with respect to positive opetopes, Proposition \ref{prop-univ-prop-cylh-pHgi} translates to universal property of the $\Cyli(P)$ with respect to resentables in $\widehat{\pOpei}$. But to have universal property of a product in $\widehat{\pOpei}$ with respect to representable functors is the same as to be a product. Thus $\Cyli(P)$ is the product of $I$ and $P$ in $\widehat{\pOpei}$.   $\Box$\vskip 2mm

\begin{theorem}\label{thm-pOpei-test-cat}
$\pOpei$ is a test category.
\end{theorem}

{\em Proof.} $\pOpei$ has terminal object and hence, by Proposition \ref{prop-test-characterization}, we need to verify that $\pOpei$ is a local test category. By Proposition \ref{prop-local-test-characterization}, it is enough to show that, for any positive opetope $P$, the product $I\times P$ in $\widehat{\pOpei}$ satisfies the assumptions of Proposition \ref{prop-aspherical-presh-characterization}. By Corollary \ref{coro-product-IxP-pOpei} it is enough to show that $\Cyli(P)$ is straight. Since the left Kan extension along $\kappa : \pOpe \ra \pOpei$ preserves colimits, representables and monomorphisms among them, it preserves straight objects. Thus it is enough to verify that $\Cylp(P)$ is a straight object. This was done in the proof of Theorem \ref{thm-cyl-aspherical-pOpe}.
$\Box$

\vskip 2mm

\begin{thebibliography}{CWMW}
\bibitem[ACM]{ACM}
D. Ara, D.-C. Cisinski and I. Moerdijk, {\em The Dendroidal Category is a  Test Category}, arXiv:1703.07098, pp. 1-14.

\bibitem[BD]{BD}
 \frenchspacing
J. Baez, J. Dolan,  {\em Higher-dimensional algebra III: n-Categories and the algebra of opetopes}. Advances in Math. 135 (1998), pp. 145-206.

\bibitem[Be]{Berger} C. Berger, {\em A Cellular Nerve for Higher  Categories}. Adv. in Mathematics 169,  (2002), pp. 118-175.

\bibitem[Bu]{Bu} A. Burroni, {\em Higher-dimensional word problems with applications to equational logic}, Theoretical
Computer Science 115 (1993), no. 1, pp. 43–62.

\bibitem[Cis]{Cis} D-C. Cisinski, {\em Les Pr\'{e}faisceaux Comme Mod\`{e}les des Types d'Homotopie}, Asterisque, vol 308, (2006), pp. 1-392.

\bibitem[Gr]{Gr} A. Grothendieck,  {\em Pursuing stacks}, Manuscript, 1983, to be published in Documents Math\'{e}matiques.

\bibitem[HMP]{HMP} C. Hermida, M. Makkai, J. Power, {\em On weak higher dimensional categories, I} Parts 1,2,3, J. Pure and Applied Alg.  153 (2000), pp. 221-246, 157 (2001), pp. 247-277, 166 (2002), pp. 83-104.

\bibitem[Mlt]{Mlt} G. Maltsiniotis, {\em La Theorie de l'Homotopie de Grothendieck}, Asterisque, vol 301, (2003).

\bibitem[P]{Palm} T. Palm, {\em Dendrotopic sets}, Hopf algebras, and semiabelian categories, Fields Inst. Commun. vol. 43 (2004), 411-461 AMS, Providence, RI.

\bibitem[SZ]{SZ} S. Szawiel, M. Zawadowski, {\em The web monoid and opetopic sets}, J. of Pure and Applied Algebra 217 (2013), pp. 1105–1140.

\bibitem[Th]{Th} R. Thomason, {\em $\Cat$ as a closed model structure}, Cahiers de Topo. et G\'eom Diff\'{e}rentielle Cat\'{e}gorique, tome 21, no 3, (1980), pp. 305-324.

\bibitem[Z1]{Z1} M. Zawadowski, {\em On positive face structures and positive-to-one computads}. Preprint, 2006, pp. 1-77.

\bibitem[Z2]{Z2} M. Zawadowski, {\em On ordered face structures and many-to-one computads}. Preprint, (2007), pp. 1-95.

\bibitem[Z3]{Z3} M. Zawadowski, {\em Lax Monoidal Fibrations}, in Models, Logics, and Higher-Dimensional Categories: A Tribute to the Work of Mihály Makkai
(B. Hart, et al. , editors) (CRM Proceedings 53, 2011), pp. 341-424.
\end{thebibliography}
\end{document}